\documentclass[12pt]{article}
\usepackage{authblk}
\usepackage{amssymb}
\usepackage{amsfonts}
\usepackage{amsmath}
\usepackage{mathrsfs}
\usepackage{graphicx}
\usepackage{subfigure}
\usepackage{color}

\newtheorem{theo}{Theorem}[section]
\newtheorem{prop}{Proposition}[section]
\newtheorem{rmk}{Remark}[section]
\newtheorem{lem}{Lemma}[section]
\newtheorem{coro}{Corollary}[section]

\numberwithin{equation}{section}

\def\cH{\mathcal H}

\def\RR{\mathbb R}
\def\R{\mathbb R}
\def\cL{\mathcal L}

\def\cV{{\mathcal V}}

\def\ds{\displaystyle}

\def\ud{\, \mathrm{d}}

\def\div{\mathrm{div}}

\def\blackbox{\leavevmode\vrule height 5pt width 4pt depth 0pt\relax}
\def\endproof{\null\hfill {$\blackbox$}\bigskip}

\long\def\symbolfootnote[#1]#2{\begingroup%
\def\thefootnote{\fnsymbol{footnote}}\footnote[#1]{#2}\endgroup}

\begin{document}

\title{Efficient numerical calculation of drift and diffusion coefficients in the diffusion approximation of kinetic equations}

\author[1]{V.~Bonnaillie-No\"el\thanks{{\tt bonnaillie@math.cnrs.fr}}}
\author[2]{J.A.~Carrillo\thanks{{\tt carrillo@imperial.ac.uk}}} 
\author[3]{T.~Goudon\thanks{{\tt thierry.goudon@inria.fr}}}
\author[2]{G.A.~Pavliotis\thanks{{\tt g.pavliotis@imperial.ac.uk}}}

\affil[1]{\small D\'epartement de Math\'ematiques et Applications UMR 8553, PSL, CNRS, ENS Paris

45 rue d'Ulm
F-75230 Paris cedex 05, France 
}

\affil[2]{Department of Mathematics, Imperial College London

180 Queen's Gate
London
SW7 2AZ, UK
}

\affil[3]{Inria,  Sophia Antipolis M\'editerran\'ee Research Centre, Project COFFEE

\& Univ. Nice Sophia Antipolis, CNRS, Labo. J. A. Dieudonn\'e, UMR 7351 

Parc Valrose, F-06108 Nice, France}
\maketitle

\begin{abstract}

In this paper we study the diffusion approximation of a swarming model given by a system of interacting Langevin equations with nonlinear friction. The diffusion approximation requires the calculation of the drift and diffusion coefficients that are given as averages of solutions to appropriate Poisson equations. We present a new numerical method for computing these coefficients that is based on the calculation of the eigenvalues and eigenfunctions of a Schr\"{o}dinger operator. These theoretical results are supported by numerical simulations showcasing the efficiency of the method. 

\end{abstract}

\section{Introduction}

In this paper we consider models from statistical physics that have the general form
\begin{equation}\label{1}
\partial_t f^\varepsilon+\ds\frac{1}{\sqrt{\varepsilon}}
(v\cdot\nabla_r f^\varepsilon-\nabla_r\Phi\cdot\nabla_vf^\varepsilon)=\ds\frac1\varepsilon Q(f^\varepsilon).
\end{equation}
Here, $f^\varepsilon(t,r,v)$ denotes the distribution in phase space of a certain population (of particles, individuals...):
$r\in\mathbb R^d$, $v\in\mathbb R^d$ stand for the space and velocity variables,  respectively.
Equation \eqref{1} is written in dimensionless form and we  consider the regime where the parameter $\varepsilon>0$, which depends on the 
typical length and time scales of the phenomena under consideration, is small.
The potential $(t,r)\mapsto \Phi(t,r)$ can be defined self--consistently, typically by a convolution with the macroscopic density $\int f_\varepsilon\ud v$.
On the right hand side,  $Q$ is a  linear operator, which can be either integral or differential operator with respect to the variable $v$; it is intended to describe some 
collision or friction--dissipation mechanisms.
Of course, the asymptotics is driven by the properties of the leading  operator $Q$:
\begin{itemize}
\item[\small$\bullet$]
We assume that $Q$ is conservative in the sense that
$$
\int Q(f^\varepsilon)\ud v=0
$$
holds.
Consequently, integrating \eqref{1} with respect to the velocity variable, we obtain the
following local  conservation law
$$
\partial_t \rho^\varepsilon+\nabla_r\cdot J^\varepsilon=0
$$
where we denote
$$
\rho^\varepsilon(t,r)=\ds\int f^\varepsilon(t,r,v) \ud v,\qquad J^\varepsilon=\ds\int \ds\frac{v}{\sqrt{\varepsilon}} f^\varepsilon(t,r,v) \ud v.
$$
\item[\small$\bullet$] We also suppose that the kernel of $Q$ is spanned by a positive normalized function $v\mapsto F(v)$.
As we shall detail below on specific examples, $Q$ also satisfies some dissipation properties: roughly speaking these properties encode the fact that $Q$ forces the distribution function to become proportional to the equilibrium $F$.
In turn, the dissipation  permits us to establish estimates
which lead to the ansatz
$f^\varepsilon(t,r,v)=\rho^\varepsilon (t,r)F(v)+\sqrt\varepsilon g^\varepsilon(t,r,v)$.
In order to have a current of order 1, we should assume 
$\int vF\ud v=0$, so that $J^\varepsilon=\int g^\varepsilon\ud v$.
\end{itemize}
Finally, let us suppose that $\rho^\varepsilon, J^\varepsilon$ and $g^\varepsilon$ have well defined limits as $\varepsilon \rightarrow 0$, denoted $\rho, J$, and $g$ respectively.
Multiplying~\eqref{1} by $\sqrt\varepsilon$ and  letting $\varepsilon$ go to 0 
 yields
$$
Q(g)=(v\cdot\nabla_r \rho)F - (\nabla_r \Phi\cdot \nabla_v F) \rho ,
$$
which provides us with a closure for  the  equation obtained by passing to the limit in the conservation law
$$
\partial _t \rho +\nabla_r\cdot\left(\int  v g\ud v\right)=0.
$$
Owing to the linearity of $Q$, we rewrite this as a convection--diffusion equation:
\begin{equation}\label{dd1}
\partial_t \rho-\nabla_r(\mathcal D\nabla_r \rho+\mathcal  K\nabla_r \Phi \rho)=0,
\end{equation}
where the coefficients are defined by the following matrices
$$\mathcal D=-\int v\otimes Q^{-1}(vF)\ud v,\qquad
\mathcal  K =  \int v\otimes Q^{-1}(\nabla_v F)\ud v,$$
with $Q^{-1}$ the pseudo-inverse of $Q$, which has to be properly defined on the orthogonal of $\mathrm{Span}\{F\}$.
In particular, the dissipative nature of $Q$ implies that $\mathcal D$ is nonnegative. 

There is a huge literature on the analysis of such asymptotic problems, motivated by various application fields (radiative transfer theory, neutron transport, the modelling of semiconductors, population dynamics, etc): we refer the reader for instance  to \cite{BookBGP,DGP,EGM,G2,GOME,MT,PS} for an overview of results and mathematical techniques used to handle this question, and for further references. The asymptotic analysis of~\eqref{1} is of great practical interest: it is clear that the numerical simulation of an equation like \eqref{dd1} is by far less costly than the one of \eqref{1}: on the one hand, we have eliminated the velocity variable, on the other hand, when $\varepsilon \to 0$, \eqref{1} contains stiff terms that induce prohibitive stability conditions. However, we are left with the difficulty of calculating the effective diffusion and drift coefficients $\mathcal D$ and $\mathcal K$ which relies on being able to calculate the inverse of the operator $Q$, i.e. on solving equations of the form $Q\phi=h$.

We address these questions in the specific case of kinetic models for swarming for which the auxiliary equations that determine $\mathcal D$ and $\mathcal K$ do not have explicit solutions.
In recent years, increasing efforts have been devoted to the development of mathematical models able to describe the self--organization of a large set of living organisms (fish, birds, bacteria...),   after the pioneering work of Vicsek {\it et al.} \cite{Vic}. Based on simple rules of information exchange between the individuals about their close environment, the whole population organizes itself in remarkable patterns.
Starting from individual-based description, kinetic equations can be derived by means of mean--field regimes, as it is usual in statistical physics \cite{Dob,Gol2,Sp,Sz}: we refer the reader to \cite{BCC11,CCR,CFRT10,HT,Laf} for a thorough study of such models. For the models we are interested in, the interaction operator $Q$ has the form 
$
Q(f)=\mathrm{div}_v(\nabla_v W(v)f+\nabla_v f)\,,
$
that involves a quite complicated potential function $v\mapsto W(v)$.
The equilibrium function simply reads $F(v)=\frac{1}{Z}e^{-W(v)}$, where $Z$ denotes the normalization constant, but inverting $Q$ is not that simple.
In contrast to the standard Fokker--Planck operator corresponding to the Langevin dynamics with linear friction, i.e. when $W(v)=v^2/2$,  in general $v\mapsto vF(v)$
is not an eigenfunction of $Q$, and we do not have explicit formulas for the effective coefficients.
Our approach for computing the coefficients is based on the spectral properties of the operator. 
In fact, a unitary transformation enables us to reformulate the Poisson equation  $Q\phi=h$ in the form $\mathcal H u = \mu$, where $\mathcal H$ is a Schr\"odinger operator associated to a certain potential $\Phi$, that depends on $W$; see Equation~\eqref{e:potential-schrod} below. In turn, the operator admits a spectral decomposition.
Then, our method is based on expanding the functions on the eigenbasis, and then computing the coefficients of $\phi$ in this basis.
The latter step does not present any  difficulty and the computational effort is concentrated on the determination of the eigenelements. In practice, we work on the discrete form of the equations, and we expect that only a few  
eigenmodes are necessary to capture the effective coefficients.
This is the strategy we shall discuss, adopting high--order Finite Elements discretizations of the Schr\"odinger equations, 
which allows us to make use of performing computational tools \cite{Melina}.

We mention now some related problems, to which the numerical method developed in this paper can also be applied, and alternative approaches for the calculation of the coefficients that appear in macroscopic equations. First, it should be noted that a similar formalism, based on the calculation of the eigenfunctions and eigenvalues of the linearized collision operator, has been developed for the calculation of transport coefficients using  the linearized Boltzmann equation~\cite{ResibDeLeen77}. 
Second, the method can be adapted to compute the effective, possibly space--dependent, coefficients that come from the homogenization process of advection-diffusion equations with periodic coefficients, see e.~g.~\cite{BLP,GP}, the diffusion approximation for fast/slow systems of stochastic differential equations~\cite{PV01, PV03, PV05}, and in connection to functional central limit theorems for additive functionals of Markov processes~\cite{OllaKomorowskiLandim12}. In all these problems, the drift and diffusion coefficients that appear in the macroscopic equation can be expressed in terms of the solution of an appropriate Poisson equation. These effective coefficients can be alternatively calculated using either Monte Carlo simulations, e.g.~\cite{PavlStBan06, Boozer1981}, the heterogeneous multiscale method~\cite{ELV05}, numerical solution of the Poisson equation using spectral methods, e.g.~\cite{MajMcL93, thesis}. Other numerical approaches include the use of linear response theory and of the Green-Kubo theory~\cite{JS12} and the expansion of the solution to the Poisson equation in appropriate orthonormal basis functions, e.g. Hermite polynomials. This technique, which is related to the continued fraction expansion~\cite{Ris84} has been applied to the calculation of the diffusion coefficient for the Langevin dynamics in a periodic potential~\cite{PavlVog08} and to simple models for Brownian motors~\cite{LatorreKramerPavliotis2014}. 

The rest of the paper is organized as follows.
In Section~\ref{Section2}, we introduce the kinetic model for swarming we are interested in.
By using the Hilbert expansion, we detail the diffusion asymptotics, and we discuss some properties of the effective coefficients. 
The convergence can be rigorously justified, and a complete proof is given in  the Appendix.
In Section~\ref{Section3}, we switch to the 
formalism of Schr\"odinger equations and we explain in detail our numerical strategy.
We also provide estimates to justify the approximation by truncating the Fourier series.
Section~\ref{Section4} is devoted to the numerical illustration of the approach, dealing with a relevant set of potentials $W$. 
In particular we bring out the role of the scaling coefficients that appear in the potential,
with difficulties related to the so--called ``tunnelling effect'' \cite{Rob87,HeSj84}.
 Section~\ref{SectionC} is reserved for conclusions. The proof of Theorem~\ref{CVTh} on the mean field limit approximation can be found in the appendix.

%%%%%%%%%%%%%%%%%%%%%%%%%%%%%%%%%%%%%%%%%%%%%%%%%5

\section{Motivation: Swarming Models\label{Section2}}

The analysis of interacting particle systems with random noise
finds applications in collective behavior and self-organization of
large individuals ensembles. These models lead typically to systems 
with non-standard friction terms. In particular, the following example 
was proposed in the literature in \cite{DorChuBerCha06}, which  includes the
effect of self-propulsion and a Rayleigh type friction to impose
an asymptotic cruising speed for individuals. This system with
noise and linear Stokes friction was considered in
\cite{CarDorPan09}. For $\mathcal N$ particles, it reads
\begin{equation}
\label{e:swarming}
\left\lbrace
\begin{array}{l}
\displaystyle \ud {r}_i = v_i \ud t, \vspace{.3cm}\\
\displaystyle \ud{v}_i 
= \left[\displaystyle (\alpha - \beta\,|{v}_i|^2) {v}_i - \frac 1 {\mathcal N} \nabla_{r_i} \sum_{j \neq i} U(|r_i - r_j|)\right]\ud t
+\sqrt{2\sigma}\ud\Gamma_i(t),
\end{array}
\right.
\end{equation}
where $\Gamma_i(t)$ are $\mathcal N$ independent Brownian motions with values in \(\RR^d\) and $\sigma >0$ is the noise
strength. Here, $\alpha>0$ is the self-propulsion strength
generated by the organisms, $\beta$ is the friction coefficient, and $\gamma=\alpha/\beta$ is the squared asymptotic cruise speed of the individuals.
Notice that $(\alpha - \beta \,|v|^2) v=-\alpha\nabla_v W(v)$ with
$W(v)=\tfrac{1}{4\gamma} |v|^4 - \tfrac{1}2 |v|^2$. Different models for friction can lead to asymptotic fixed speeds for individuals, see for instance \cite{LevRapCoh00}. It is interesting, therefore, to study properties of more general potentials in velocity $W(v)$, with the basic behavior of having a confinement when the speed is large and that zero speed is a source fixed point of the speed dynamics. The class of potentials for the velocity that we will consider in this paper includes, for instance, potentials of the form $W(v)=\tfrac{1}{a\gamma} |v|^a - \tfrac{1}b |v|^b$, with $a>b\geq 1$.

The mean-field limit of the stochastic particle system~\eqref{e:swarming} above under
suitable assumptions on the interaction potential $U$ is given by
the following kinetic Fokker-Planck equation, see \cite{BCC11}:
\begin{equation*}
\frac{\partial f}{\partial t} + v \cdot \nabla_r f - \div_v
\left[(\nabla_r U \star \rho)  f \right]=\div_v [\alpha \nabla_v
W(v) f+\sigma\nabla_v \,f ],
%\label{e:fp}
\end{equation*}
where 
\begin{equation*}\label{e:density}
\rho(t,r)=\ds\int f(t,r,v)\ud v,
\end{equation*}
and $\star$ stands for the usual convolution product with respect to the space variable.

Let us first remark that ${\mathcal T}_{F}=\frac1{\alpha}$ is the
natural relaxation time for particles to travel at asymptotic
speed $\sqrt{\alpha/\beta}$.  We introduce the time and length scales $\mathrm T$ and $\mathrm L$, which
are determined by the time/length scales of observation. They define the  speed  unit $\mathrm U=\mathrm L/\mathrm T$ that will be compared to $\mathrm V$, the typical particle speed
and ${\mathcal V}_{th}=\sqrt{\sigma/\alpha}$, the typical value of fluctuations
in particle velocity, called the thermal speed. Then we can define dimensionless
variables, denoted by primed quantities, as
$$
t=\mathrm T \, t',\qquad  r=\mathrm L \,  r',\qquad v=
%{\mathcal V}_{th} 
\mathrm V \, v',
$$
$$
f'(t',r',v')= \mathrm L^d \ %{\mathcal V}_{th}
{\mathrm V}
^d\ f(\mathrm T
t',\mathrm L r',%{\mathcal V}_{th}
{\mathrm V}
 v'),\qquad \mbox{and} \qquad
U'(r')= C^2\ U(\mathrm L r').
$$
Note that $C$ has the dimension of velocity. With this rescaling, we obtain the following
dimensionless kinetic equation
\begin{align}
\frac{\partial f}{\partial t} + \eta_1 \, v \cdot \nabla_r f -
\eta_2 \div_v \left[(\nabla_r U \star \rho)  f \right]=
\frac1\varepsilon Q (f), \label{mainfp3}
\end{align}
where primes have been eliminated for notational simplicity and where the  operator $Q$ is defined as 
\begin{align}
Q (f) \equiv \div_v [
\nabla_v W(v) f+\theta \nabla_v \,f ],\label{e:fp2}
\end{align}
with $W(v)=\tfrac{1}{4\gamma} |v|^4 - \tfrac{1}2 |v|^2$ for the problem corresponding to \eqref{e:swarming}. Here, $\eta_1$, $\eta_2$, $\gamma$, $\theta$ and $\varepsilon$ are dimensionless parameters given by
$$
\eta_1 =\frac{%{\mathcal V}_{th}}
{\mathrm V}}{\mathrm{U}}, \qquad \eta_2 =
\frac{C^2}{\mathrm{U}{\mathrm V}%\mathcal V}_{th}
}, \qquad
\gamma=\frac{\alpha}{\beta {\mathrm V}%\mathcal V}_{th}
^2}, 
\qquad
\theta=\left(\ds\frac{{\mathcal V}_{th}}{\mathrm V}\right)^2,
\qquad \mbox{and}
\qquad \varepsilon=\frac{{\mathcal T}_{F}}{T}=\frac1{T\alpha} .
$$
The operator $Q$ defined in~\eqref{e:fp2} is the Fokker-Planck operator corresponding to the stochastic differential equation
\begin{equation}\label{e:sde-vel}
\ud v(t) = - \nabla_v W(v(t)) \ud t + \sqrt{2 \theta} \ud \Gamma (t),
\end{equation}
where $\Gamma(t)$ denotes standard $d$-dimensional Brownian motion. The generator of this diffusion process is 
\begin{equation}\label{eq.defL}
\cL = - \nabla_v W(v) \cdot \nabla_v + \theta \Delta_v.
\end{equation}

Now, assume the following relation between the dimensionless
parameters, with a  finite asymptotic dimensionless speed,  
$$
\eta_1\simeq \eta_2 \simeq \varepsilon^{-1/2},\qquad 
\gamma,\theta\simeq O(1).
$$ 
With these assumptions, Equation~\eqref{mainfp3} becomes
\begin{equation}\label{eq.nnn}
\frac{\partial f}{\partial t} + \frac 1{\sqrt\varepsilon}\Big(v \cdot \nabla_r f - \div_v \left[(\nabla_r U \star \rho)  f \right]\Big)
= \frac1\varepsilon Q (f).
\end{equation}
In this regime, the dominant mechanisms are the noise and the nonlinear friction. This scaling is the so-called diffusion scaling for kinetic equations, see \cite{BookBGP,DGP} and the references therein. We remark that for the particular application of swarming, other distinguished limits can also be considered, see \cite{CarDorPan09} for details.
Different interesting features of the model arise, depending on the relative magnitudes of the scaling parameters $\gamma$, $\theta$.  

In order to obtain a closed macroscopic equation for the density $\rho$ in the limit $\varepsilon\to 0$, we use the standard Hilbert expansion method. Inserting the following Hilbert expansion
\begin{equation*}%\label{hilb} 
f^\varepsilon=f^{(0)}+\sqrt\varepsilon f^{(1)}+\varepsilon
f^{(2)}+\dots \qquad \mbox{and} \qquad
\rho^\varepsilon=\rho+\sqrt\varepsilon \rho^{(1)}+\dots
\end{equation*}
into \eqref{mainfp3} and identifying terms with equal power of
$\sqrt\varepsilon$, we get:
\begin{itemize}
\item $\varepsilon^{-1}$ terms: $Q (f^{(0)})=0$ which implies
that $f^{(0)}(t,r,v)=\rho(t,r)\ M(v)$, where $$M(v)=Z^{-1} e^{-W(v)/\theta},\qquad
Z=\ds\int e^{-W(v)/\theta} \ud v, $$ is
the Maxwellian distribution associated to the
Fokker-Planck operator $Q $ and $Z$ denotes the normalization constant. This is particularly clear by rewriting the Fokker-Planck operator $Q$ as 
\begin{equation}\label{formfp}
Q (f) =\theta\nabla_v \cdot\left[ M
\nabla_v\left(\frac{f}{M}\right)\right].
\end{equation}

\item $\varepsilon^{-1/2}$ terms:
\begin{align*}
Q  (f^{(1)}) &= v\cdot \nabla_r f^{(0)} - \div_v
\left[(\nabla_r U \star \rho)  f^{(0)} \right]\\
&=\left[ v\cdot \nabla_r \rho+\ds\frac1\theta\nabla_v W(v)\cdot(\nabla_r U \star
\rho)\rho \right]\,M(v) .
\end{align*}
To invert this equation, we need to solve the following problems:
\begin{subequations}\label{e:invert}
\begin{equation}\label{invert-1}
Q(\chi_i)=v_i M(v)  \; , %\sqrt{M(v)} 
\end{equation}
\begin{equation}\label{invert-2}
Q(\kappa_i)=\ds\frac1\theta\frac{\partial W}{\partial v_i}(v) M(v).% \sqrt{M(v)}.
\end{equation}
\end{subequations}
Note that $\int h\ud v=0$ is clearly a  necessary condition for the problem $Q(f)=h$ to admit a solution. Assuming that the potential $W$ increases sufficiently fast as $|v| \rightarrow +\infty$ and given that the righthand side of~\eqref{invert-2} is equal to $-\frac{\partial M(v)}{\partial v_{i}}$, from the divergence theorem we deduce that the solvability condition is  satisfied for~\eqref{invert-2}. The solvability condition is satisfied for~\eqref{invert-1} under, for example, the assumption that the velocity potential $W$ is spherically symmetric. In Section~\ref{Sec.CASC}, we will see how the case of nonsymmetric potentials where the compatibility condition~\eqref{invert-1} is not fulfilled can be dealt with. 

Existence of a solution for equations of the form $Q(f)=h$ relies on the possibility to apply the Fredholm alternative. For this it is sufficient to show that $Q$ has compact resolvent in the space $L^{2}(\R^{d} ; M^{-1}(v))$. This follows, for example, by assuming that the potential satisfies 
\begin{equation}\label{e:poincare}
\lim_{|v| \rightarrow +\infty} \left(\frac{1}{2} |\nabla_{v} W(v)|^{2} - \Delta_{v} W(v) \right) = +\infty.
\end{equation}
See, for example,~\cite[Thm.~A.19]{Vil04HPI}. Under this assumption we can apply the Fredholm alternative to obtain $f^{(1)}(t,r,v)=\chi \cdot \nabla_r
\rho+\kappa\cdot(\nabla_r U \star \rho)\rho$ with
$\chi=(\chi_1,\dots,\chi_d)$ and  $\kappa=(\kappa_1,\dots,\kappa_d)$.
\item  $\varepsilon^{0}$ terms:
\begin{align*}
Q (f^{(2)}) = \, &\partial_t f^{(0)}+ v\cdot \nabla_r
f^{(1)} - \div_v \left[(\nabla_r U \star \rho)  f^{(1)} \right]-
\div_v \left[(\nabla_r U \star \rho^{(1)}) f^{(0)} \right],
\end{align*}
with $\rho^{(1)} =\int f^{(1)}\, dv$. However, using again the compatibility condition, we conclude that
\begin{equation}\label{DD}
\partial_t \rho-\nabla_r \cdot\big(\mathcal{D} \nabla_r\rho+ \mathcal{K} (\nabla_r U \star
\rho)\rho\big)=0\,,
\end{equation}
where the diffusion and drift matrices are given by:
\begin{equation}\label{DD2}
\mathcal{D} = -\int_{\RR^d} v\otimes\chi \ud v \qquad \mbox{and}
\qquad \mathcal{K} = -\int_{\RR^d} v\otimes\kappa \ud v\,.
\end{equation}
\end{itemize}
Therefore, in the $\varepsilon\to 0$ limit regime we expect the
macroscopic density to be well approximated by the solution to the
aggregation-diffusion equation \eqref{DD}. We remark that by setting $\widehat{\chi} = \chi M$ and $\widehat{\kappa} = \kappa M$, the Poisson equations~\eqref{e:invert} become
\begin{equation}\label{e:invert-L}
\cL (\widehat{\chi}_i)=v_i   \; , \qquad\mbox{and}\qquad
\cL (\widehat\kappa_i) = \ds\frac1\theta\frac{\partial W}{\partial v_i}(v)\; ,
\end{equation}
where $\cL$, defined by  \eqref{eq.defL}, denotes the generator of the diffusion process $t\mapsto v(t)$  defined in~\eqref{e:sde-vel}.
The drift and diffusion coefficients~\eqref{DD2} are also given by the formulas
\begin{equation*}%\label{DD2L}
\mathcal{D} = -\int_{\RR^d} v\otimes\widehat{\chi} \,M \ud v  \qquad \mbox{and}
\qquad \mathcal{K} = -\int_{\RR^d} v \otimes \widehat{\kappa}\, M \ud v\, .
\end{equation*}
This is the form of the Poisson equation that appears in functional central limit theorems and in the diffusion approximation for fast/slow systems of stochastic differential equations~\cite{EthKur86,OllaKomorowskiLandim12}.

Furthermore, we note that the diffusion matrix is positive definite
and thus, we can talk properly about a diffusion matrix. To show this, we use \eqref{formfp} to deduce
\begin{equation}\label{coer}
-\int_{\RR^d} Q(f) \frac{g}{M}\ud v =\theta \int_{\RR^d} M
\nabla_v\left(\frac{f}{M(v)}\right) \cdot
\nabla_v\left(\frac{g}{M(v)}\right) \ud v\,.
\end{equation}
Now, given any vector $\xi\in\RR^d\setminus\{0\}$, we can compute
\begin{align*}
 \mathcal D \xi \cdot\xi&= -\int_{\RR^d} \left(\left[v M(v)\right]\cdot\xi
\right)\left( \chi\cdot\xi \right) \frac{1}{M(v)}\ud v \\
&=-\int_{\RR^d} Q(\chi\cdot\xi)\frac{\chi\cdot\xi}{M(v)}\ud v =
\theta\int_{\RR^d} M
\left|\nabla_v\left(\frac{\chi\cdot\xi}{M(v)}\right)\right|^2
\ud v> 0\,.
\end{align*}
The strict equality follows from the fact that $\chi \cdot \xi/M \neq \mbox{const}$, as we can immediately deduce from~\eqref{invert-1}. 

The analysis of the asymptotic regime remains technically close to the derivation of the diffusion regimes for the Vlasov-Poisson-Fokker-Planck equation \cite{EGM, G2, PS, MT}. We detail in the Appendix the proof of the following statement.

\begin{theo}\label{CVTh}
Let us consider a sequence of  initial data $f^\varepsilon_{\mathrm{Init}}\geq 0$ that satisfies
\[
\ds\sup_{\varepsilon>0}\ds\iint \big(1+|r|+ W+|\ln (f^\varepsilon_{\mathrm{Init}}|)\big)f^\varepsilon_{\mathrm{Init}}\ud v\ud r= M_0<0.
\]
We suppose that  the potentials $U$ and $W$ satisfy the technical requirements listed in Appendix \ref{app.A}.
Let $0<T<
\infty$. Then up to  a subsequence, still labelled by $\varepsilon$,  the associated solution $f^\varepsilon$ to \eqref{eq.nnn} converges weakly in $L^1((0,T)\times\RR^d\times
\RR^d)$ to $\rho(t,r)M(v)$, $\rho^\varepsilon$ converges to $\rho$ in $C^0([0,T]; L^1(\RR^d)-\textrm{weak})$, with $\rho$ being the solution to \eqref{DD}
and initial data $\rho(t=0,r)$ given by the weak limit in $L^1(\RR^d)$ of $\int f^\varepsilon_{\mathrm{Init}}\ud v$.
\end{theo}

 For the time being, let us discuss the numerical evaluation of the effective coefficient matrices $\mathcal D, \mathcal K$. In the particular case considered here, the right hand side of both equations in \eqref{e:invert} (or equivalently \eqref{e:invert-L}) 
is 
 of the form $v$ times a radial function, then we can simplify the diffusion and drift matrices by taking into
account the symmetries of the problem. We leave the reader to
check the following result.

\begin{lem}\label{symmetry}
Given $v\in \RR^d$, let us define for any indices $i$, $j$ the
linear operators
$$
\mathcal{T}_{ij} (v) : v_i\leftrightarrows v_j \qquad \mbox{and}
\qquad \hat{\mathcal{T}}_{i} (v) : v_i\leftrightarrows -v_i
$$
that exchange $v_i$ and $v_j$ and change $v_i$ by $-v_i$
respectively. Then, the following relations hold:
$$
\chi_i(\mathcal{T}_{ij} (v))=\chi_j(v) \qquad\mbox{and}\qquad
\chi_i(\hat{\mathcal{T}}_{i} (v))=-\chi_i(v)\,,
$$
and
$$
\kappa_i(\mathcal{T}_{ij} (v))=\kappa_j(v) \qquad\mbox{and}\qquad
\kappa_i(\hat{\mathcal{T}}_{i} (v))=-\kappa_i(v)\,.
$$
As a consequence, we deduce that there exist reals $D>0$ and $\kappa$
such that
$$
\mathcal{D}_{ij}=-\int_{\RR^d} v_i \chi_j(v) \ud v = D \delta_{ij}
\qquad \mbox{and} \qquad \mathcal{K}_{ij}=-\int_{\RR^d} v_i
\kappa_j(v) \ud v = \kappa \delta_{ij}\,,
$$
and the macroscopic aggregation-diffusion equation~\eqref{DD} becomes
\begin{equation}\label{cont1}
\partial_t \rho = D \Delta_r\rho +\kappa \nabla_r\cdot((\nabla_r U \star
\rho)\rho\big)\,.
\end{equation}
\end{lem}

Therefore, in order to compute the effective  macroscopic equation \eqref{cont1}, we need to find the solutions to one component of \eqref{e:invert} (or \eqref{e:invert-L}). They are given by explicit formulae  in very few cases. A particular example is given by  the  quadratic potential $W(v)=\tfrac{|v|^2}2$:
then $\chi(v)=\theta\kappa(v)=-v M(v)$. This is due to the fact that $-v M(v)$ is an eigenfunction of the Fokker-Planck operator $Q$ with quadratic potential. 
However,  for more general, non-quadratic potentials such as the one used in the swarming model, it is not possible to obtain explicit formulas for the coefficients of the limiting equation~\eqref{cont1}. In the next section we will study this problem by eigenfunction expansion of the Fokker-Planck operator and we will discuss how to accurately approximate those coefficients.

%%%%%%%%%%%%%%%%%%%%%%%%%%%%%%%%%%%%%%%%%%%%%%%%%%%%

\section{Approximation of the Diffusion and Drift Coefficients}
\label{Section3}

We start by recalling a well-known (unitary) equivalence between Fokker-Planck (in $L^2(\mathbb{R}^d ; M^{-1}\ud v)$) and Schr\"odinger operators (in $L^2(\mathbb{R}^d)$). 
By setting 
$$
\mathcal H(u)=-\ds\frac{\theta}{\sqrt M}\nabla_v\cdot\Big(M\nabla_v \Big(\ds\frac{u}{\sqrt M}\Big)\Big)
=-\ds\frac{1}{\sqrt M}Q\Big(u\sqrt M\Big)
,\qquad
u=\ds\frac{f}{\sqrt M},
$$
it is easy to check that $\mathcal H$ reduces to the Schr\"odinger operator 
$$
 \mathcal{H}(u)= -\theta\Delta u +\Phi(v) u,
$$
with the potential
\begin{equation}\label{e:potential-schrod}
\Phi(v)=-\frac12\Delta_v W (v)+\frac1{4\theta}|\nabla_v W (v)|^2.
\end{equation}
We remark that $\Phi$ is precisely the potential that appears in Assumption~\ref{e:poincare}.  
The  operator $\mathcal{H}=-\theta\Delta+\Phi(v)$
is defined on the domain 
$$
D(\mathcal{H})=\{u\in L^2(\R^d),\ \Phi u\in L^2(\R^d),\ \Delta u\in L^2(\R^d)\}.
$$
We point out that this also defines the
domain of the operator $Q$.
When working with the operator $\mathcal{H}$, the  Lebesgue space
$L^2(\R^d)$ is a natural framework; it what follows, we shall  denote by $(\cdot|\cdot)$ the standard inner product in $L^2(\RR^d)$.
Using classical results for Schr\"odinger operators, see for
instance \cite[Theorem XIII.67]{RS}, we have a spectral
decomposition of the operator $\mathcal{H}$ 
under suitable confining assumptions. 

\begin{lem}
We suppose that
$\Phi\in
L^1_{loc}(\RR^d)$ is bounded from below, and satisfies $\Phi(v)\to+\infty$ as
$|v|\to\infty$. Then,
$\mathcal{H}^{-1}$ is a self-adjoint compact operator in $L^2(\RR^d)$ and
$\mathcal{H}$ admits a spectral decomposition: there exist a non decreasing
sequence of real numbers $\{\lambda_n\}_{n\in\mathbb{N}}\to \infty$, and a
$L^2(\RR^d)$-orthonormal basis $\{\Psi_n\}_{n\in\mathbb N}$ such that
$\mathcal{H}(\Psi_n)=\lambda_n\Psi_n$, $\lambda_0=0$, $\lambda_1\geq
\Lambda>0$. 
\end{lem}

We remark that the spectral gap $\Lambda$ of the Schr\"{o}dinger operator $\cH$ is the Poincar\'{e} constant in the Poincar\'{e} inequality associated to the Fokker-Planck operator $Q$, i.e. the Poincar\'{e} inequality for the probability measure $M(v) \, dv=Z^{-1} e^{-W(v)/\theta} \, dv$.  

\noindent
Notice that
the property \eqref{coer} implies that the operator $\mathcal{H}$
is positive and that the kernel is spanned by $\sqrt{M}$.
We wish to solve the equations 
$$
Q(\chi)=vM \quad\mbox{and}\quad Q(\kappa)=\frac1\theta\nabla_v W M\,,
$$  
or equivalently,
$$
\mathcal{H}(u_\chi)=-v\sqrt{M}\, \mbox{ with } u_\chi=\frac{\chi}{\sqrt M} \qquad \mbox{ and } \qquad 
\mathcal{H}(u_\kappa)=-\frac1\theta\nabla_v W \,\sqrt{M}\, \mbox{ with } u_\kappa=\frac{\kappa}{\sqrt M}\,. 
$$
The diffusion and drift coefficients are defined by the quadratic quantities
\begin{equation}\label{diffmatrix}
\mathcal{D}=-\ds\int v\otimes \chi \ud v= \ds\int \mathcal{H}(u_\chi)\otimes
u_\chi\ud v
\end{equation}
and
\begin{equation}\label{driftmatrix}
\mathcal{K}=\ds-\int v\otimes \kappa \ud v= \ds\int \mathcal{H}(u_\chi)\otimes
u_\kappa\ud v.
\end{equation}

The Schr\"{o}dinger operator is selfadjoint in $L^{2}(\R^{d})$ and under Assumption~\eqref{e:poincare} it has compact resolvent. We denote its eigenvalues and eigenfunctions by $\{\lambda_{n}, \, \Psi_{n} \}_{n=1}^{+\infty}$. A lot of information on the properties of the eigenvalues and eigenfunctions of Schr\"{o}dinger operators is available~\cite{RS}. Using the spectral decomposition of $\cH$ we can obtain a formula for the 
effective
drift and diffusion coefficients. These formulas are similar to the Kipnis-Varadhan formula for the diffusion coefficient in the functional central limit theorem for additive functionals of reversible Markov processes~\cite{kipnis}. We can use these formulas to develop a numerical scheme for the approximate calculation of the drift and diffusion coefficients. Indeed, if we develop $\mathcal{H}(u_\chi)$, $u_\chi$, $\mathcal{H}(u_\kappa)$ and $u_\kappa$ in the eigenbasis, we get
$$
\mathcal{H}(u_\chi) = \sum_{k=1}^\infty \eta_k \Psi_k\,,\quad u_\chi=(u_\chi | \Psi_0)
\Psi_0+\sum_{k=1}^\infty \frac{\eta_k}{\lambda_k} \Psi_k\,,
$$
and
$$
\mathcal{H}(u_\kappa) = \sum_{k=1}^\infty \omega_k \Psi_k\,,\quad u_\kappa=(u_\kappa | \Psi_0) \Psi_0+\sum_{k=1}^\infty \frac{\omega_k}{\lambda_k} \Psi_k\,,
$$
with
$$
\eta_k=\ds\int \mathcal{H}(u_\chi) \Psi_k\ud v \quad \mbox{and} \quad \omega_k=\ds\int \mathcal{H}(u_\kappa) \Psi_k\ud v\,, \quad\mbox{for all } k\in\mathbb N.
$$
Substituting in \eqref{diffmatrix} and \eqref{driftmatrix}, we obtain the following formulas for the diffusion and drift matrices:
\begin{equation}\label{diffdriftmatrix}
\mathcal{D}=\sum_{k=1}^\infty \frac{\eta_k\otimes\eta_k}{\lambda_k}
\qquad
\mbox{and}
\qquad
\mathcal{K}=\sum_{k=1}^\infty \frac{\eta_k\otimes\omega_k}{\lambda_k}\,.
\end{equation}
We can obtain approximate formulas for the drift and diffusion coefficients by truncating the data: given $h\in L^2(\RR^d)$ and denoting
$$
h=\sum_{k=0}^\infty \zeta_k \Psi_k\,,\qquad \zeta_k=\ds\int h\Psi_k\ud v\,, \quad\mbox{for all } k\in\mathbb N\,,
$$
we set
$$
h^N=\sum_{k=0}^N \zeta_k \Psi_k.
$$
We denote by $u^N$ the solution of $\mathcal{H}(u^N)=h^N$, and
compare it to $u$, the solution of $\mathcal{H}(u)=h$. We
use the Sobolev and Cauchy-Schwarz  inequalities as follows
\[
\Lambda \|u-u^N\|^2_{L^2} \leq \big( \mathcal{H}(u-u^N) |
u-u^N\big)=(h-h^N | u -u^N)\leq \|h-h^N\|_{L^2}\ \|u-u^N\|_{L^2}.
\]
Hence, we get
\begin{equation*}%\label{errorestim}
\|u-u^N\|_{L^2} \leq \frac1{\Lambda} \|h-h^N\|_{L^2} \,.
\end{equation*}
The accuracy of the approximation is therefore driven by the accuracy of 
the approximation of the data $h$ by its truncated Fourier series: we need 
information on the behavior of  the eigenvalues $\lambda_n$ for large $n$ and on the 
accuracy of the spectral projection of $h$. We proceed by analogy to the 
standard theory of Fourier series, where the behavior of the Fourier 
coefficients is related to the regularity of the function. More precisely, 
the estimate
\begin{equation}\label{decay}
 \|h-h^N\|_{L^2}\leq C\lambda_{N+1}^{-k}
\end{equation}
holds for some $k>0$. Assume $\mathcal{H}^k (h)\in
L^2(\RR^d)$; by the spectral decomposition, we obtain
$$
\mathcal{H}^k (h) = \sum_{n=0}^\infty (\mathcal{H}^k (h) | \Psi_n)
\Psi_n \qquad \mbox{with } \sum_{n=0}^\infty |(\mathcal{H}^k (h) |
\Psi_n)|^2 <\infty \,.
$$
Now, we estimate the difference as
$$
\|h-h^N\|_{L^2}^2 = \sum_{n>N} |(h | \Psi_n)|^2 = \sum_{n>N}
\frac1{\lambda_n^{2k}} |(\mathcal{H}^k (h) | \Psi_n)|^2,
$$
where we used that $\mathcal{H}(\Psi_n)=\lambda_n\Psi_n$ and the
fact that $\mathcal{H}$ is self-adjoint. Therefore, we deduce that
$$
\|h-h^N\|_{L^2}^2 \leq \frac1{\lambda_{N+1}^{2k}} \sum_{n>N}
|(\mathcal{H}^k (h) | \Psi_n)|^2 \leq \frac1{\lambda_{N+1}^{2k}}
\|\mathcal{H}^k (h)\|_{L^2}^2,
$$
since the eigenvalues are in increasing order, leading to the
desired estimate \eqref{decay}. A similar argument shows that if
$\mathcal{H}^{k+1} (h)\in L^2(\RR^d)$, then
$$
\|\mathcal{H}(h-h^N)\|_{L^2} \leq \frac1{\lambda_{N+1}^{k}}
\|\mathcal{H}^{k+1} (h)\|_{L^2}\,.
$$
A direct application of the strategy above to 
$$
v\mapsto h(v)=-v\sqrt{M(v)} \qquad \mbox{ and to }
\qquad v\mapsto h(v)=-\frac1\theta \nabla_v M \sqrt{M(v)}\,,
$$ 
which satisfy $\mathcal{H}^k (h)\in L^2(\RR^d)$ for all
$k>0$, together with the symmetry of the potential in Lemma \eqref{symmetry}
leads to the main result of this section estimating the
error due to the truncation in \eqref{diffdriftmatrix}.

\begin{theo}
Given $\mathcal{D}^N$ and $\mathcal{K}^N$ the truncated diffusion and drift coefficients defined
by
$$
\mathcal{D}^N= \ds\int  \mathcal{H}(u_\chi^N)\otimes u_\chi^N\ud
v\,= D^N\mathbb I,\qquad  \mbox{and} \qquad \mathcal{K}^N= \ds\int  \mathcal{H}(u_\chi^N)\otimes u_\kappa^N\ud
v\,= K^N\mathbb I,
$$
with $\mathcal{H}(u_\chi^N)(v)=(-v\sqrt{M(v)})^N$ and $\mathcal{H}(u_\kappa^N)(v)=(-\tfrac1\theta \nabla_v M \sqrt{M(v)})^N$, or equivalently 
\begin{equation}\label{diffdriftcoef}
D^N=\frac1d \sum_{k=1}^N  \frac{|\eta_k|^2}{\lambda_k}
\qquad
\mbox{and}
\qquad
K^N=\frac1d \sum_{k=1}^N \frac{\eta_k\cdot\omega_k}{\lambda_k}\,,
\end{equation}
then the following error
estimate holds: for all $k>0$ and all $N\in\mathbb{N}$, there
exists $C_k>0$ (depending on $k$ but not on $N$) such that
$$
|D-D^N|+|K-K^N|\leq C_k\lambda_{N+1}^{-k}\,.
$$
\end{theo}

%%%%%%%%%%%%%%%%%%%%%%%%%%%%%%%%%%%%

\section{Numerical Approximation and Simulations}
\label{Section4}

The numerical method in practice works as follows:
\begin{itemize}
\item {\bf Step 1.- $\mathcal H\rightarrow \mathcal H^R$}: We consider the problem set on $[-R,+R]^d$, with $R\gg 1$ and completed with homogeneous Dirichlet boundary conditions:
owing to the functional framework, we expect that the eigenfunctions are localized around the origin and are exponentially decreasing far away the wells of the potential, see \cite{Agmon}. In particular, it holds under our assumptions on the behavior of the potential as $|v| \rightarrow + \infty$, from \cite[Thm.~3.4, Thm.~3.10]{HislopSigal1996}. Then, we choose $R$ large enough  to reduce the truncation error. In our examples, we fix $R=10$. 

\item {\bf Step 2.- $\mathcal H^R\rightarrow \mathcal H^{R,\mu}$}: 
We use Finite Elements methods to discretize the operator $\mathcal H^R$. In practice, we use the library {\sc M\'elina} \cite{Melina}, a uniform mesh of $[-R,R]$ with 1000 uniform $\mathbb P_{10}$ elements and a quadrature of degree 21.
We have chosen this method since we need an approximation of the solution of the PDE with high order accuracy, see Remark~\ref{rem.toto}. 
Here and  in what follows we denote by  $\mu>0$ a measure of the accuracy of the underlying discretization method. It thus contains the information both on the  refinement of  the mesh, and the degree of the piecewise reconstruction.

\item {\bf Step 3.- $(\lambda_n,\Psi_n)\rightarrow (\lambda_n^{R,\mu},\Psi_n^{R,\mu})$}:
Having at hand the discrete operator on the truncated domain, denoted  $\mathcal H^{R,\mu}$, we determine its first $N$ eigenelements
$\big((\lambda_1^{R,\mu},\Psi_1^{R,\mu}),...,(\lambda_N^{R,\mu},\Psi_N^{R,\mu})\big)$.
In the Finite Elements framework, the eigenvectors $\Psi_n^{R,\mu}$ are piecewise polynomials functions approximating  the $n$th eigenfunctions.
We recall that the eigenvectors form an orthonormal family.

\item {\bf Step 4.- $(\eta_n,\omega_{n})\rightarrow (\eta_n^{R,\mu},\omega_n^{R,\mu})$}: Given the data 
$$\mathcal{H}(u_\chi)(v) = -v\sqrt{M(v)}\qquad\mbox{ and }\qquad \mathcal{H}(u_\kappa)(v) = -\frac{1}{\theta} \nabla_v W(v) \sqrt{M(v)},$$
we compute the corresponding $N$ Fourier coefficients 
by using an appropriate quadrature formula, depending on the approximation framework, for the discrete analogue $(\eta_n^{R,\mu},\omega_n^{R,\mu})$ of
$$
\eta_n=\ds\int \mathcal{H}(u_\chi)(v) \Psi_n(v)\ud v,\qquad
\omega_n=\ds\int \mathcal{H}(u_\kappa)(v) \Psi_n(v)\ud v.
$$
Our results here are computed using a simple composite rectangular rule.

\item {\bf Step 5.- $(D,K)\rightarrow (D^{R,\mu,N},K^{R,\mu,N})$}: We now approximate the diffusion and drift coefficients using \eqref{diffdriftcoef} to conclude
\begin{equation}\label{eq.DKapprox}
D^{R,\mu,N}=\ds\frac 1d \ds\sum_{n=1}^N \ds\frac{|\eta_n^{R,\mu}|^2}{\lambda_n^{R,\mu}}
\qquad
\mbox{and}
\qquad
K^{R,\mu,N}=\ds\frac 1d \ds\sum_{n=1}^N \ds\frac{\eta_n^{R,\mu}\cdot\rho_n^{R,\mu} }{\lambda_n^{R,\mu}}.
\end{equation}
\end{itemize}
Once the eigenelements are known, the computational cost of the evaluation of the coefficients $(\eta_n^{R,\mu},\omega_n^{R,\mu})$ is linear with respect to the size of the linear problem to be solved (that depends directly on $\mu$). Hence, the main source of the computational cost relies on the determination of the $N$ eigenpairs. 

For the potentials that we consider in this paper, Lanczos-like algorithms can be used. As an iterative method, its computational cost cannot be estimated a priori. Nevertheless, we expect that only a few eigenpairs can provide an accurate result (for the quadratic case, the problem is exactly solved with  the first eigenpair associated with a  positive eigenvalue), so that the resolution would be far less costly than solving the linear system, a problem that, for small $\mu$'s,  would also require iterative methods. Here we use standard Lanczos techniques; we refer the reader to~\cite{LS, Sa, So} for further information on these methods and to \cite{BDMV07} for the computation of the first few eigenpairs of complicated Schr\"odinger operators, based on Finite Elements approximations.

We will show numerical simulations for three different potentials in one dimension given by:
\begin{itemize}
\item {\bf Case A.-} The symmetric smooth potential given by
$W(v)=\tfrac{1}{4\gamma} v^4 - \tfrac{1}2 v^2\ $ with $\gamma>0$. 
\item {\bf Case B.-} The symmetric singular potential given by
$W(v)=\tfrac{1}{4\gamma} v^4 - \tfrac{1}3 |v|^3\ $ with $\gamma>0$.
\item {\bf Case C.-} The tilted smooth potential given by
$W(v)=\tfrac{1}{4\gamma} v^4 - \tfrac{1}2 v^2 -\delta v\ $ with $\gamma,\delta>0$.
\end{itemize}
We can gather all of them in a single potential
\begin{equation}\label{eq.defWgeneral}
W(v)=\frac{1}{4\gamma} v^4 - \frac{\sigma}3 |v|^3 - \frac{1-\sigma}2 v^2 -\delta v \qquad \mbox{with }\quad \gamma>0,\ \delta\geq 0,\ \sigma\in\{0,1\}\,.
\end{equation}
Note that in Case C (or $\sigma=0$, $\delta>0$ in \eqref{eq.defWgeneral}), the potential $W$ is not symmetric and the compatibility condition for solving the auxiliary equation \eqref{invert-1} is not satisfied. We shall see how the theory can be adapted to this case (see Section~\ref{Sec.CASC}).

\begin{rmk}
Note that in the one dimensional case, the drift coefficient can be expressed in a simpler form. This is due to the fact that we can solve explicitly the one dimensional Poisson equation, up to quadratures~\cite[Sec. 13.6]{PavlSt08}. Indeed, according to \eqref{e:invert-L} and the expression of the operator $Q$ in \eqref{formfp}, we notice that 
\begin{align*}
Q(\psi)
&=\frac1\theta \frac{\ud W}{\ud v}(v) M(v) = -\frac{\ud M}{\ud v}(v) 
=\theta \frac{\ud }{\ud v}\left(M\frac{\ud }{\ud v}\left(\frac{\psi}M\right) \right)(v). 
\end{align*}
Then, direct integration yields 
$$\psi(v)= \frac1\theta\left(-v+\int_{\RR}v M(v)\ud v\right)M(v)\ud v.$$
Therefore the drift coefficient defined in \eqref{DD2} becomes
\begin{equation}\label{eq.drift1D}
K = \frac1\theta\int_{\RR}\left(-v+\int_{\RR}v M(v)\ud v\right)^2M(v)\ud v.
\end{equation}
In cases A and B, we have $\int_{\RR}v M(v)\ud v=0$ by symmetry and $K$ is given by $K= \frac1\theta\int_{\RR}v^2M(v)\ud v$.
This explicit formula will be used to check the accuracy of the method. 
\end{rmk}

In order to reduce the number of free parameters, we rescale the velocity by defining $v=\sqrt{\gamma}\tilde v$ into the Fokker-Planck operator in \eqref{formfp}, to get
$$
Q (f) = \div_v \left[ \nabla_v \widetilde{W}(v) f+\frac{\theta}{\gamma} \nabla_v \,f \right],
$$
where we have dropped the tildes for notational simplicity, with the rescaled potential
$$
\widetilde{W}(v) = \frac{1}{4} v^4 - \frac{\sigma\sqrt{\gamma}}3 |v|^3 - \frac{1-\sigma}2 v^2 -\frac{\delta}{\sqrt{\gamma}} v \qquad \mbox{with } \gamma>0,\ \delta\geq 0,\ \sigma\in\{0,1\}\,.
$$
In this way, $\theta\to 0$ and $\gamma\to \infty$ play the same role. In fact, for the symmetric smooth potential of Case A ($\sigma=0$), all terms involving $\gamma$ disappear in the rescaled potential and we can remove one parameter by setting $\theta=1$.
In our simulations, we consider the Schr\"odinger operator 
$$ \widetilde{\mathcal{H}}(u)= -\frac{\theta}{\gamma} \frac{\partial^2 u}{\partial v^2} +\widetilde{\Phi}(v) u,$$ 
with the potential
\begin{align*}
\widetilde{\Phi}(v)
&=-\frac12\frac{\partial^2 \widetilde{W} }{\partial v^2} (v)+\frac{\gamma}{4\theta}\left(\frac{\partial \widetilde{W} }{\partial v}(v)\right)^2 \nonumber\\
&= -\frac12 \left( 3v^2 -2 \sigma\sqrt{\gamma} |v| -(1-\sigma) \right)  +\frac{\gamma}{4\theta}\left( v^3- \sigma\sqrt{\gamma} v|v|-(1-\sigma)v-\frac{\delta}{\sqrt{\gamma}} \right)^2\,,
%\label{potentialschr}
\end{align*}
in each of the different cases above and denote by $\lambda_{k}(\gamma)$ the $k$-th positive eigenvalue.

The first eigenvalue for $\widetilde{\cal H}$ is simple, equal to 0. 

In Cases A and B, the potentials are spherically symmetric. Then the second eigenvalue $\lambda_{1}(\gamma)$ tends exponentially to 0 as $\gamma$ tends to $\infty$~\cite{Rob87,HeSj84} and~\cite[Thm. 1.5]{Simon1984}.  In fact, a careful reading of these references gives that $\lambda_{1}(\gamma)\simeq e^{-c \, \gamma}$ for some positive constant $c$ and $\lambda_{2}(\gamma)\simeq O(1)$. This is a manifestation of the tunnelling effect. 
This behavior leads to numerical difficulties. Indeed, we have to capture two simple eigenvalues $(0,\lambda_{1}(\gamma))$ but with $\lambda_{1}(\gamma)\to 0$ exponentially fast as $\gamma\to\infty$. The first eigenfunction is symmetric and the second one antisymmetric. Numerically, when $\gamma$ is very large, the gap between $0$ and $\lambda_{1}(\gamma)$ becomes negligible compared to the order of the accuracy of the method or even compared to machine precision. Then numerically it appears as if the problem has a double eigenvalue. Then the computation breaks down. Similar difficulties appear for the magnetic tunnelling effect, see \cite{BHR}. A good way to determine whether or not the computation is accurate is to look at the eigenfunction: as soon as the symmetry is broken for the first two eigenfunctions, the computation is wrong.

\begin{rmk}\label{rem.toto}
We have also tested the method by using the standard Finite Difference discretization. We roughly obtain similar results for small values of $\gamma$'s and with the same number of numerical unknowns as for the Finite Elements algorithm (which means with a very refined grid for the finite discretization method).
Discrepancies appear as $\gamma$ increases: the loss of symmetry of the eigenfunctions is sensitive earlier. 
This is reminiscent to the well known fact for similar problems that increasing the degree of polynomials involved in the approximation ($p$-extension) is more efficient than refining the mesh ($h$-extension), see \cite{Ains04,BDMV07}.
\end{rmk}

%\newpage
\subsection{Case A}  
\begin{figure}[h!t]
\begin{center}
\begin{tabular}{ccccc}
$\gamma=1$ & 10 & 50 & 100 & 120\\
\includegraphics[height=2.5cm]{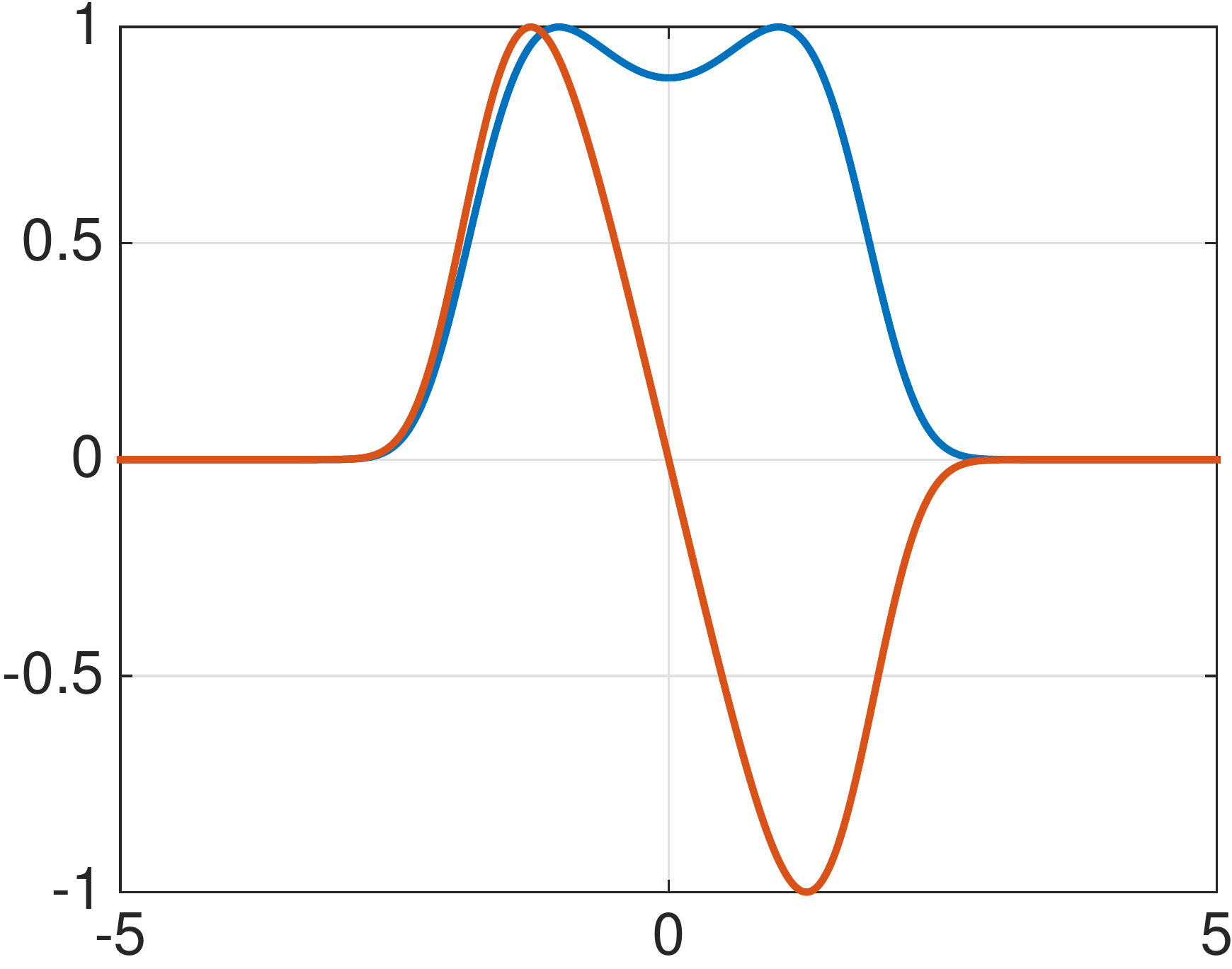}
&\includegraphics[height=2.5cm]{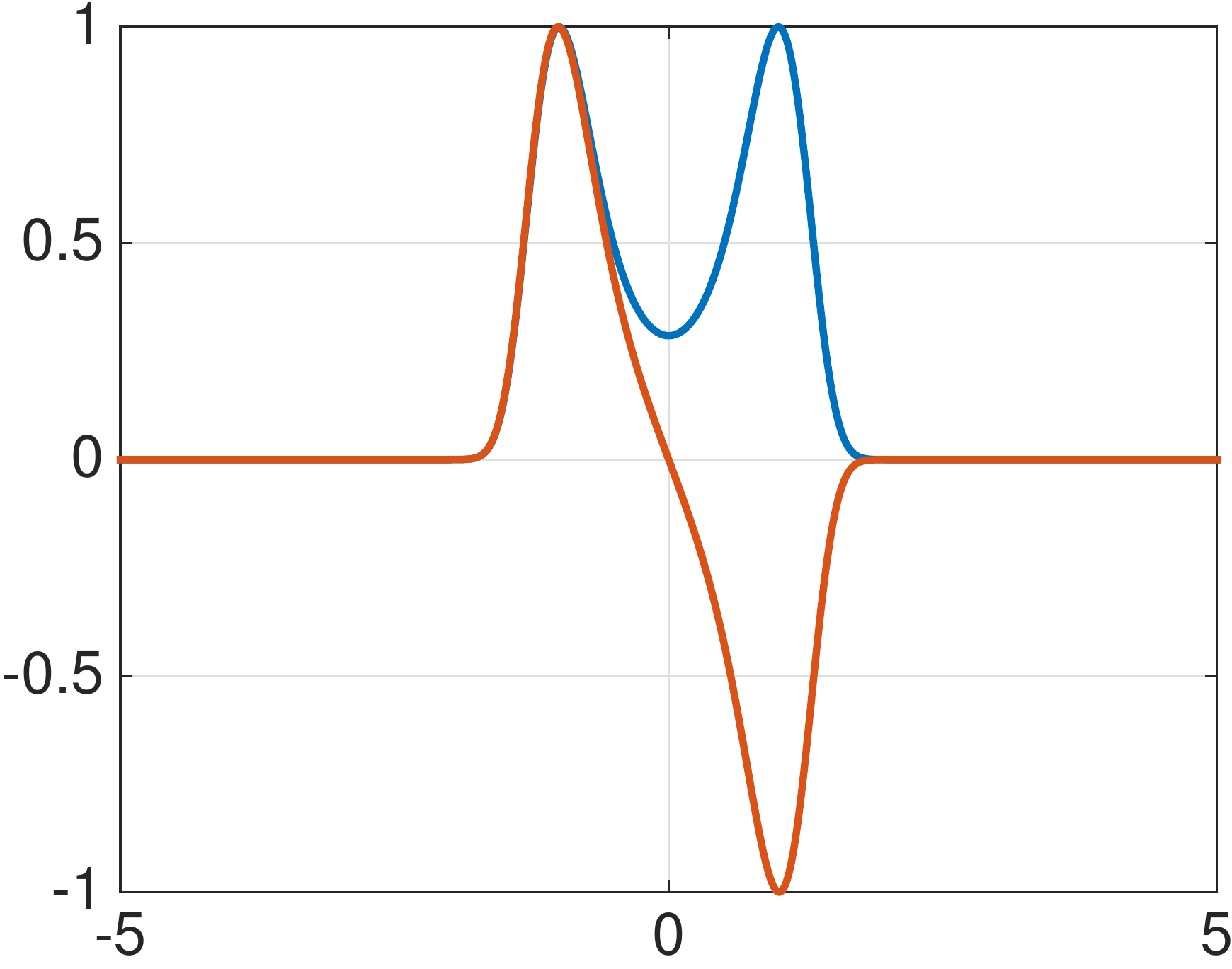}
&\includegraphics[height=2.5cm]{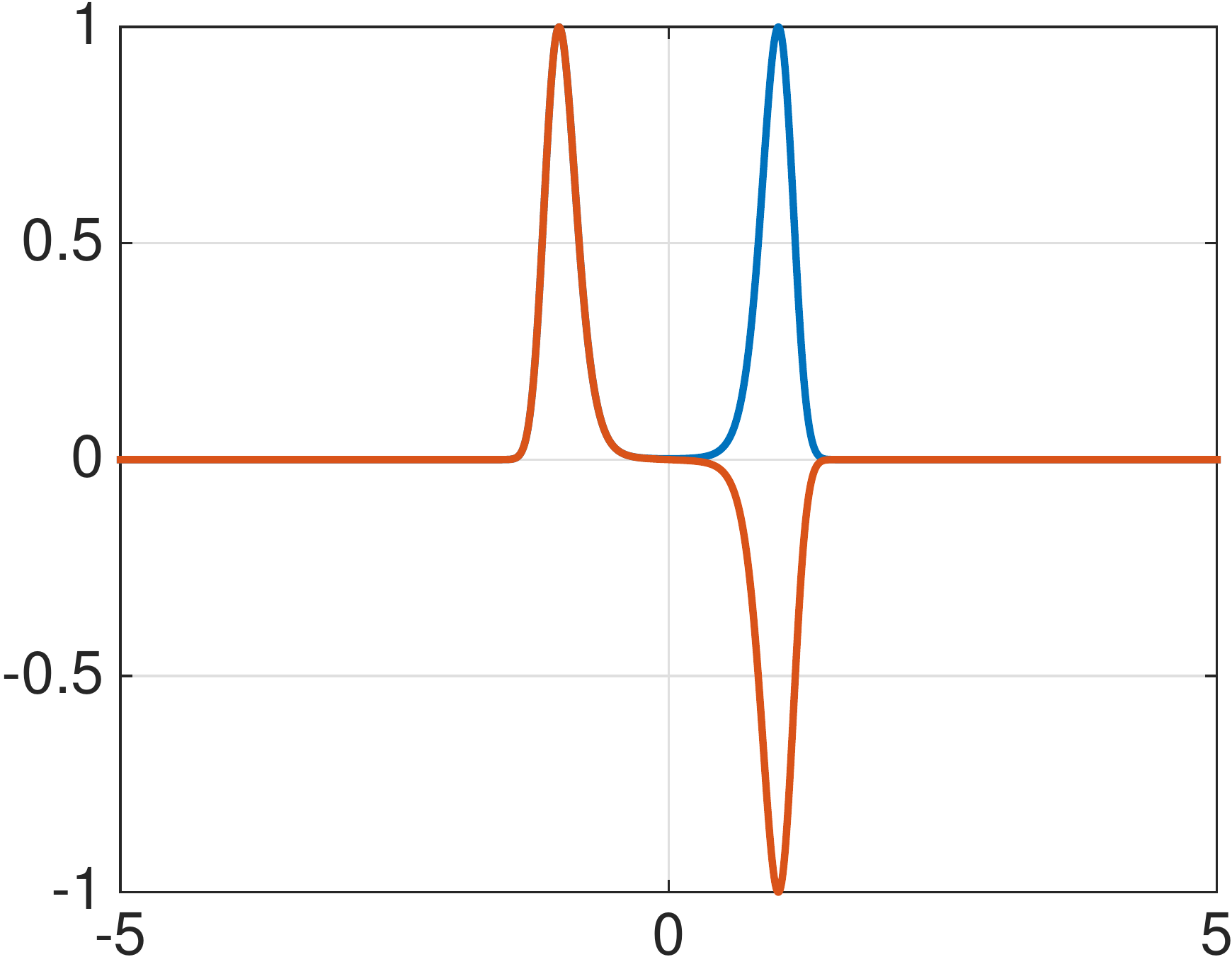}
&\includegraphics[height=2.5cm]{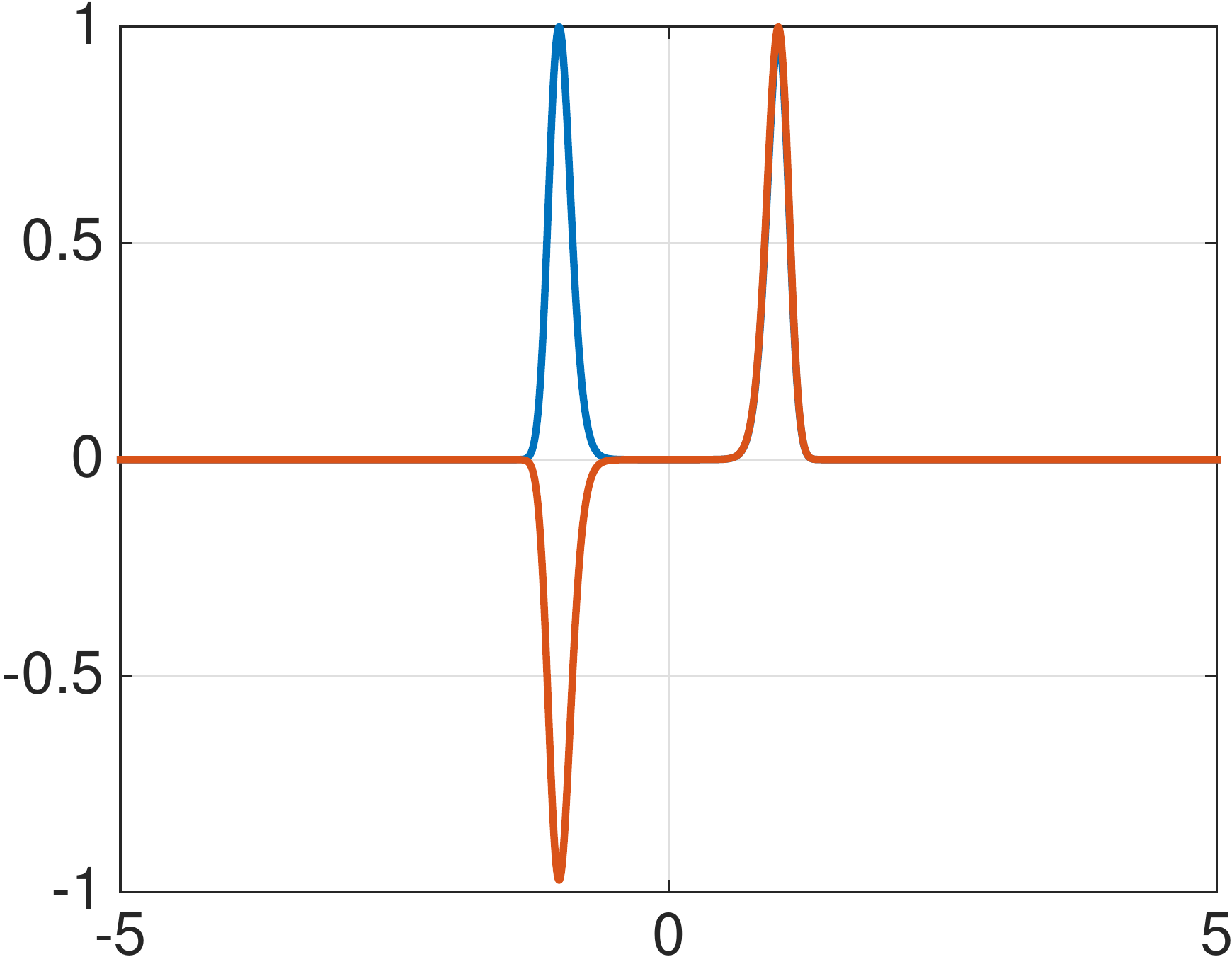}
&\includegraphics[height=2.5cm]{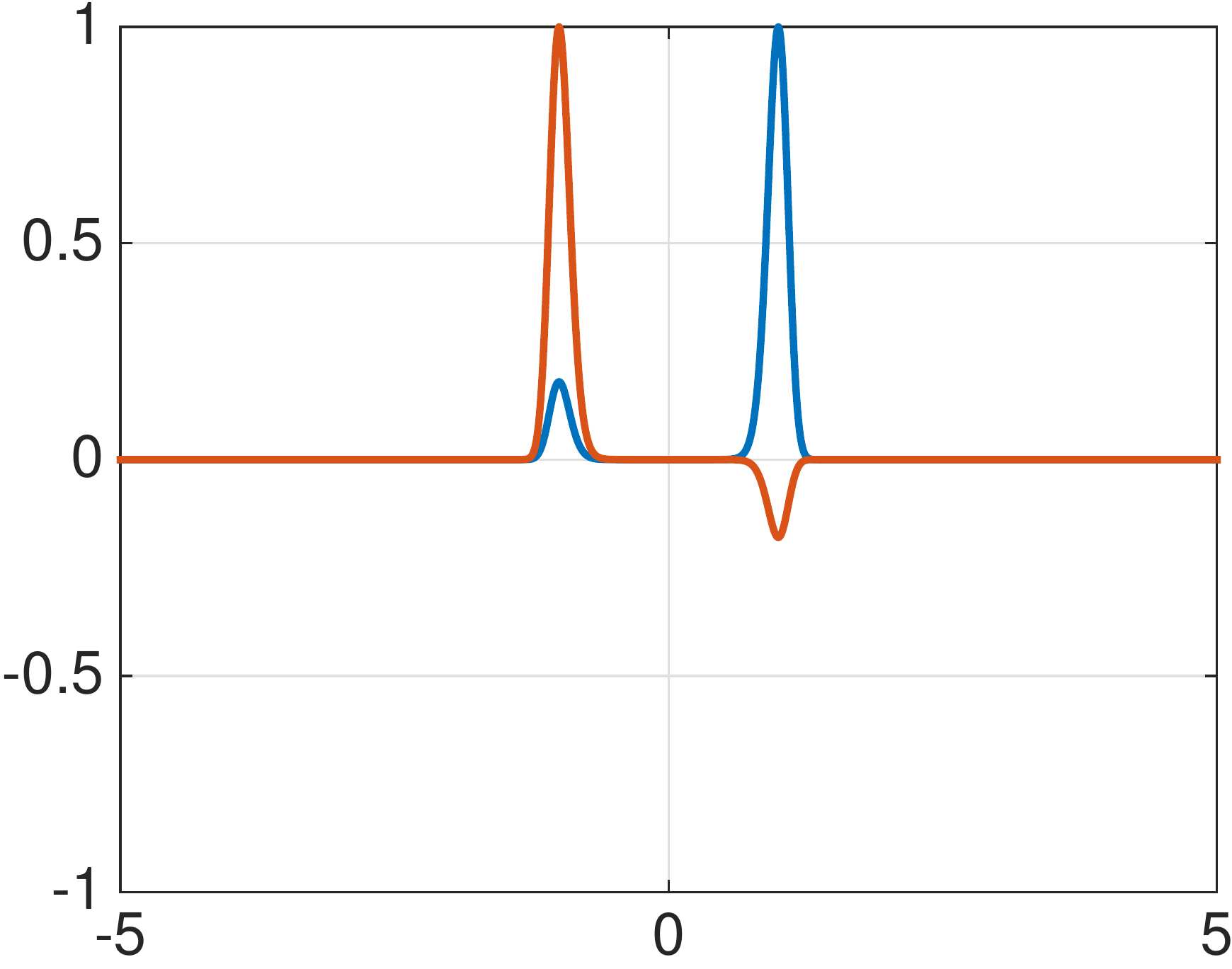}
\end{tabular}
\caption{First two eigenfunctions on $[-R/2,R/2]$.\label{fig.casAVecP}}
\end{center}
\end{figure}

In Figure~\ref{fig.casAVecP} we present the first two eigenfunctions for $\gamma=1, 10, 50, 100, 120$ on $[-R/2,R/2]$. 
We can observe the localization of the eigenfunction and the exponential decay far away from the wells of the potential.
Looking at the symmetry of the eigenfunction, we see that the computation becomes problematic for $\gamma>100$ (the eigenfunctions for $\gamma=120$ are neither symmetric nor antisymmetric). 

\begin{figure}[h!t]
\begin{center}
\begin{tabular}{cc}
\includegraphics[height=3cm]{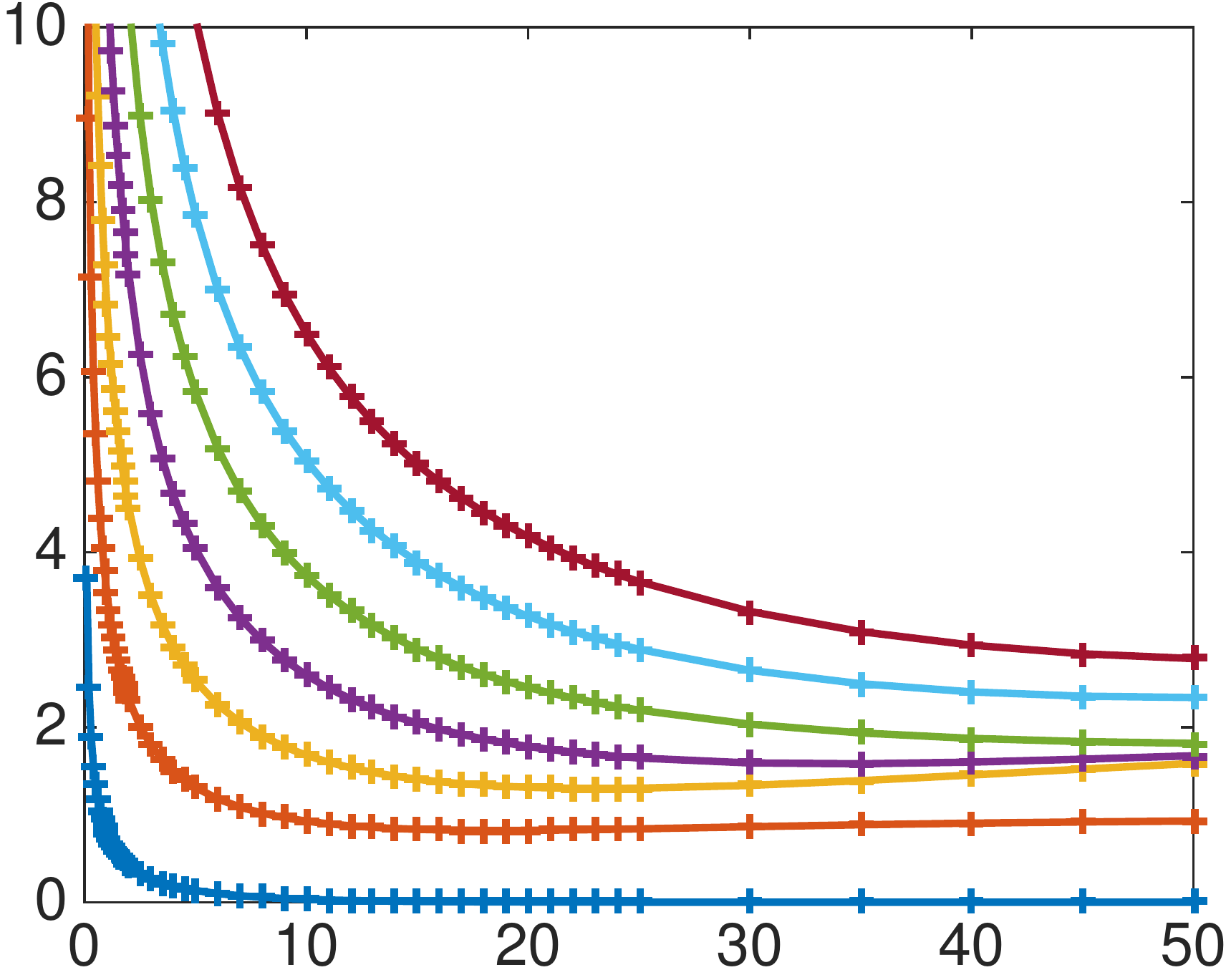}
&$\qquad$\includegraphics[height=3cm]{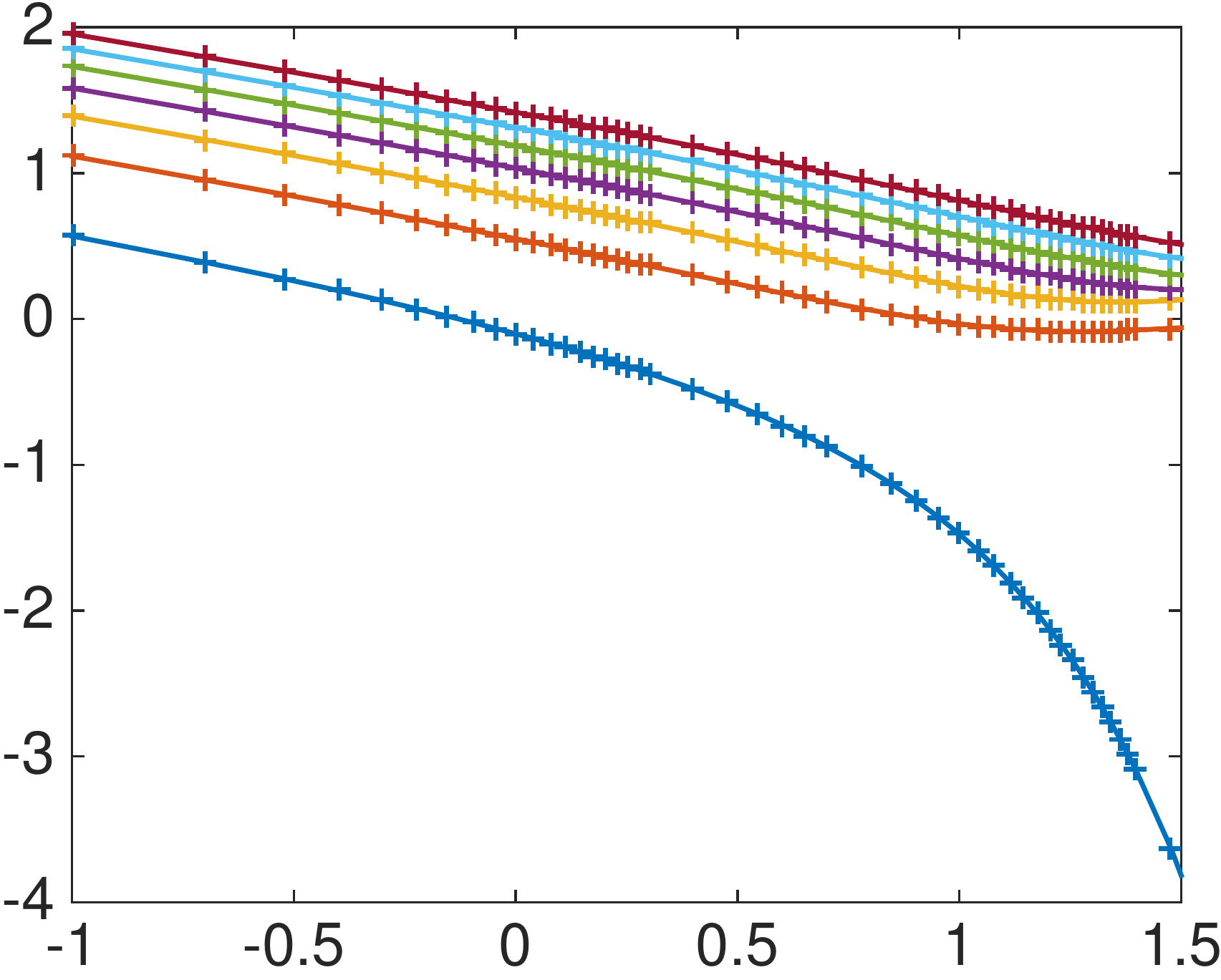}\\
$\gamma\mapsto\lambda_{j}(\gamma)$ & $\qquad\log\gamma\mapsto\log_{10}\lambda_{j}(\gamma)$ 
\end{tabular}
\caption{Convergence of the lowest eigenvalues with respect to $\gamma$.\label{fig.casAVP}}
\end{center}
\end{figure}

In Figure~\ref{fig.casAVP}, we plot the convergence of the positive eigenvalues: $\gamma\mapsto \lambda_j(\gamma)$ and $\log\gamma\mapsto \log_{10}\lambda_j(\gamma)$. We clearly observe the tunnelling effect, with the first eigenvalue $\lambda_{1}(\gamma)$ converging exponentially fast to $0$ while the higher eigenvalues remain of order 1.

\begin{figure}[h!t]
\begin{center}
\begin{tabular}{cc|c}
\includegraphics[height=3.2cm]{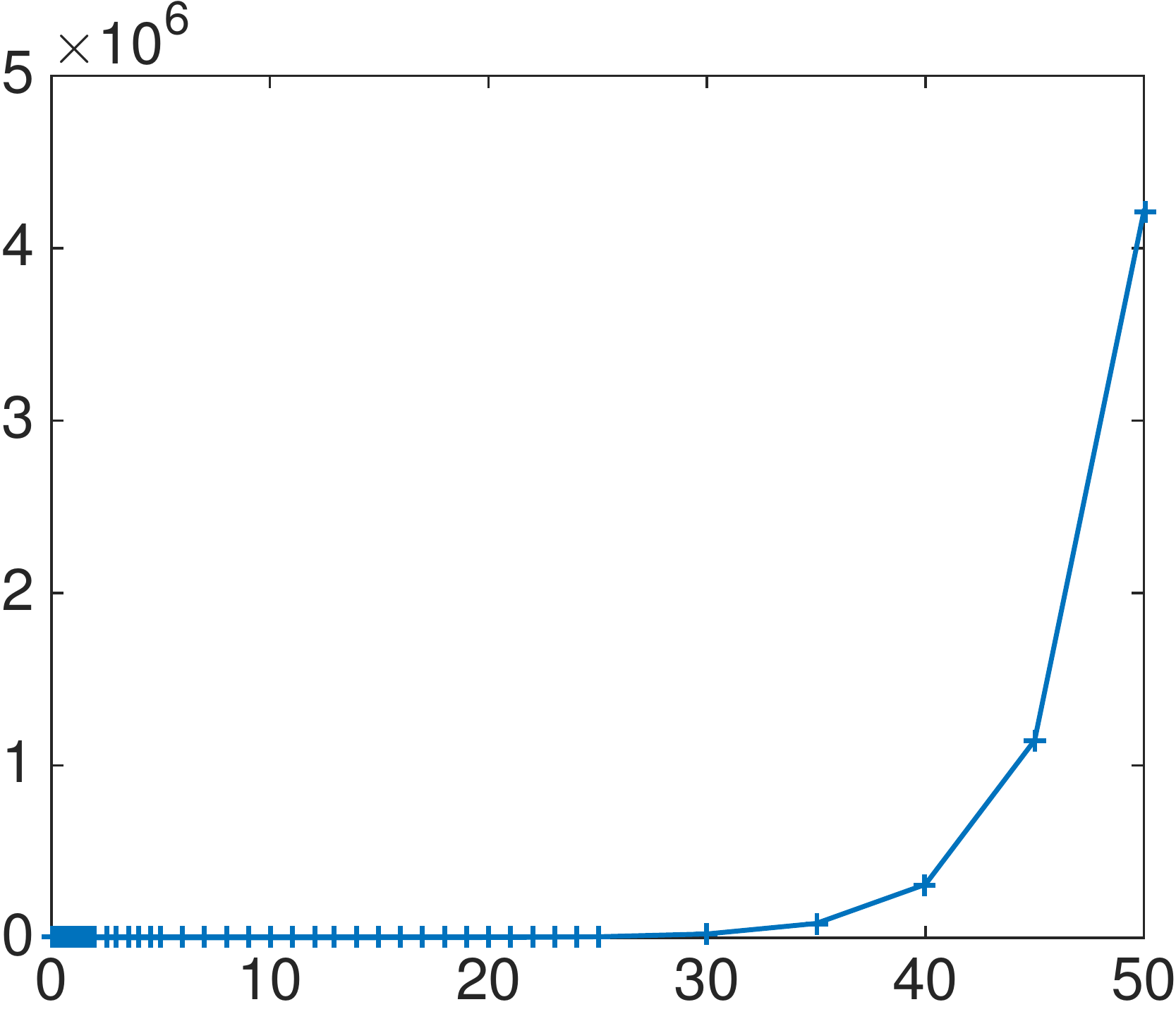}
\qquad&\includegraphics[height=3.2cm]{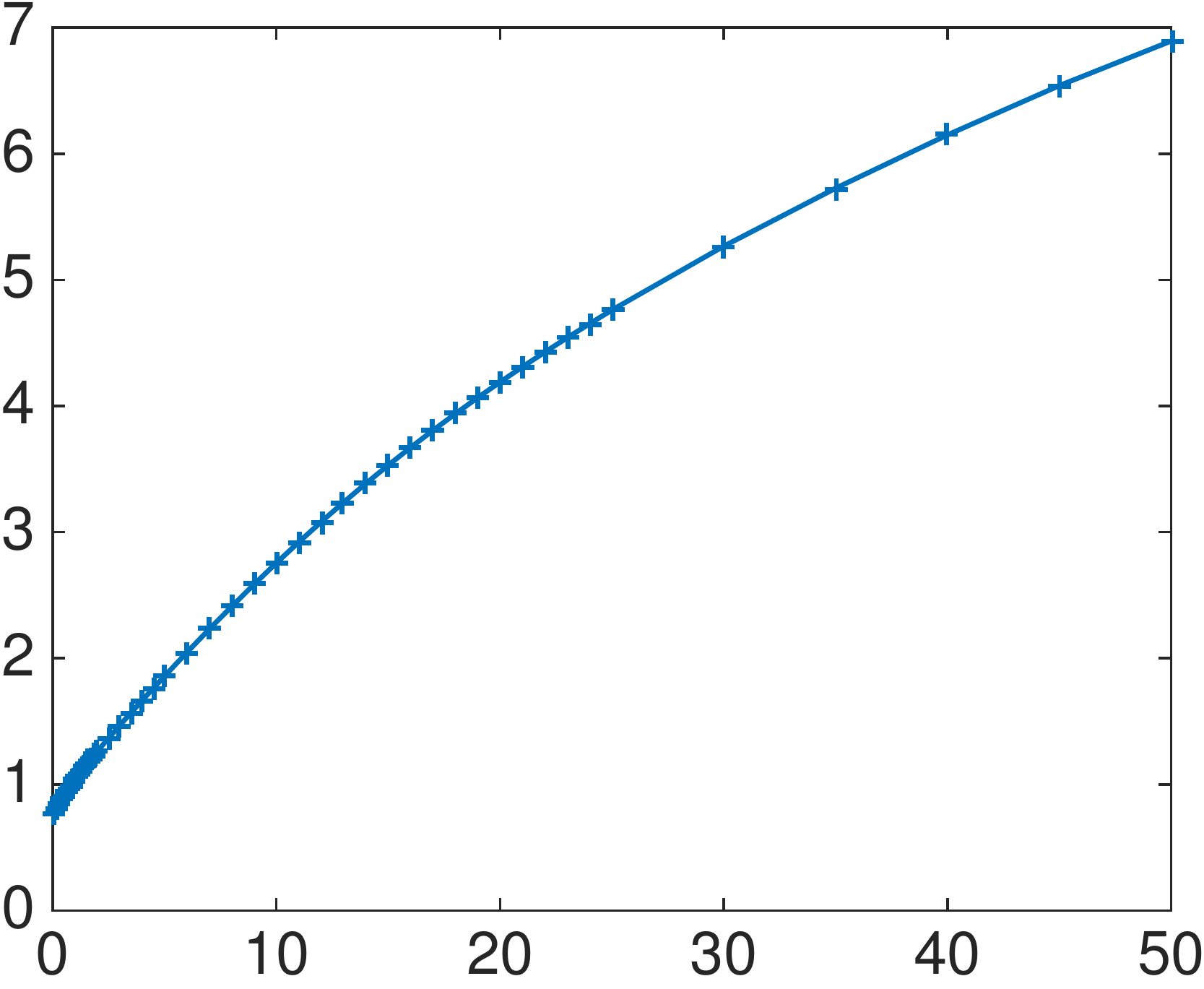}
\qquad&\quad\includegraphics[height=3.2cm]{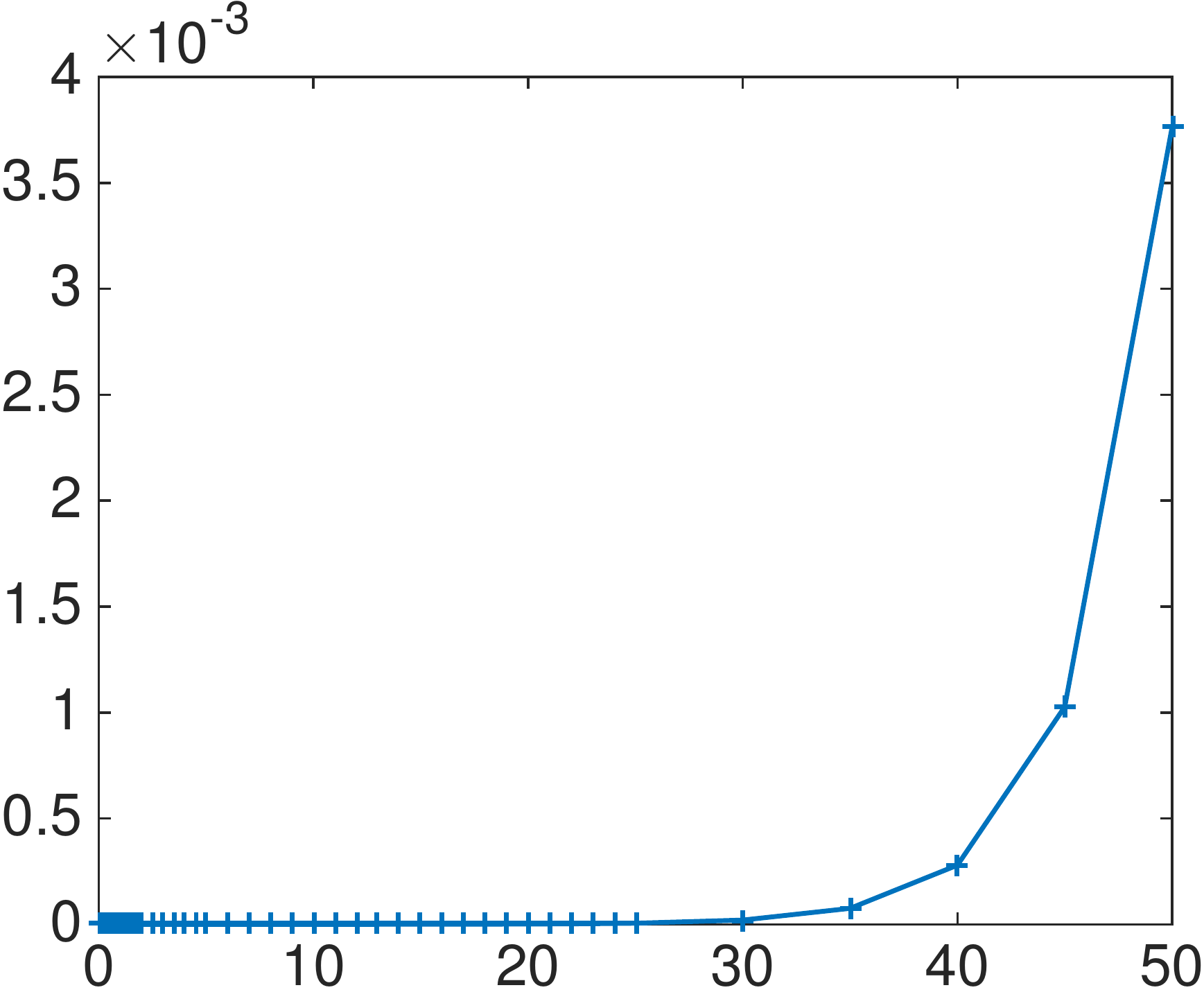}\\
$\gamma\mapsto {\mathcal D}(\gamma)$ & $\gamma\mapsto{\mathcal K}(\gamma)$ & $\gamma\mapsto \frac{|\mathcal K_{*}(\gamma)-\mathcal K(\gamma)|}{\mathcal K_{*}(\gamma)}$ \\[5pt]
\includegraphics[height=3.2cm]{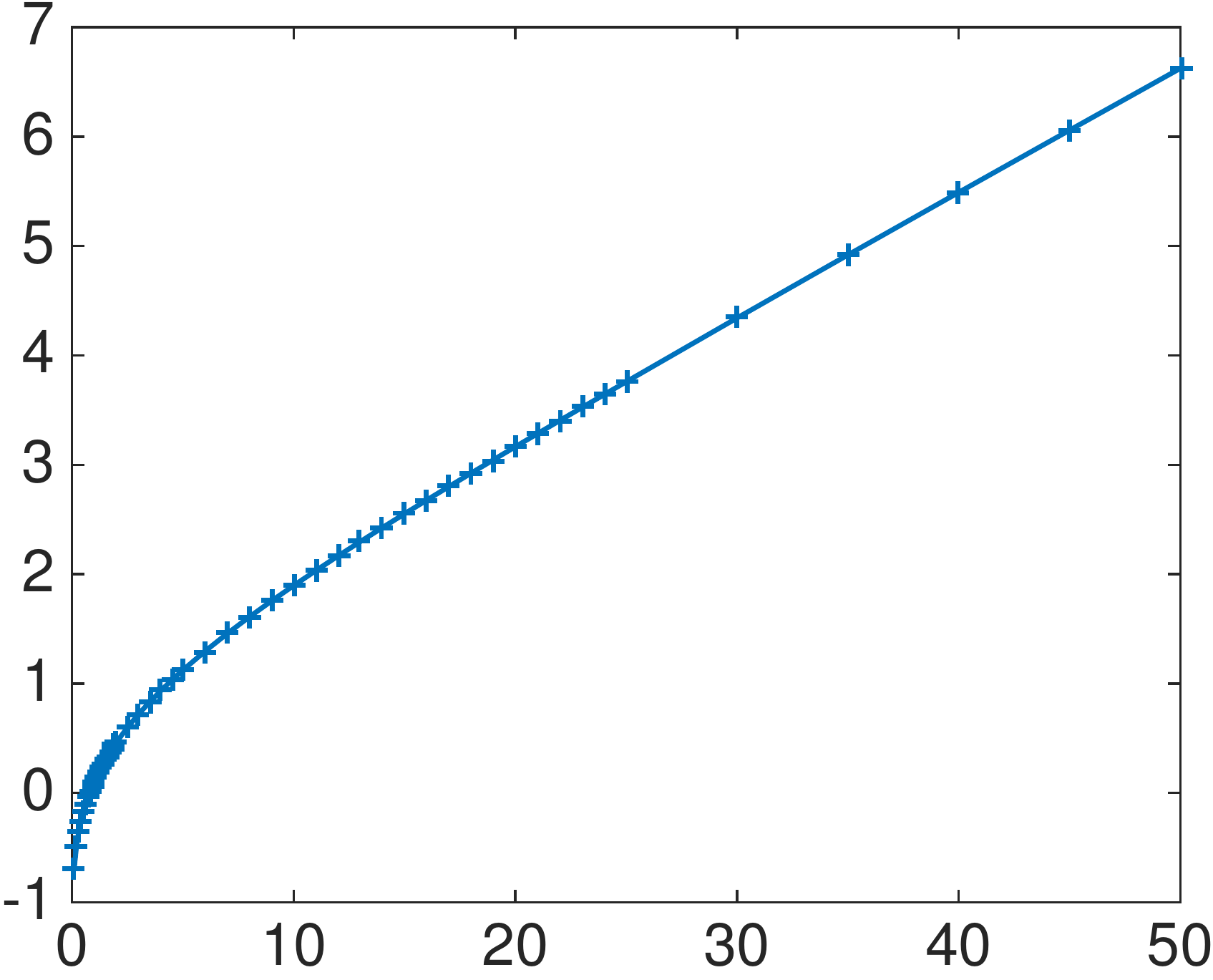}
\qquad&\includegraphics[height=3.2cm]{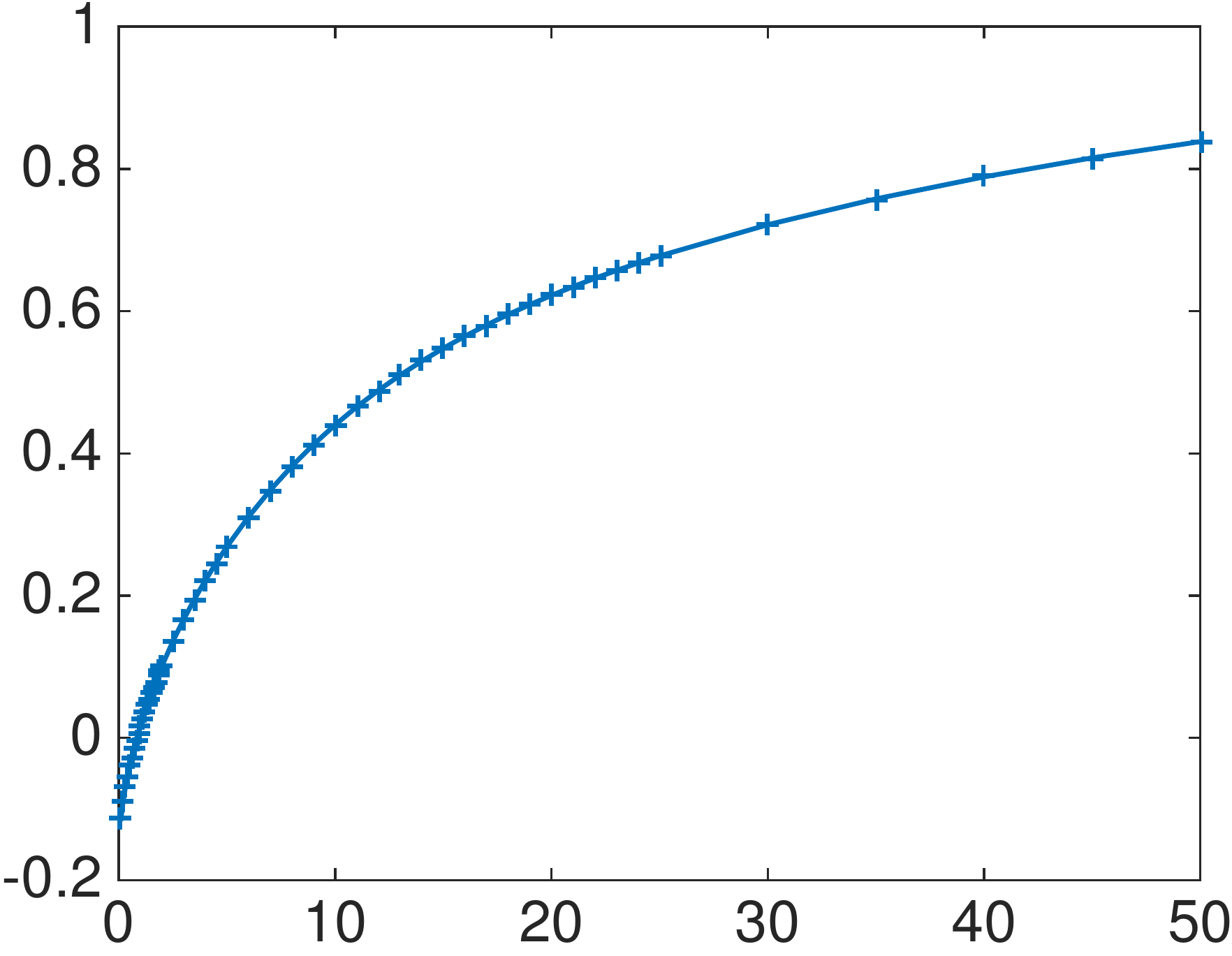}
\qquad&\quad\includegraphics[height=3.2cm]{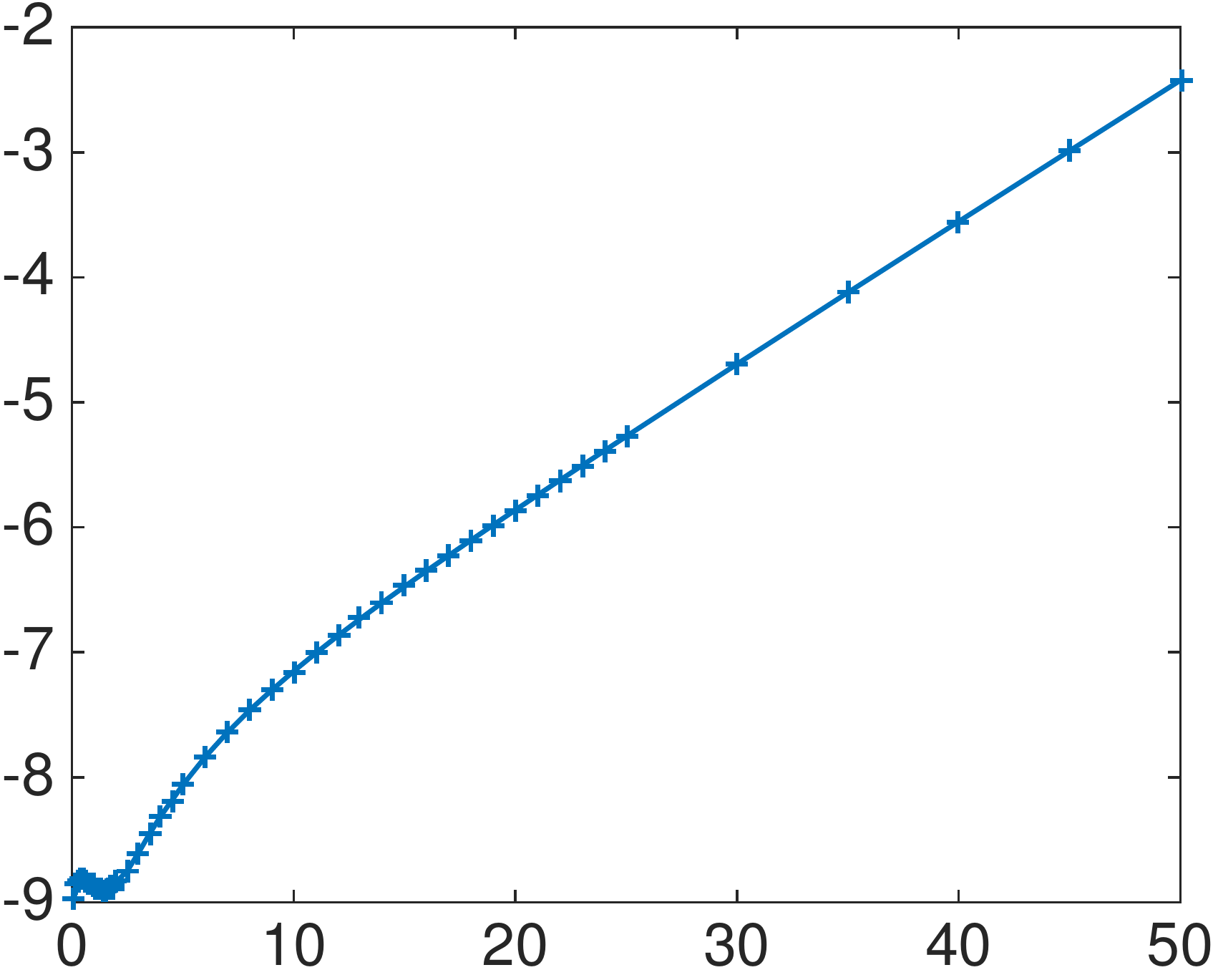}\\
$\gamma\mapsto  \log_{10}{\mathcal D}(\gamma)$ & $\gamma\mapsto  \log_{10}{\mathcal K}(\gamma)$ & $\gamma\mapsto  \log_{10}\frac{|\mathcal K_{*}(\gamma)-\mathcal K(\gamma)|}{\mathcal K_{*}(\gamma)}$ 
\end{tabular}
\caption{The diffusion and drift coefficients as functions of  $\gamma$.\label{fig.CasADK}}
\end{center}
\end{figure}
In Figure~\ref{fig.CasADK} we present the behavior of the diffusion and drift coefficients as  functions of $\gamma$. The first two columns give the result of our algorithm {\bf Step 1.-Step 5.} 
In this last column, we compare the numerical drift coefficient with the value ${\mathcal K}_{*}$ given by formula \eqref{eq.drift1D}. More precisely, we display the relative error $\frac{|\mathcal K_{*}(\gamma)-\mathcal K(\gamma)|}{\mathcal K_{*}(\gamma)}$ as a function of $\gamma$. It validates the accuracy of the algorithm. As expected from the results on the exponential decay of the second eigenvalue, the existence of a spectral gap and the formulas for $\mathcal{D}$ and $\mathcal{K}$, the diffusion and the drift coefficients grow exponentially fast as $\gamma \rightarrow \infty$.

\begin{figure}[h!t]
\begin{center}
\subfigure[$\gamma=1$]{\includegraphics[height=3cm]{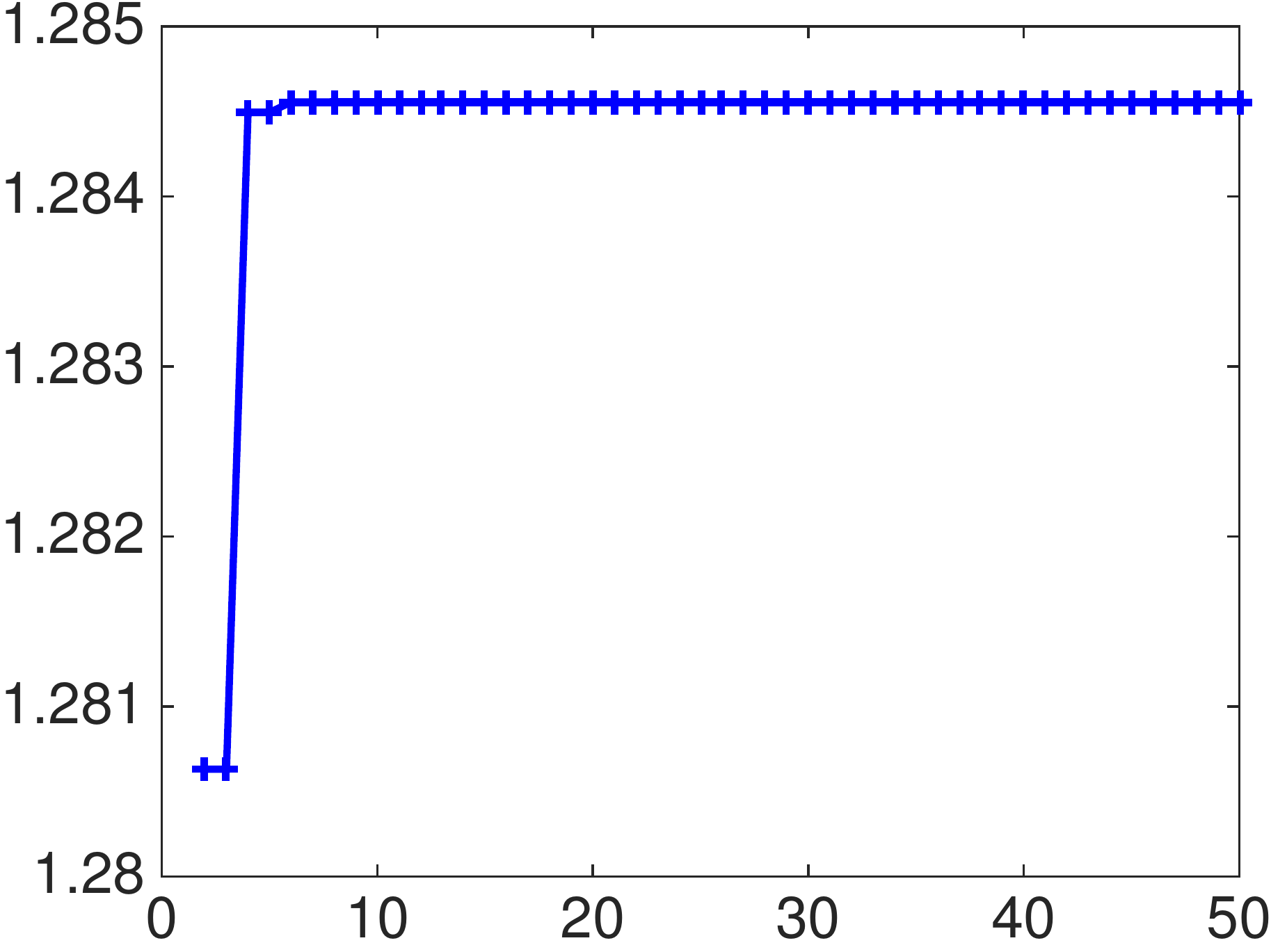}
\quad\includegraphics[height=3cm]{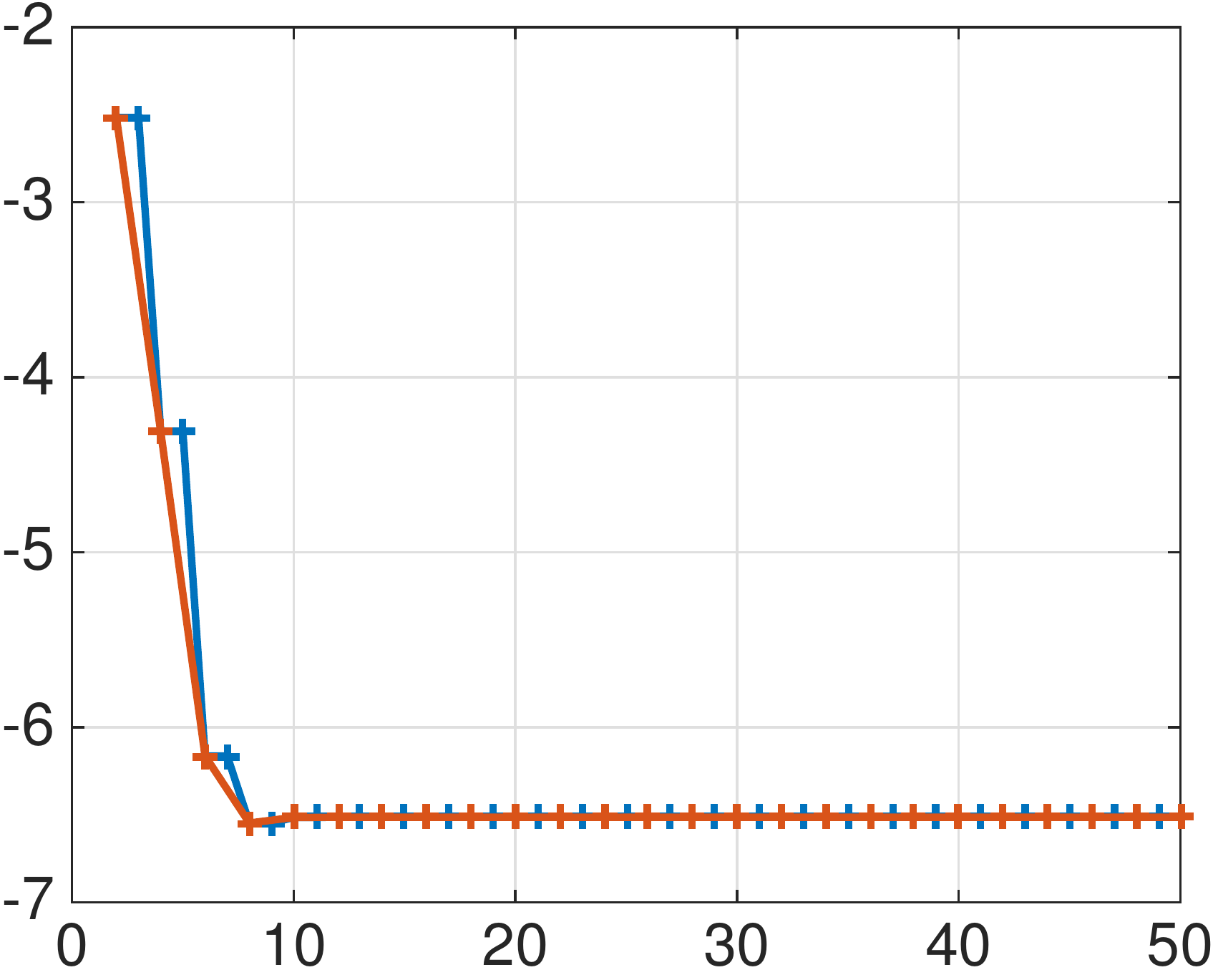}
\quad\includegraphics[height=3cm]{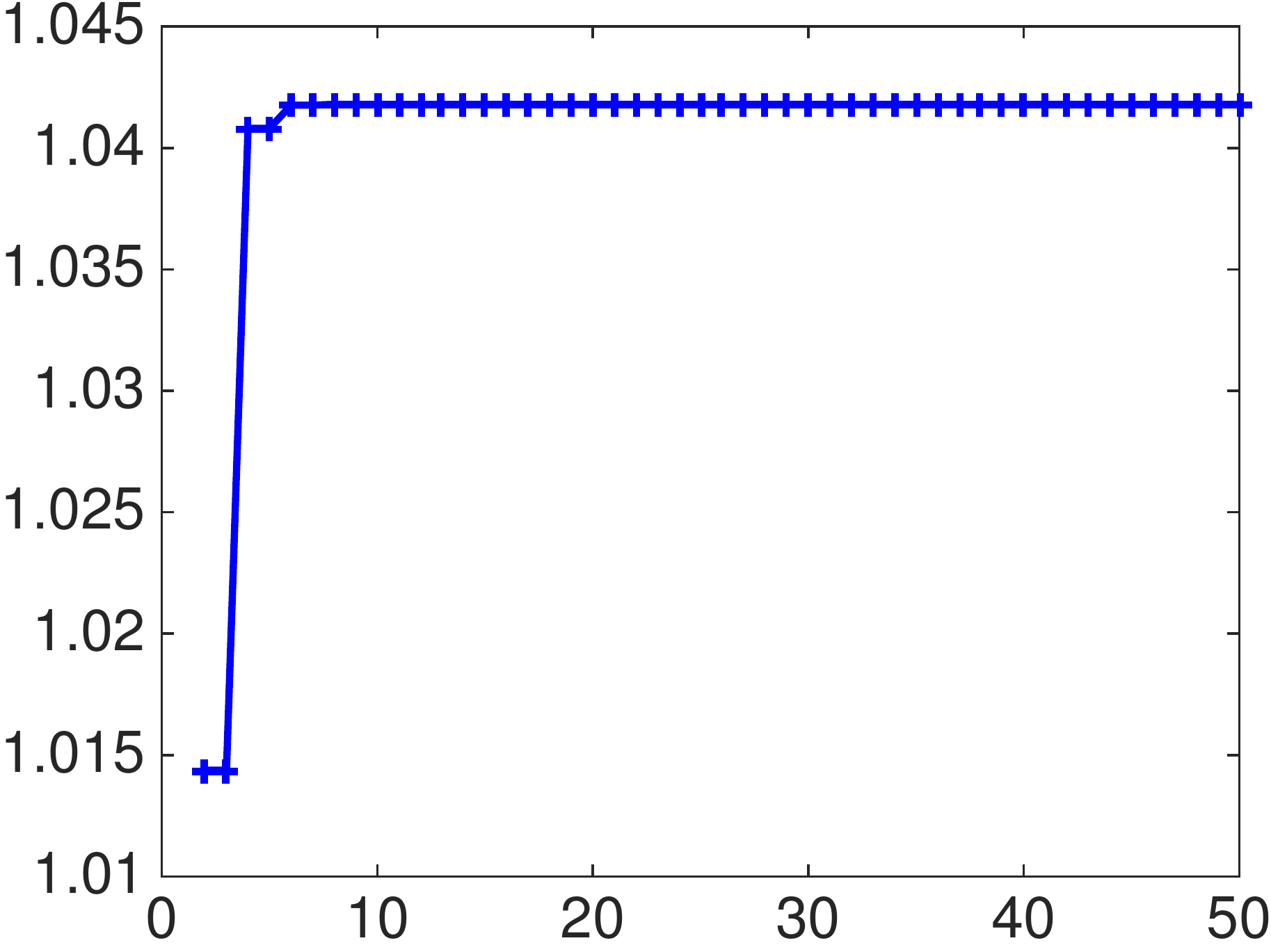}
\quad\includegraphics[height=3cm]{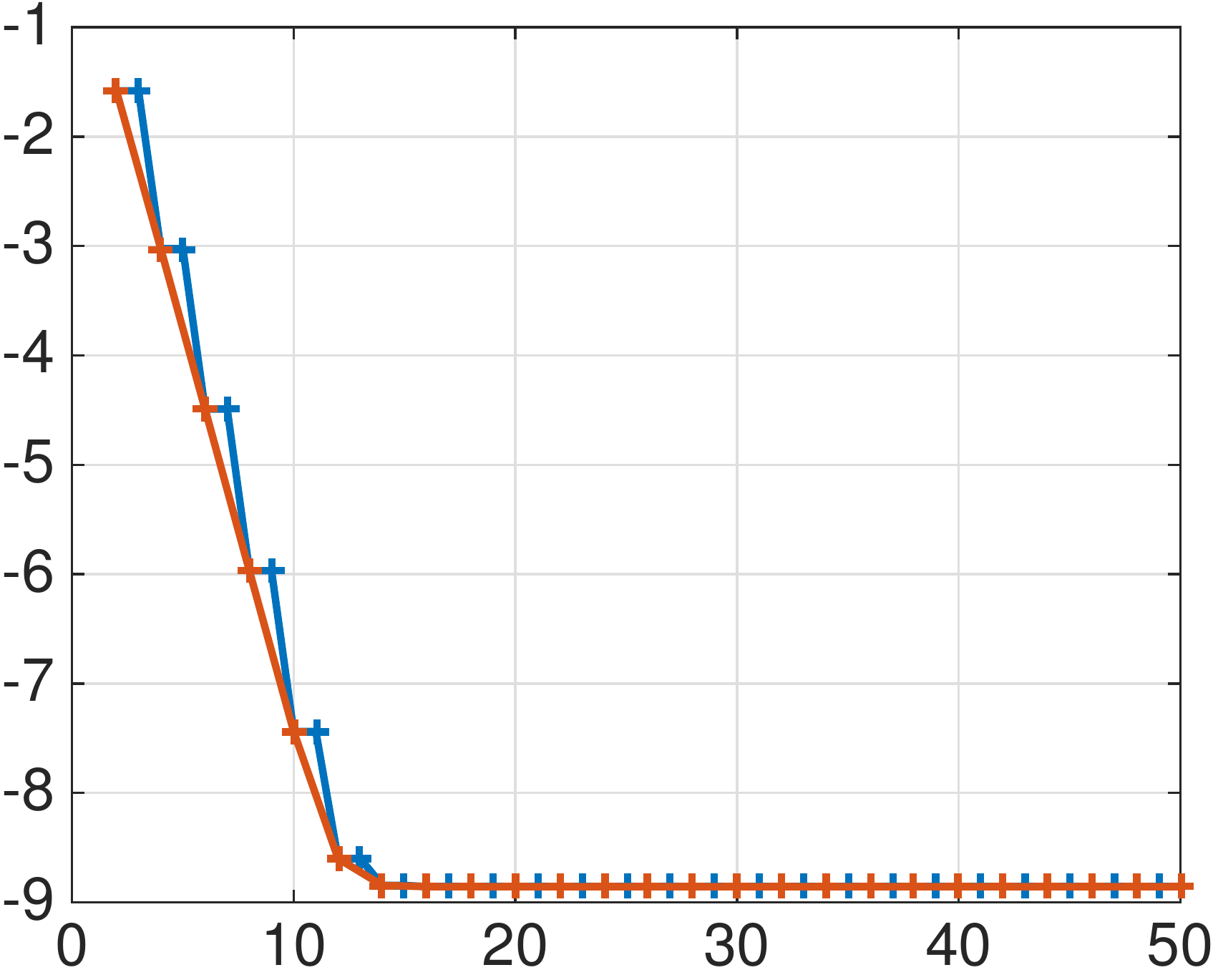}}
\subfigure[$\gamma=10$]{\includegraphics[height=3cm]{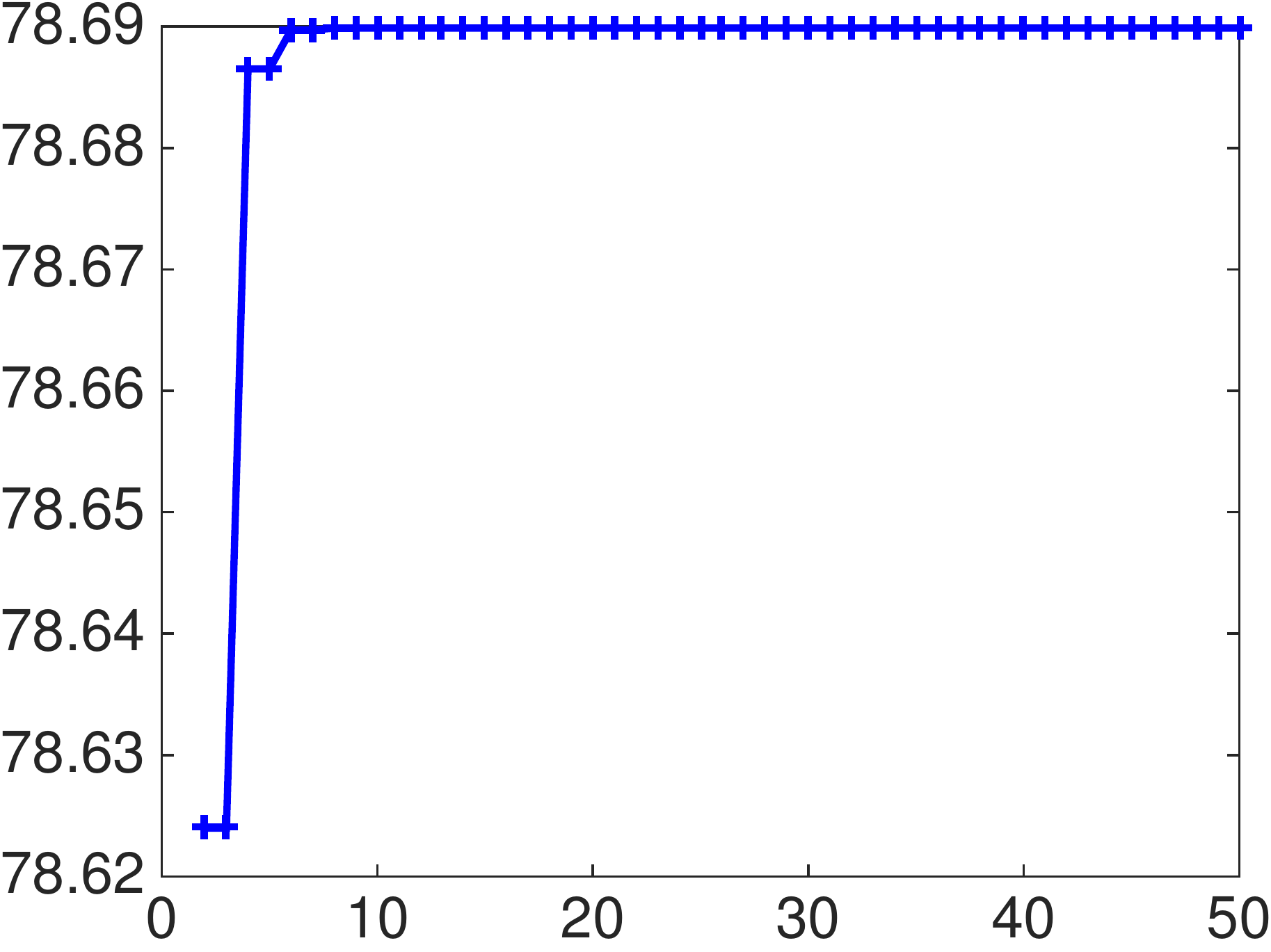}
\quad\includegraphics[height=3cm]{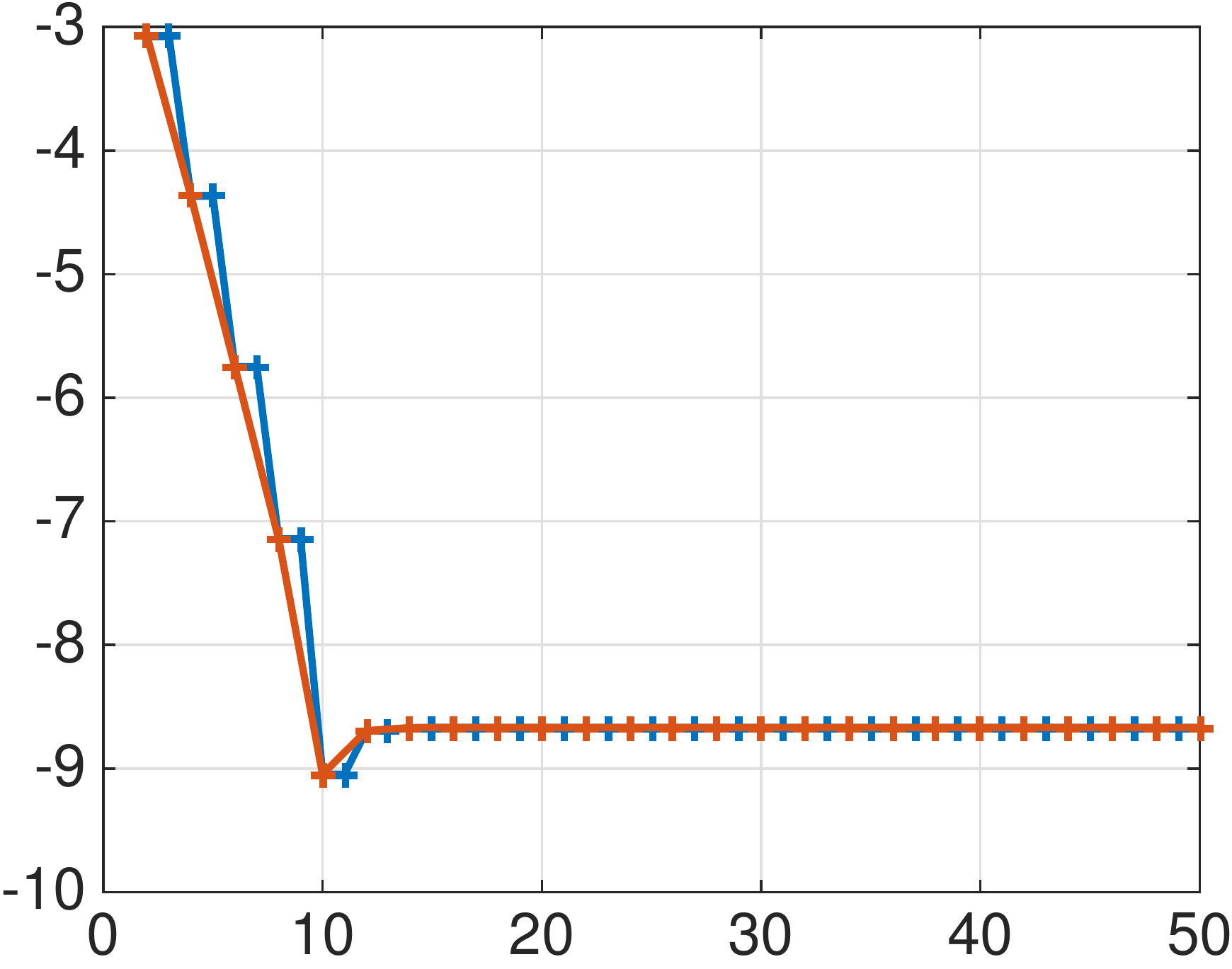}
\quad\includegraphics[height=3cm]{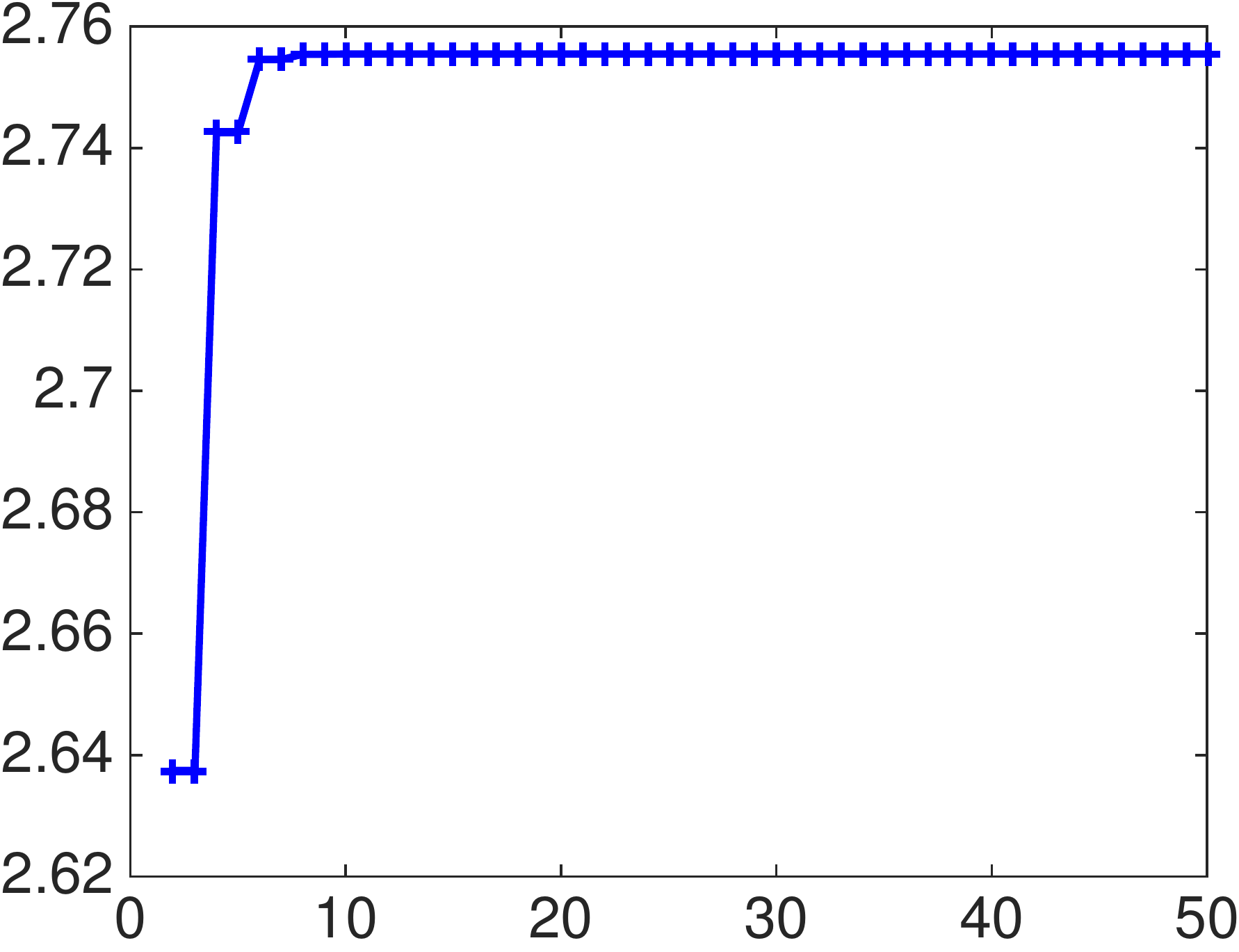}
\quad\includegraphics[height=3cm]{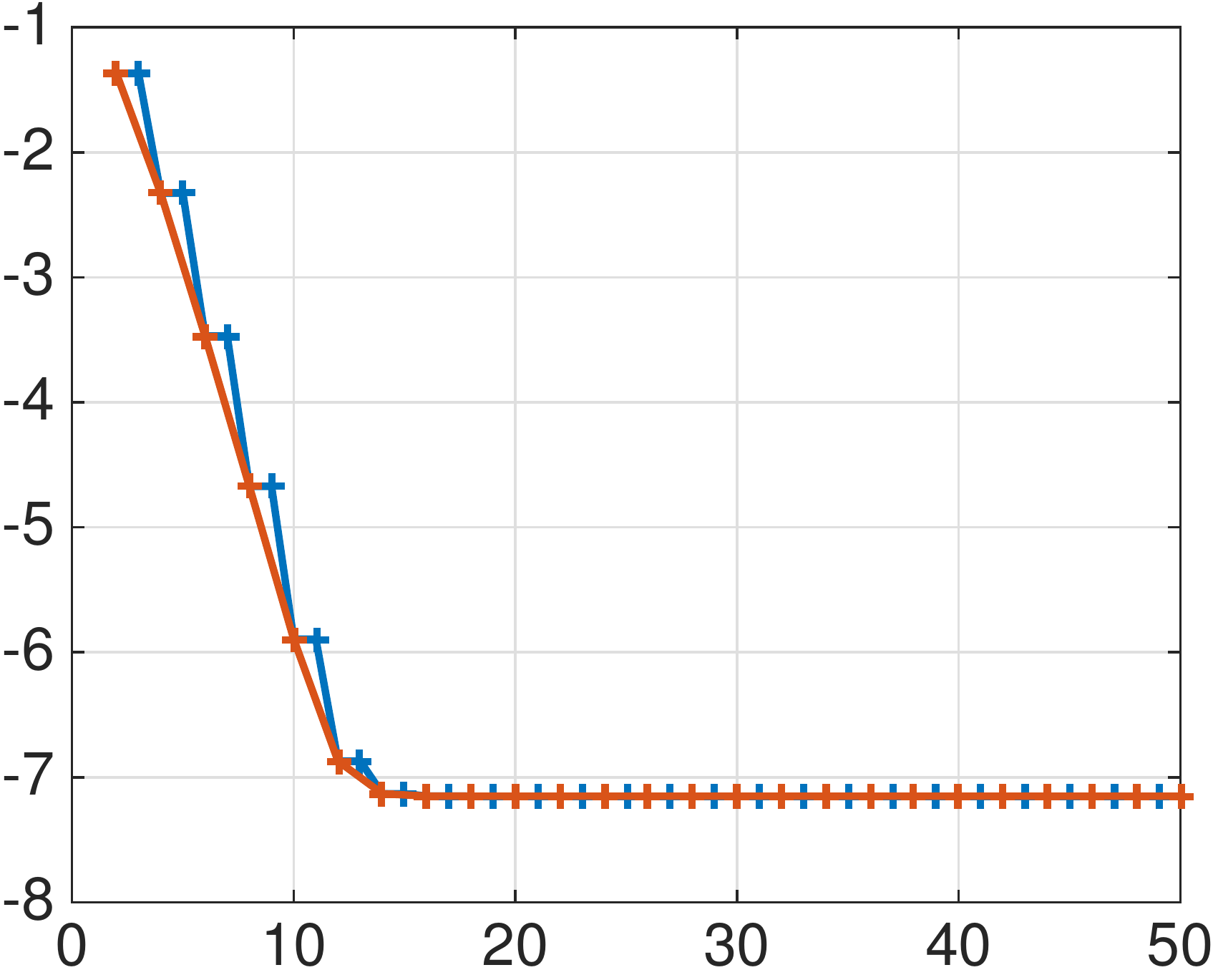}}
\subfigure[$\gamma=50$]{\includegraphics[height=3cm]{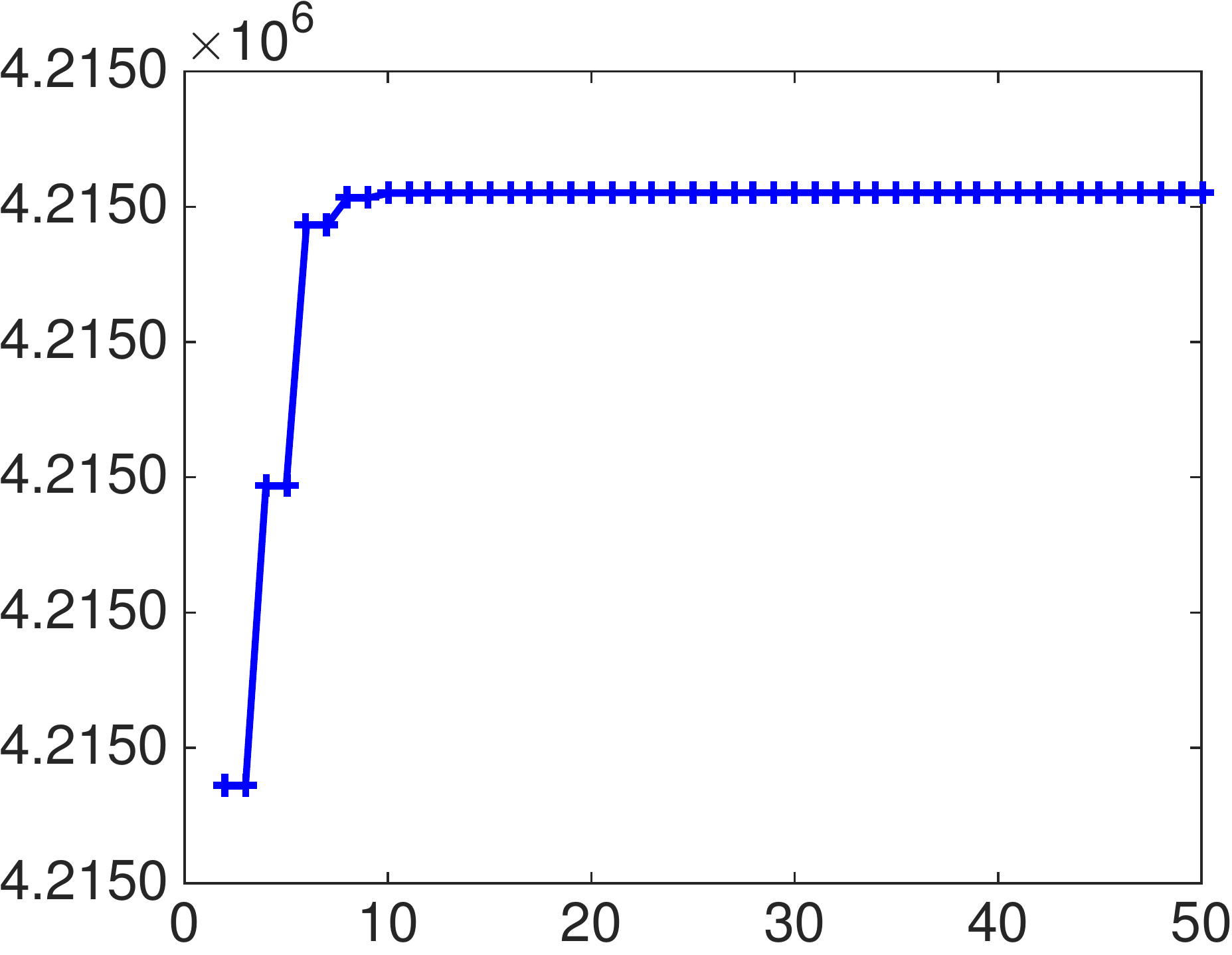}
\quad\includegraphics[height=3cm]{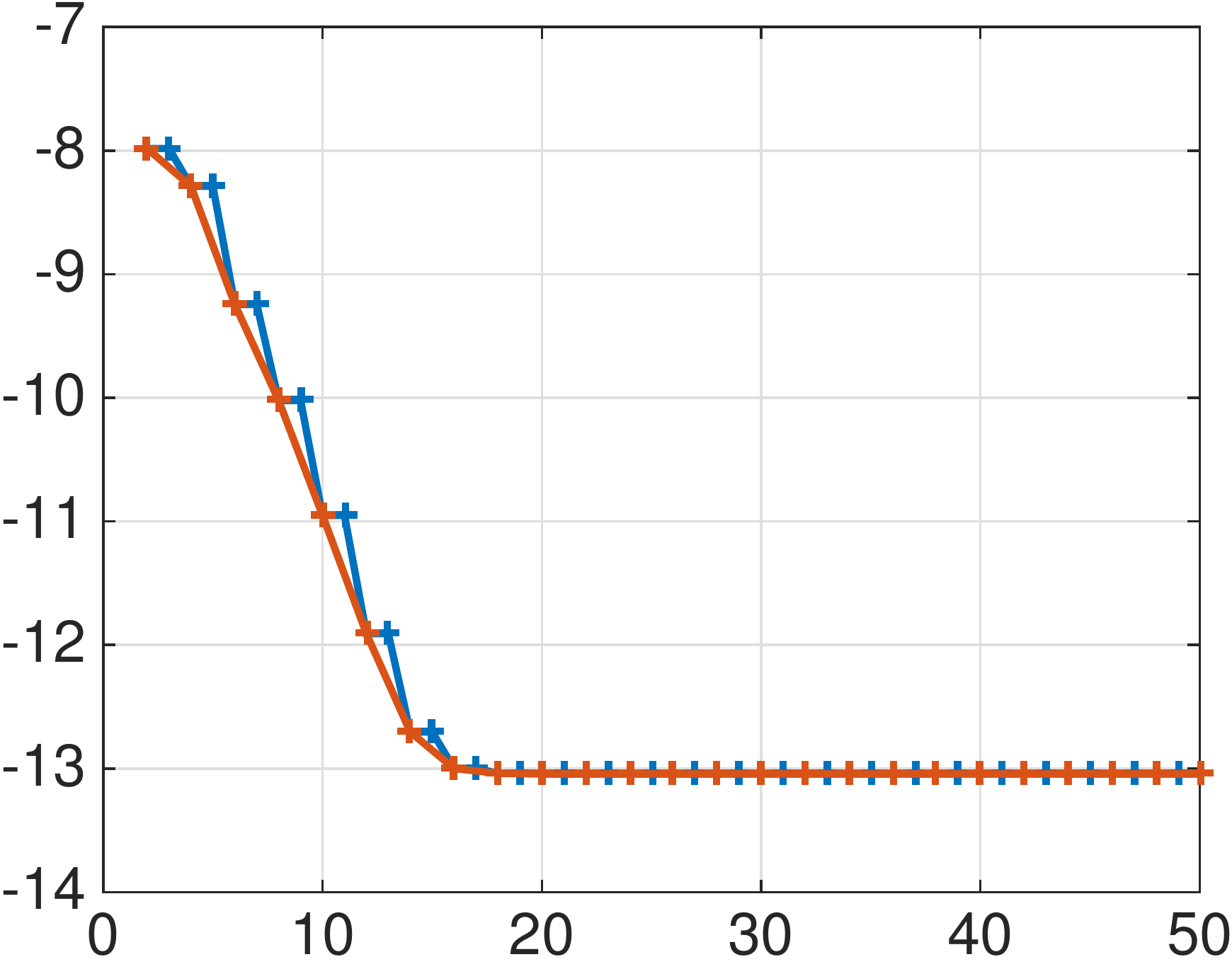}
\quad\includegraphics[height=3cm]{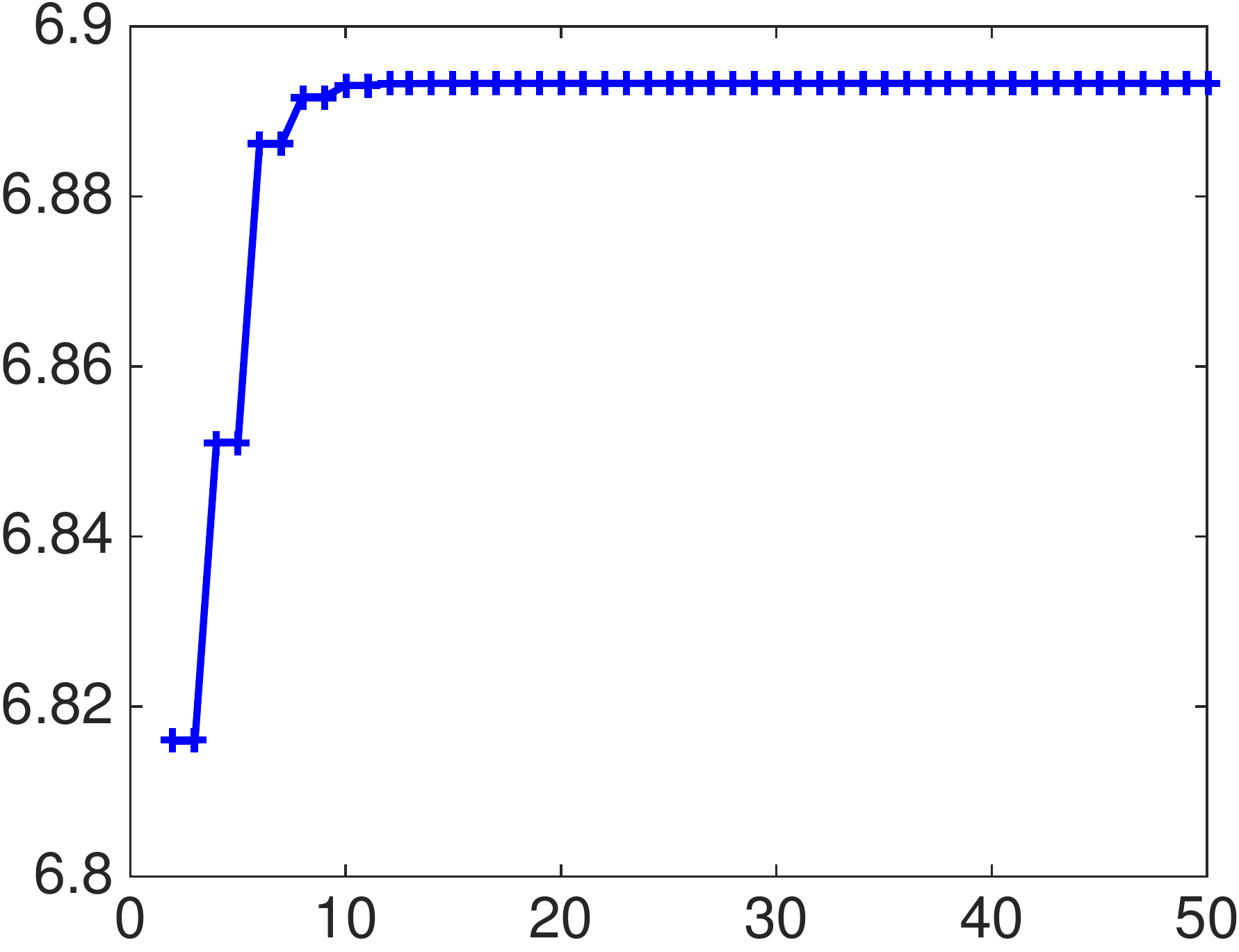}
\quad\includegraphics[height=3cm]{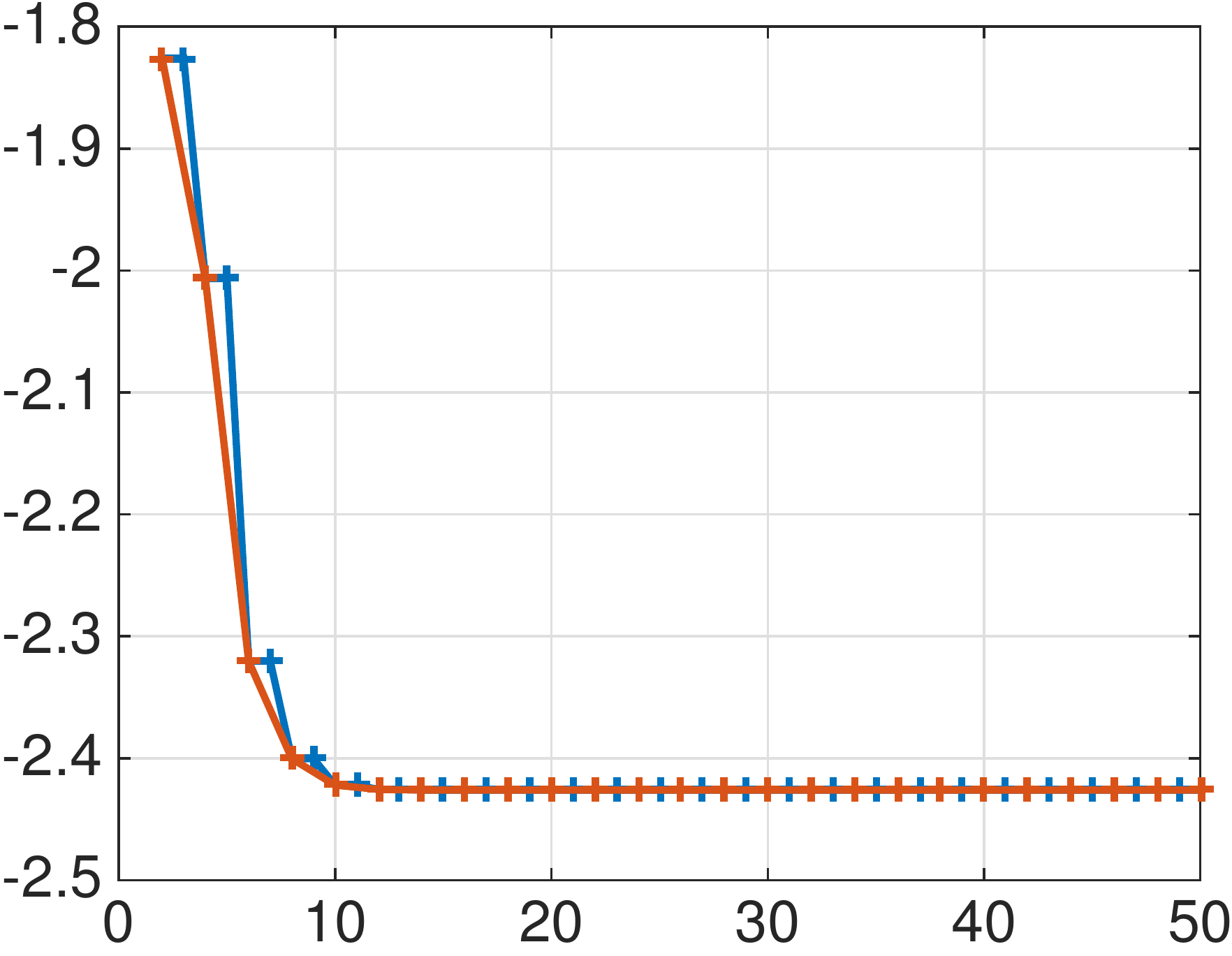}}
\caption{Convergence as a function of $N$ : $D(N)$, $\log_{10}\frac{|{\cal D}_{*}-D(N)|}{{\cal D}_{*}}$, $\kappa(N)$, $\log_{10}\frac{|{\cal K}_{*}-\kappa(N)|}{{\cal K}_{*}}$ vs. $N$. \label{fig.CasADKn}}
\end{center}
\end{figure}
Figure~\ref{fig.CasADKn} presents the convergence of the algorithm with respect to the number of eigenmodes for several values of $\gamma$ : $\gamma=1, 10,50$. Using relations~\eqref{eq.DKapprox}, we represent 
\begin{align*}
N\mapsto D(N), \qquad N\mapsto\log_{10}\frac{{\cal D}_{*}-D(N)}{{\cal D}_{*}},\qquad
N\mapsto\kappa(N),\qquad N\mapsto\log_{10}\frac{{\cal K}_{*}-\kappa(N)}{{\cal K}_{*}},
\end{align*}
where we use the shorthand notation $D(N)=D^{R,\mu,N}$ for fixed $R$ and $\mu$ (respectively $\kappa(N)=\kappa^{R,\mu,N}$)
and ${\cal D}_{*}$ is a reference value for the diffusion coefficient obtained for $N=50$ and ${\cal K}_{*}$ is as above a numerical approximation of \eqref{eq.drift1D} with a composite rectangular rule.
We observe that only 10 eigenmodes are enough to calculate accurately the diffusion coefficient and 15 eigenmodes for the drift coefficient. We can conclude that our numerical method leads to an efficient and accurate calculation of the drift and diffusion coefficients.

\subsection{Case B -- nonsmooth potential}  
\begin{figure}[h!t]
\begin{center}
\begin{tabular}{ccc}
$\gamma=1$ & $\gamma=5$ & $\gamma=7$\\
\includegraphics[height=2.5cm]{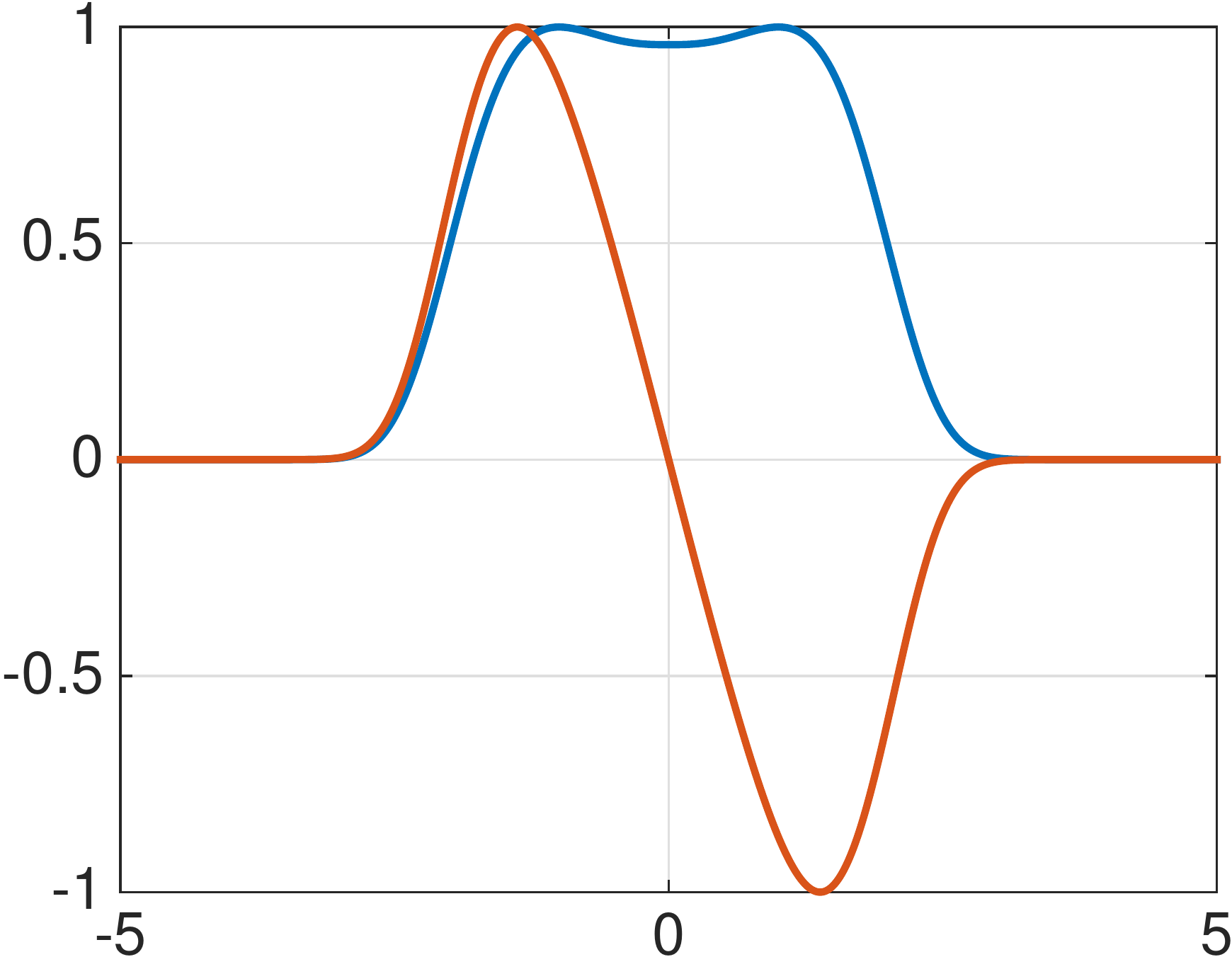}
&$\qquad$\includegraphics[height=2.5cm]{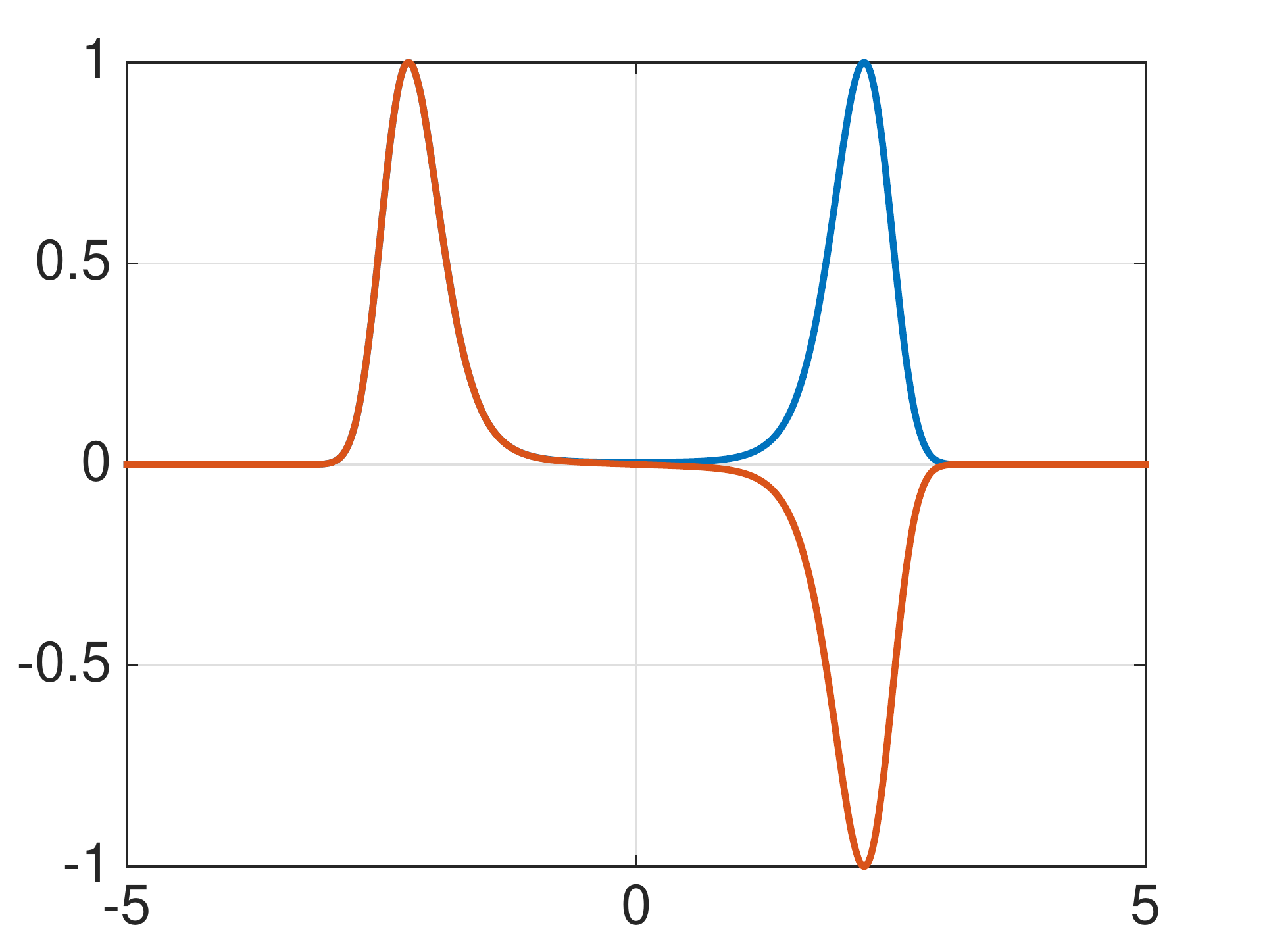}
&$\qquad$\includegraphics[height=2.5cm]{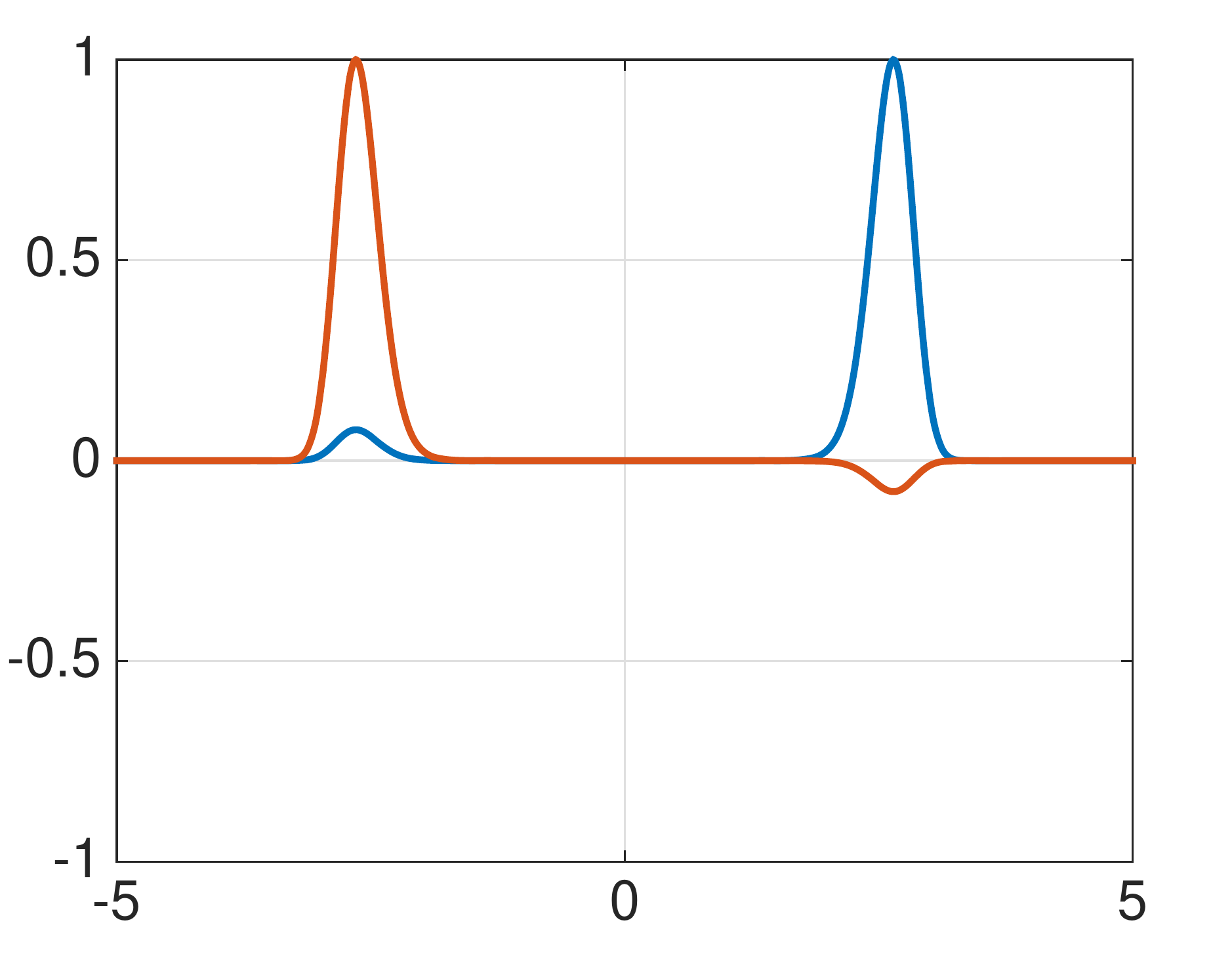}\\
\end{tabular}
\caption{First two eigenfunctions on $[-R/2,R/2]$ for $\gamma=1, 5, 7$.\label{eq.CasBVecP}}
\end{center}
\end{figure}

For the Case B, we recall that 
$$
W(v)=\tfrac{1}{4\gamma} v^4 - \tfrac{1}3 |v|^3 \qquad \mbox{with } \gamma>0\,.
$$
We show similar computations as presented in the previous subsection. Figure~\ref{eq.CasBVecP} shows that the computations are no longer accurate beyond $\gamma\geq 7$: the numerical eigenfunctions have lost their symmetry. For the nonsmooth potential, we can not take  a larger value of $\gamma$. 
\begin{figure}[h!t]
\begin{center}
\begin{tabular}{cc}
\includegraphics[height=3cm]{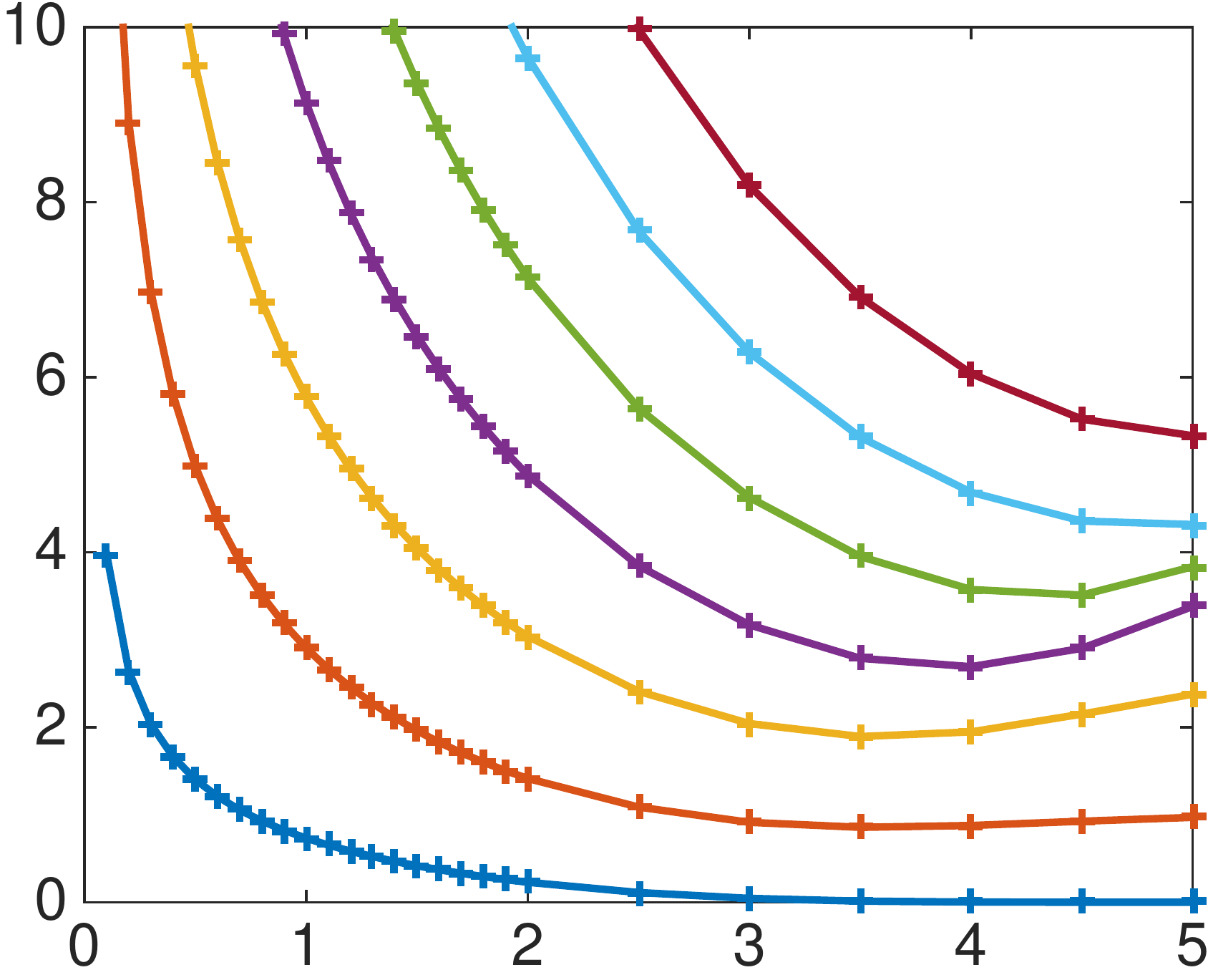}
&$\qquad$\includegraphics[height=3cm]{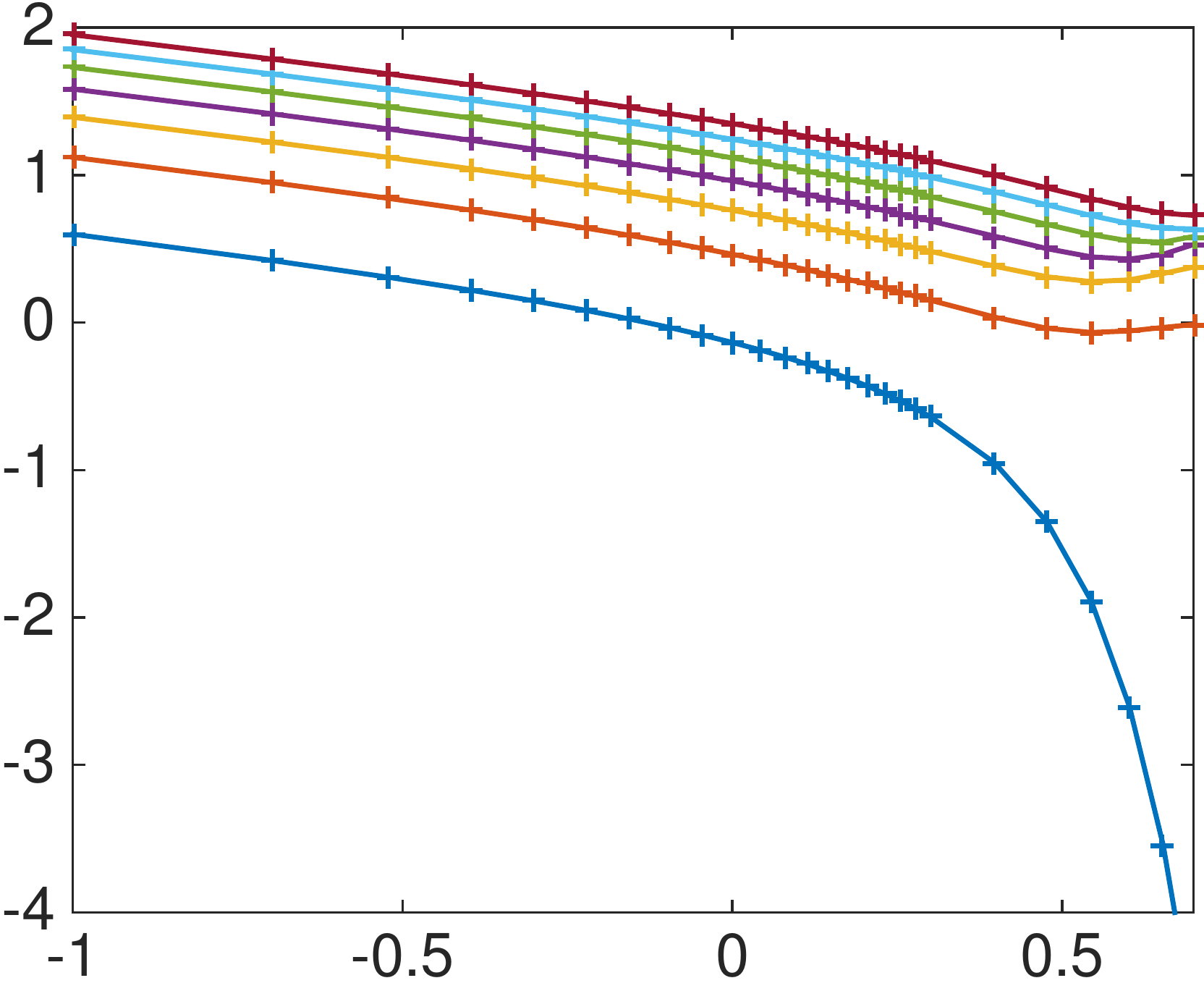}\\
$\gamma\mapsto\lambda_{j}(\gamma)$ & $\qquad\log\gamma\mapsto\log_{10}\lambda_{j}(\gamma)$ 
\end{tabular}
\caption{Convergence of the eigenvalues.\label{eq.CasBVP}}
\end{center}
\end{figure}

\begin{figure}[h!t]
\begin{center}
\begin{tabular}{cc|c}
\includegraphics[height=3.5cm]{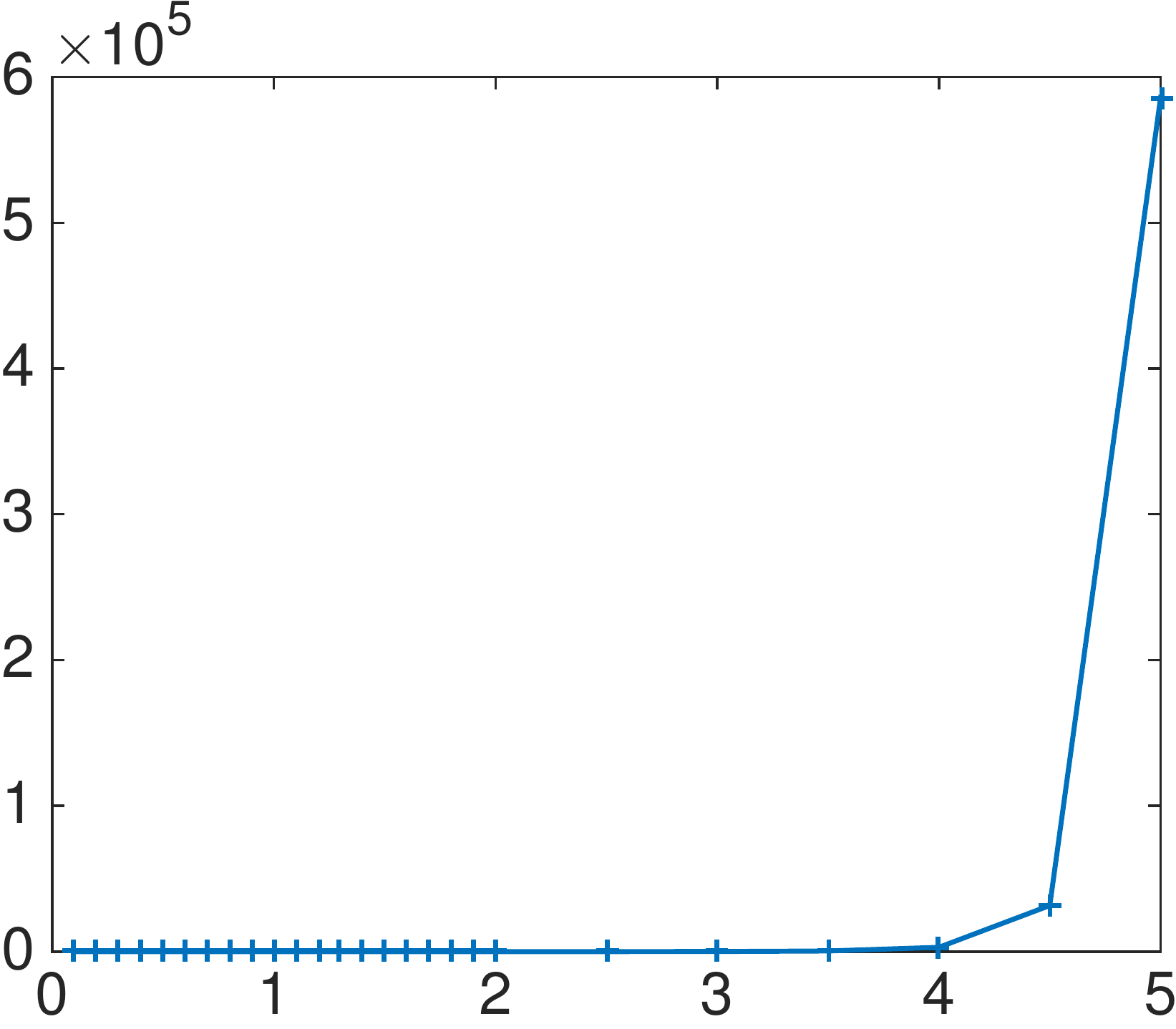}
\qquad&\includegraphics[height=3.5cm]{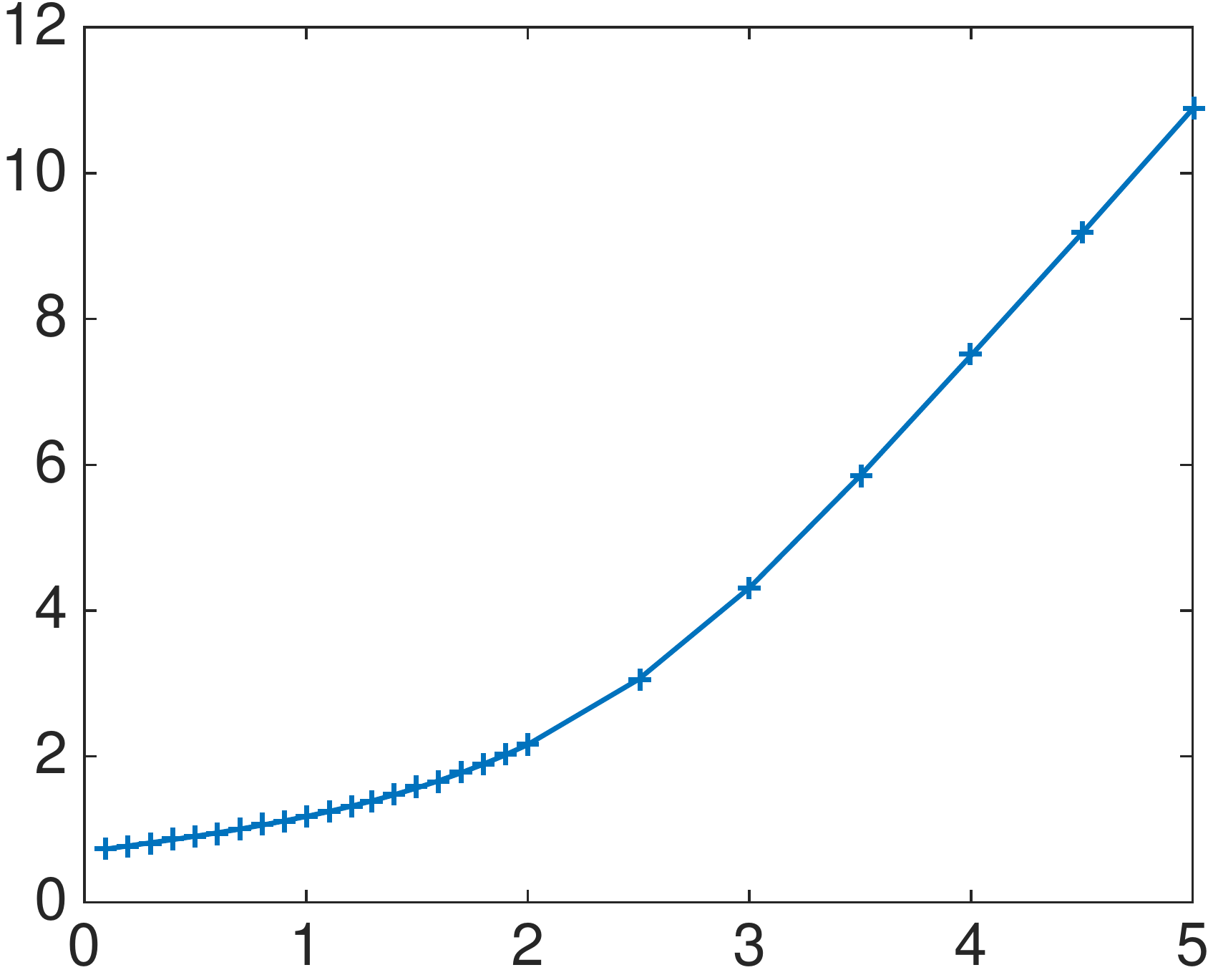}
\qquad&\quad\includegraphics[height=3.5cm]{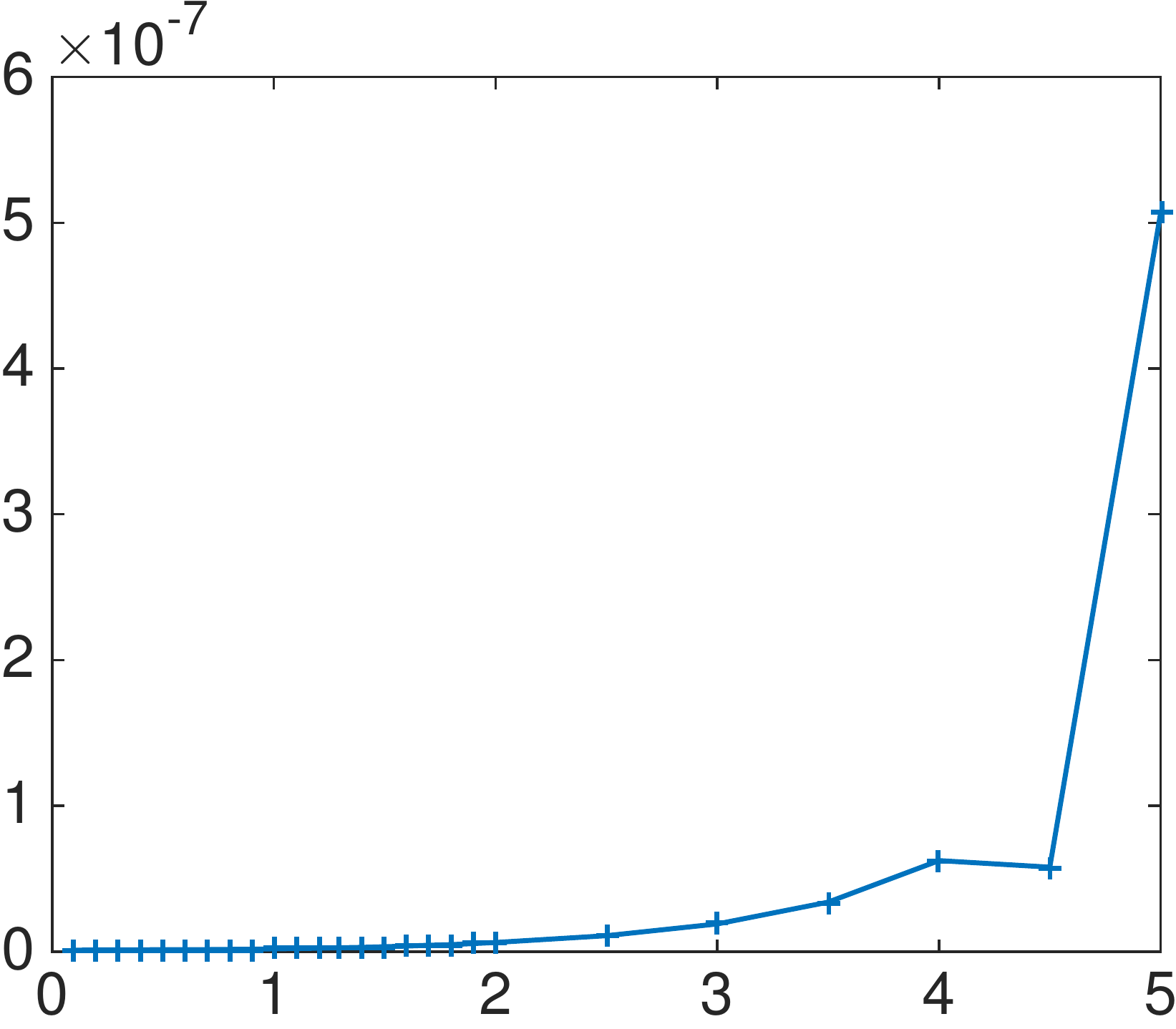}\\
$\gamma\mapsto{\mathcal D}(\gamma)$ & $\gamma\mapsto{\mathcal K}(\gamma)$  &$\gamma\mapsto \frac{|\mathcal K_{*}(\gamma)-\mathcal K(\gamma)|}{\mathcal K_{*}(\gamma)}$ \\[5pt]
\includegraphics[height=3.5cm]{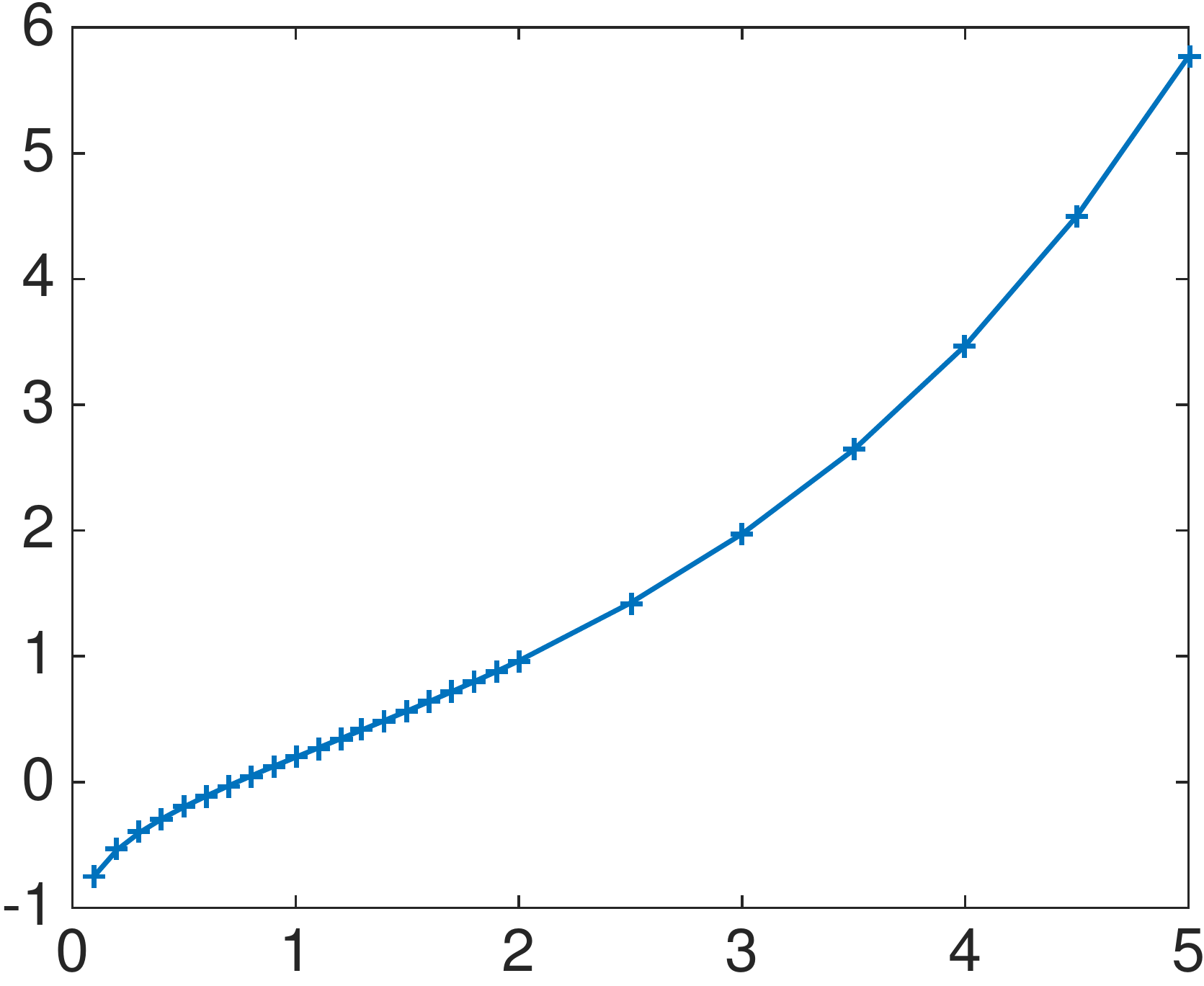}
\qquad&\includegraphics[height=3.5cm]{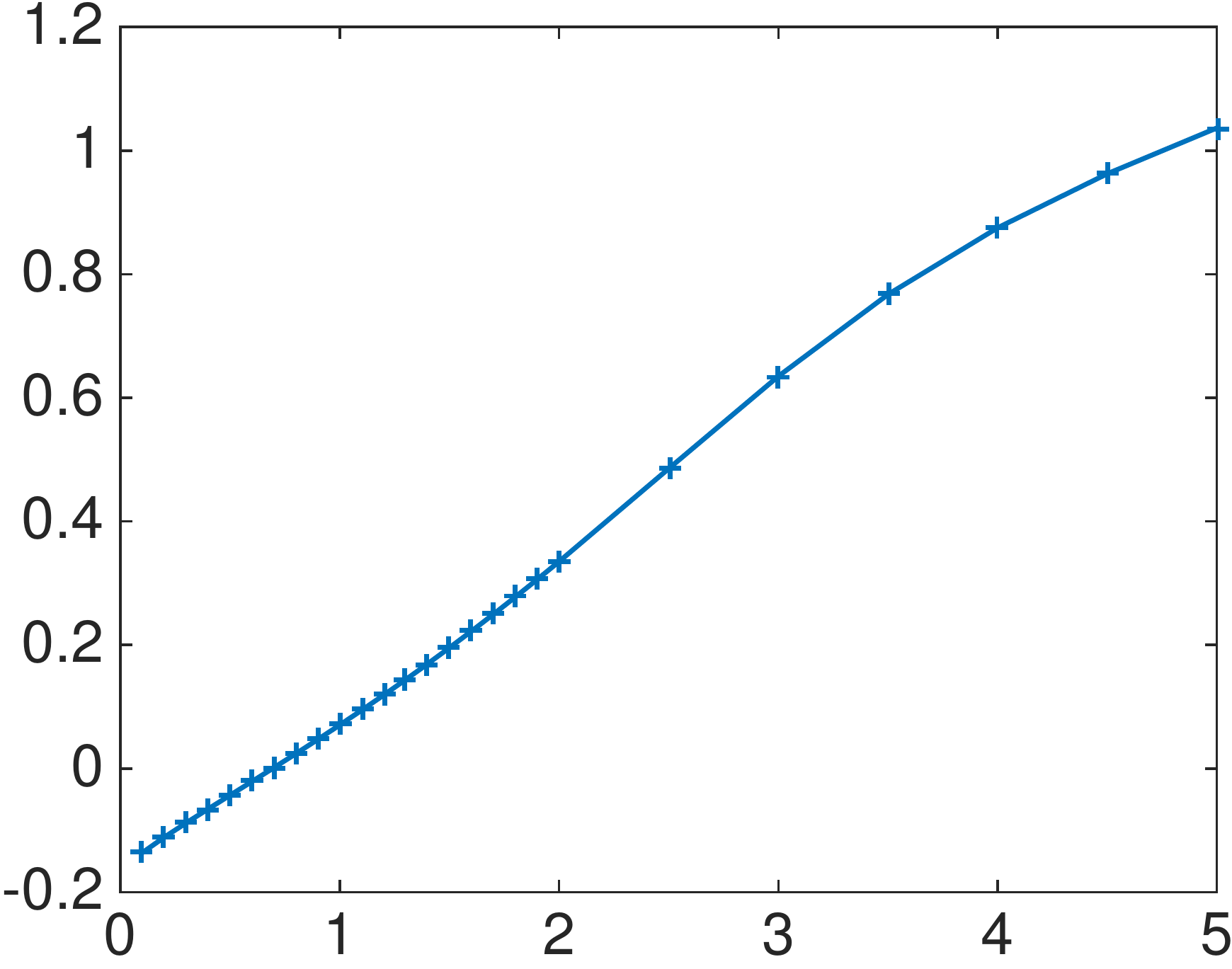}
\qquad&\quad\includegraphics[height=3.5cm]{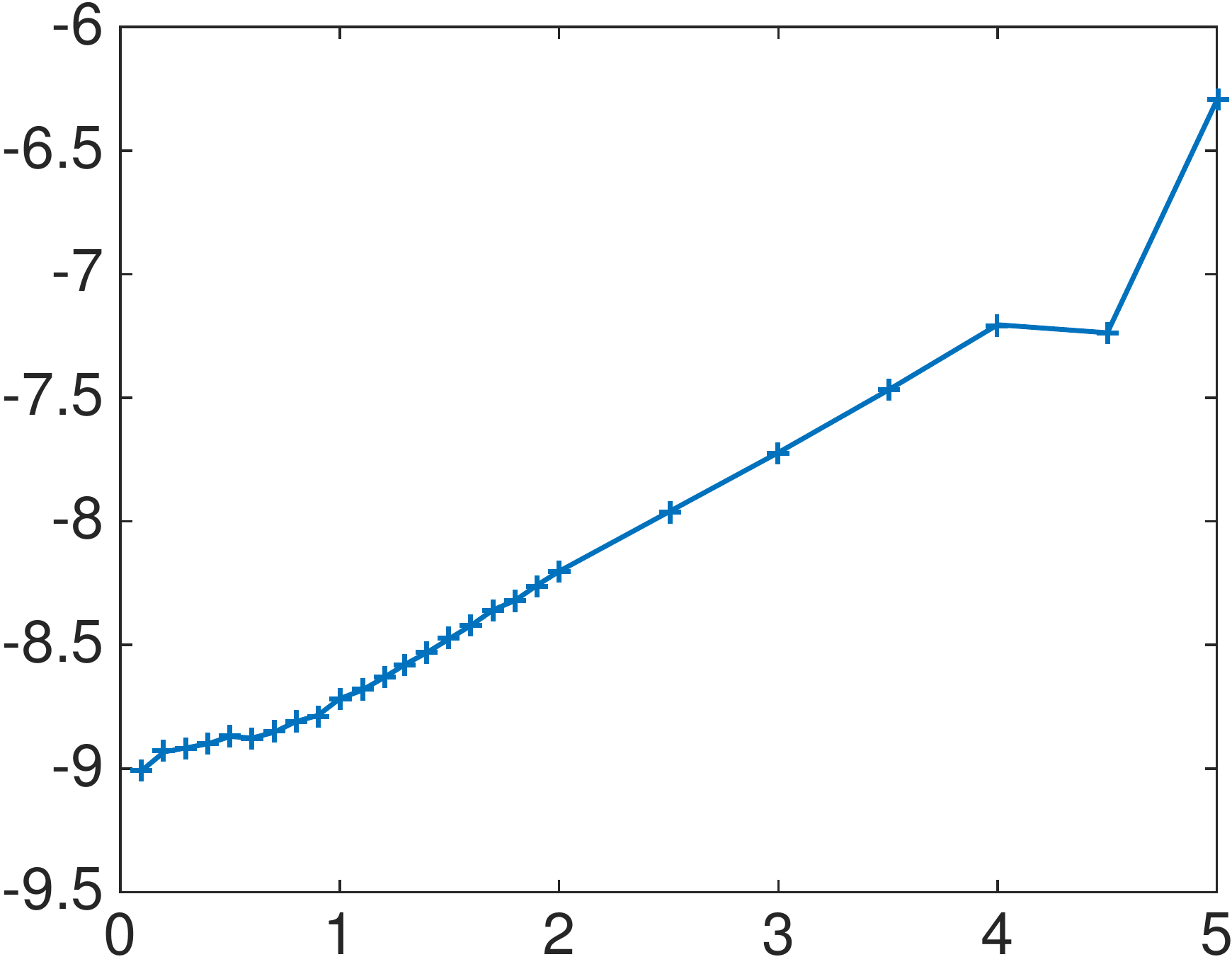}\\
$\gamma\mapsto  \log_{10}{\mathcal D}(\gamma)$ & $\gamma\mapsto \log_{10}{\mathcal K}(\gamma)$ & $\gamma\mapsto  \log_{10}\frac{|\mathcal K_{*}(\gamma)-\mathcal K(\gamma)|}{\mathcal K_{*}(\gamma)}$ 
\end{tabular}
\caption{Behavior of the diffusion and drift coefficients as a function of  $\gamma$.\label{fig.CasBDK}}
\end{center}
\end{figure}
\begin{figure}[h!]
\begin{center}
\subfigure[$\gamma=1$]{\includegraphics[height=3cm]{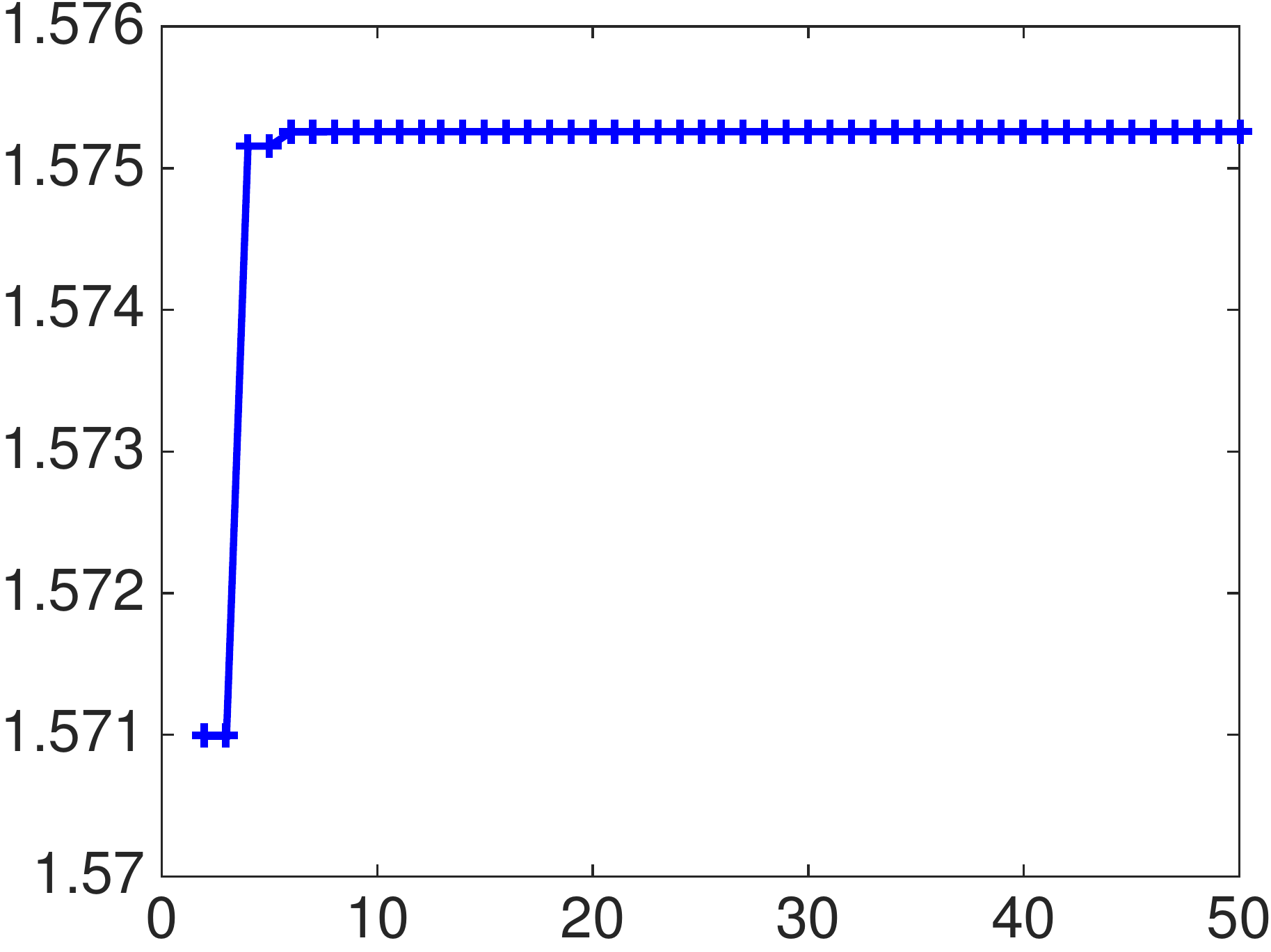}
\quad\includegraphics[height=3cm]{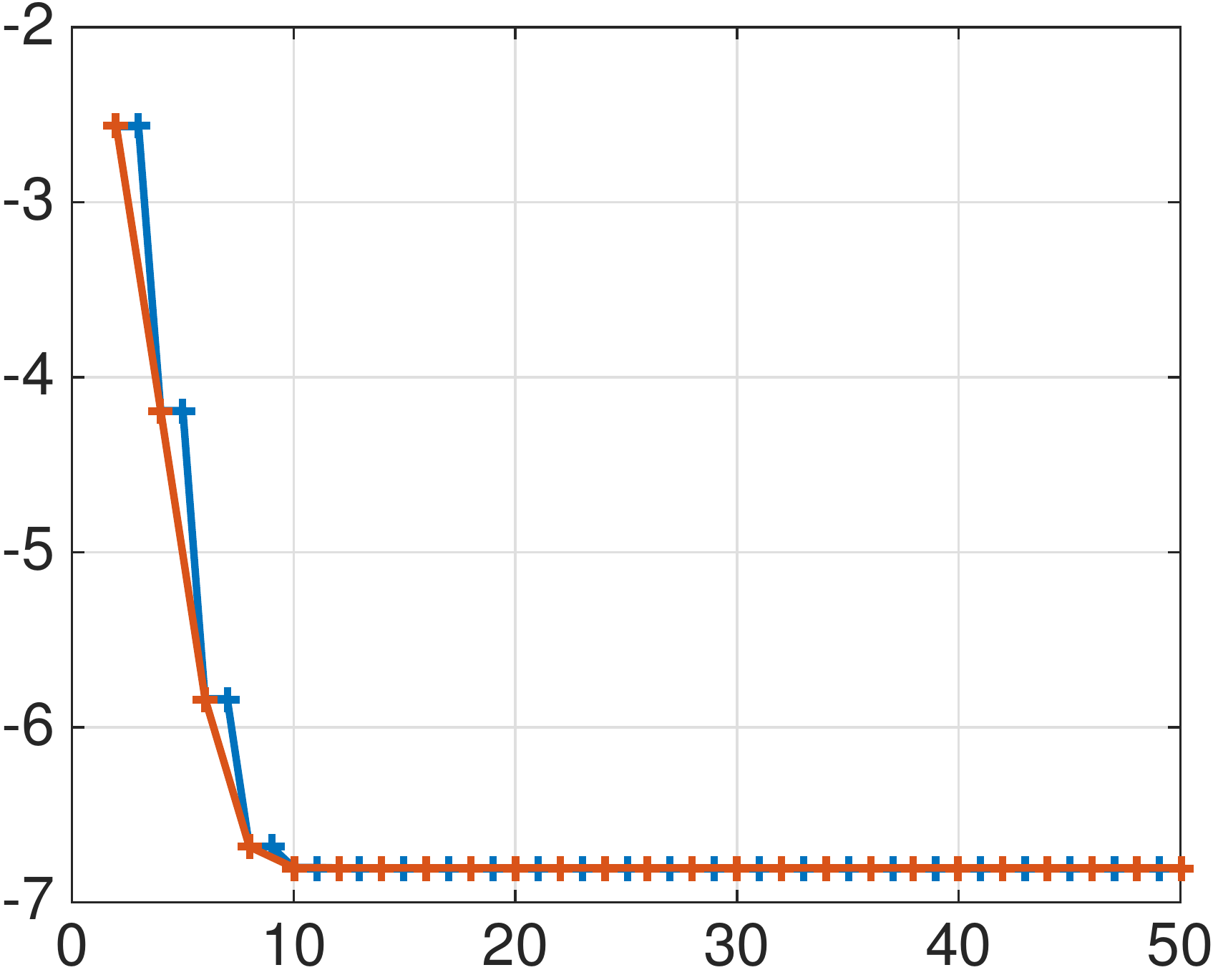}
\quad\includegraphics[height=3cm]{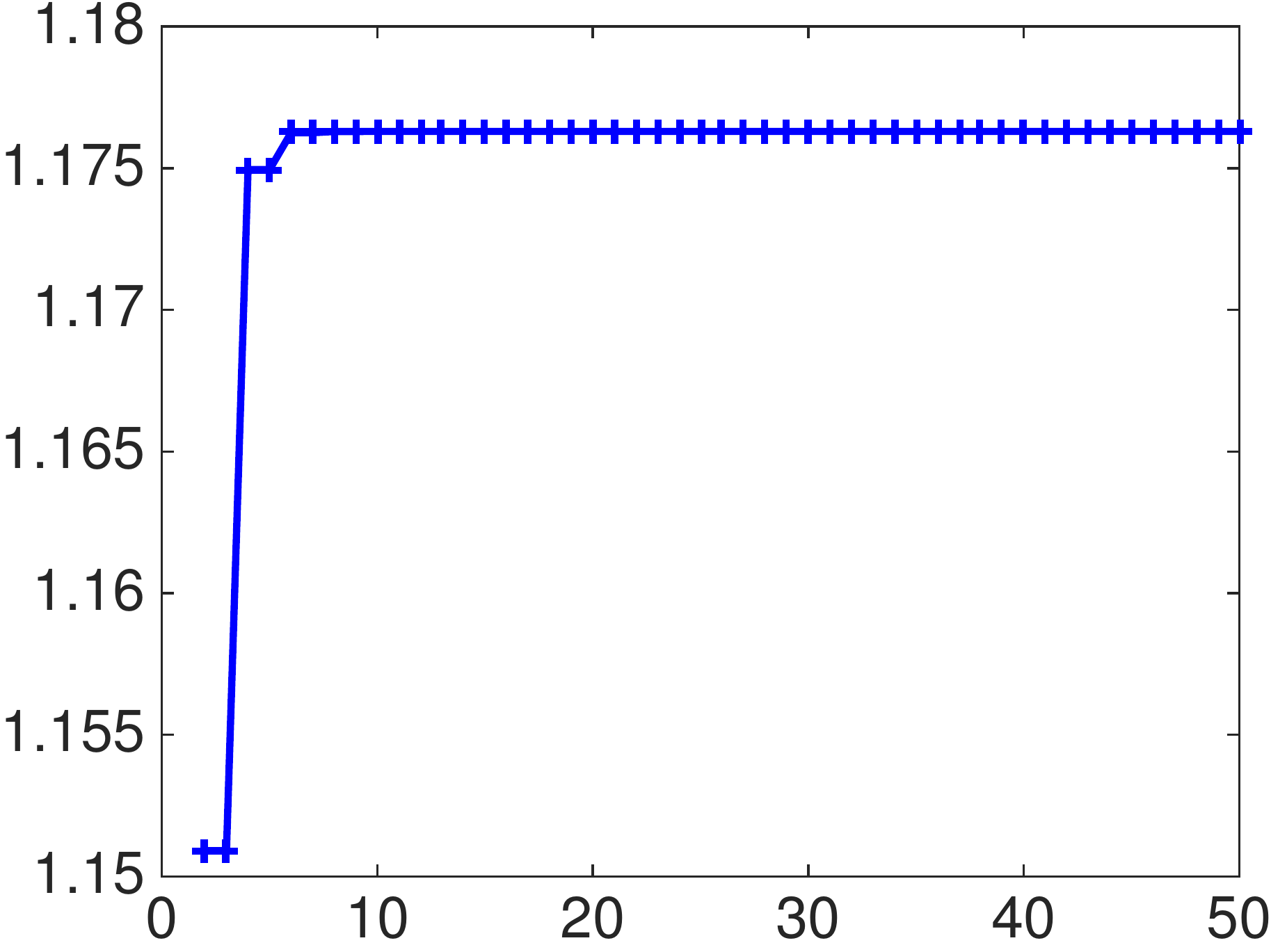}
\quad\includegraphics[height=3cm]{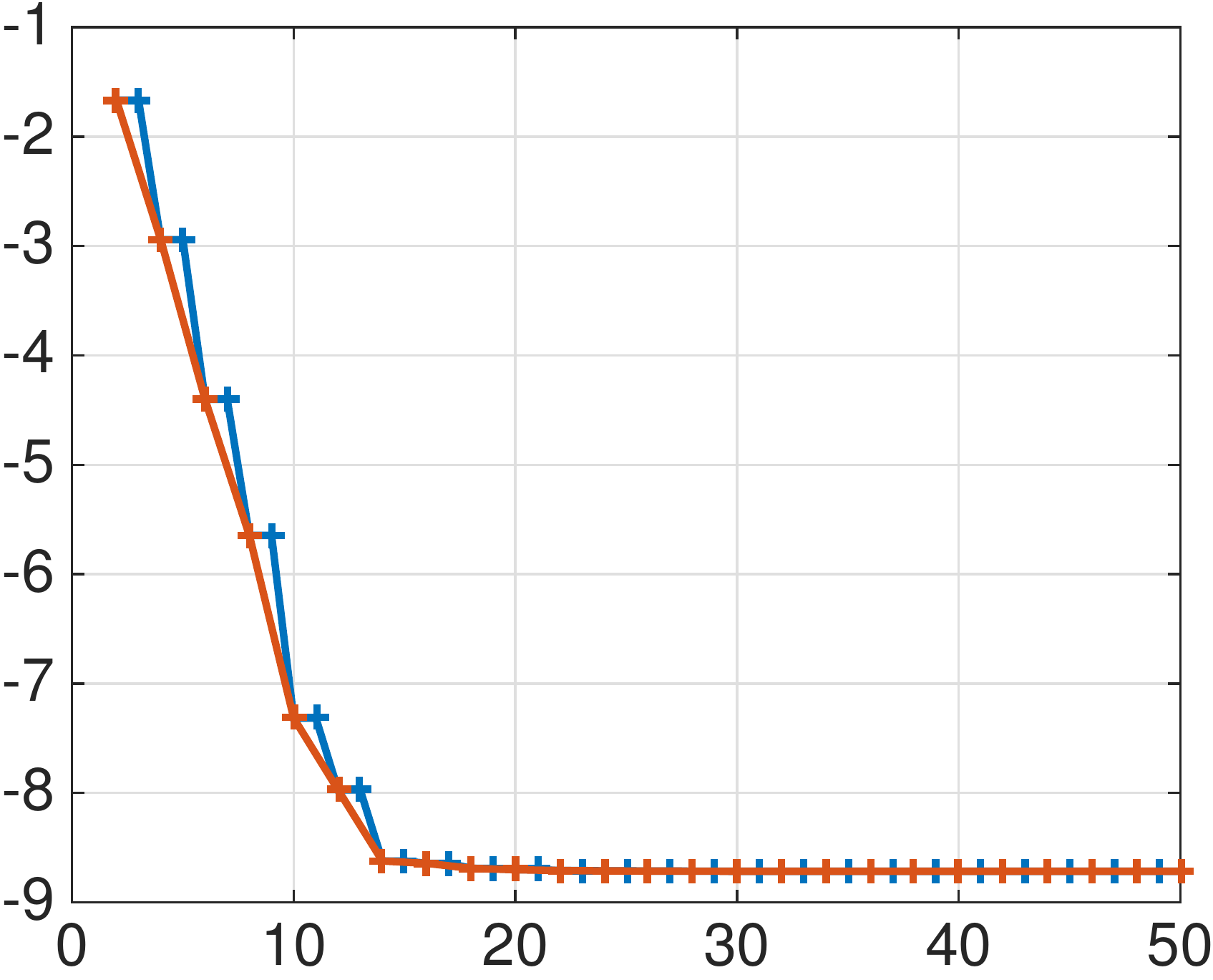}}
\subfigure[$\gamma=5$]{\includegraphics[height=3cm]{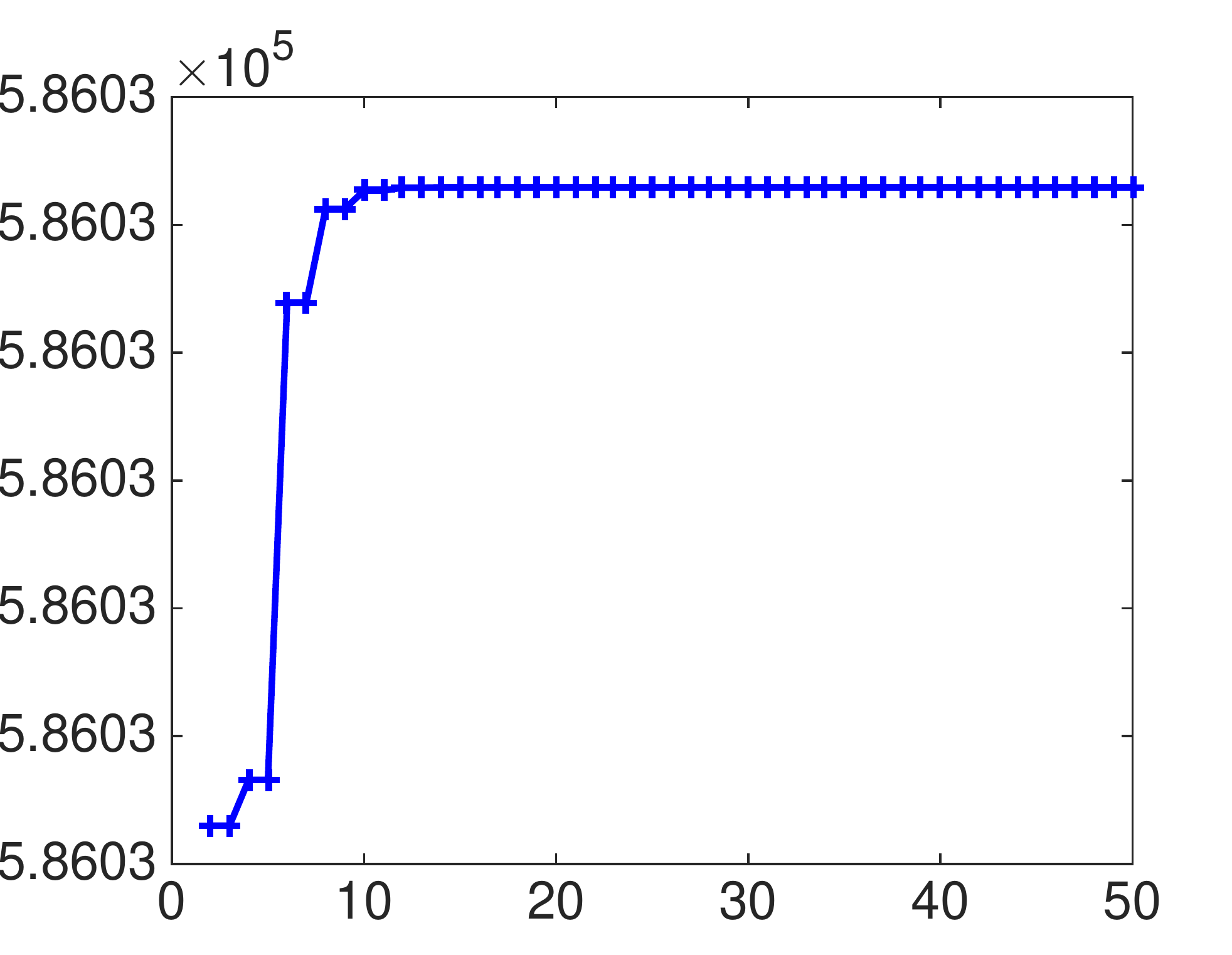}
\quad\includegraphics[height=3cm]{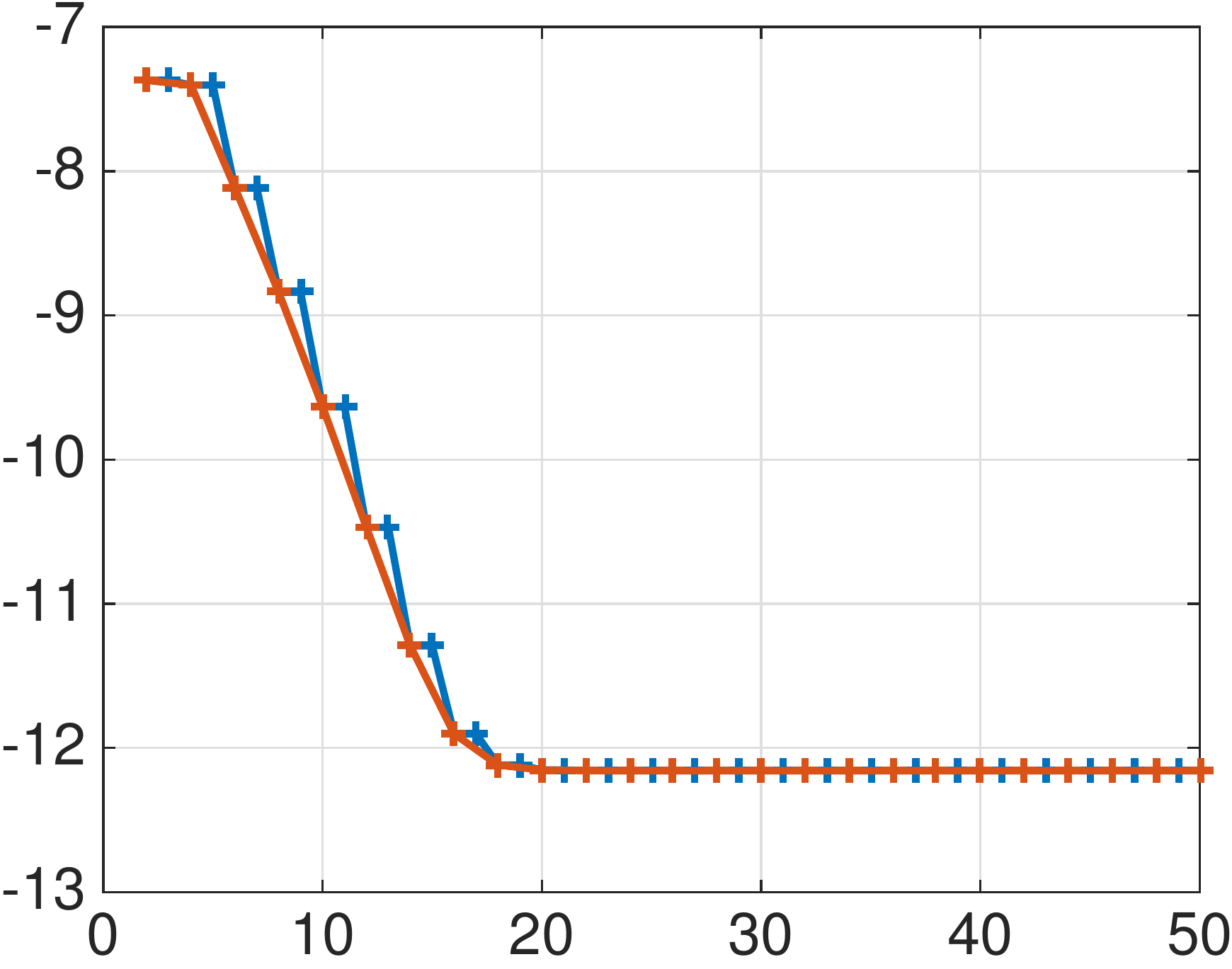}
\quad\includegraphics[height=3cm]{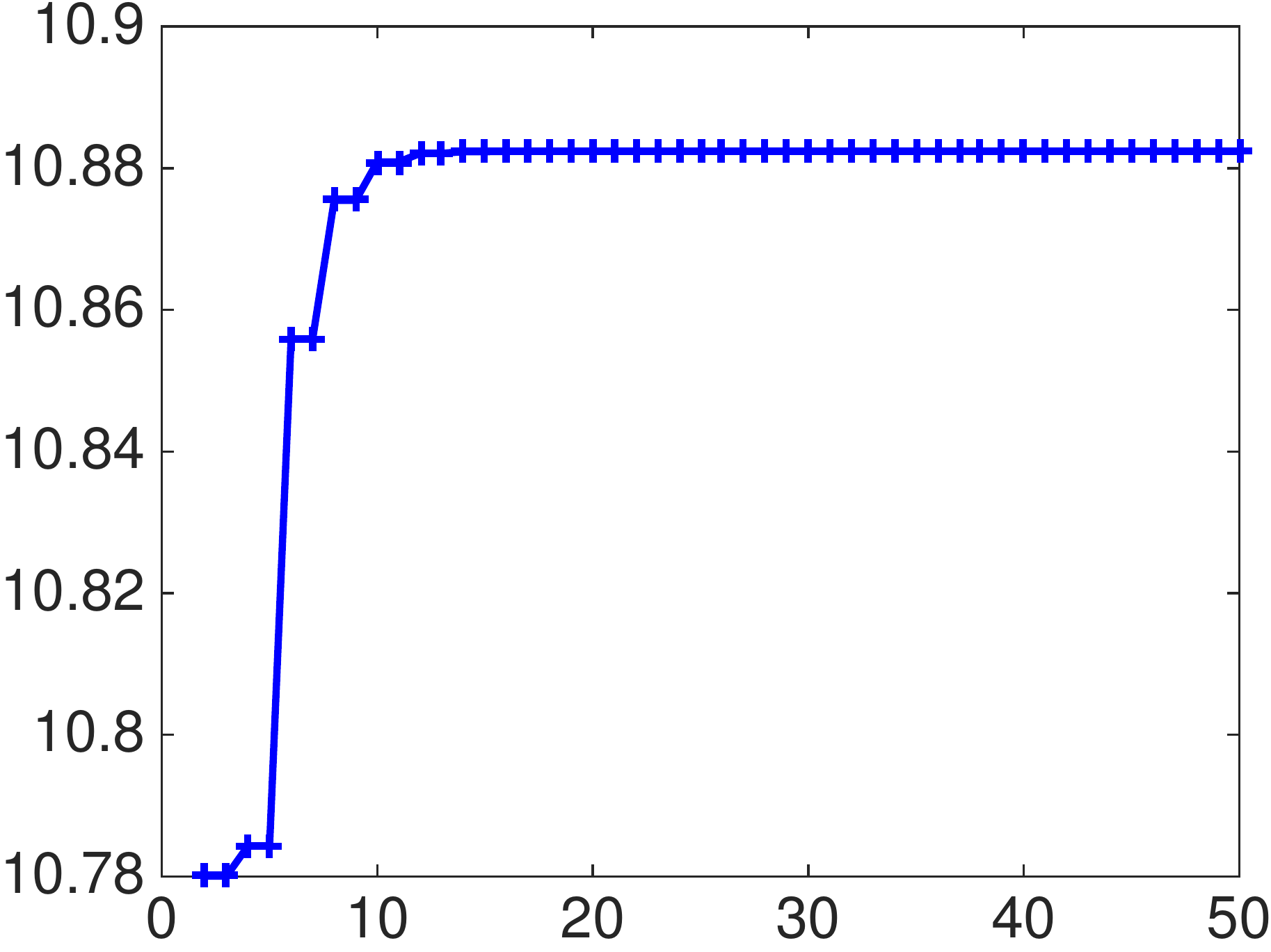}
\quad\includegraphics[height=3cm]{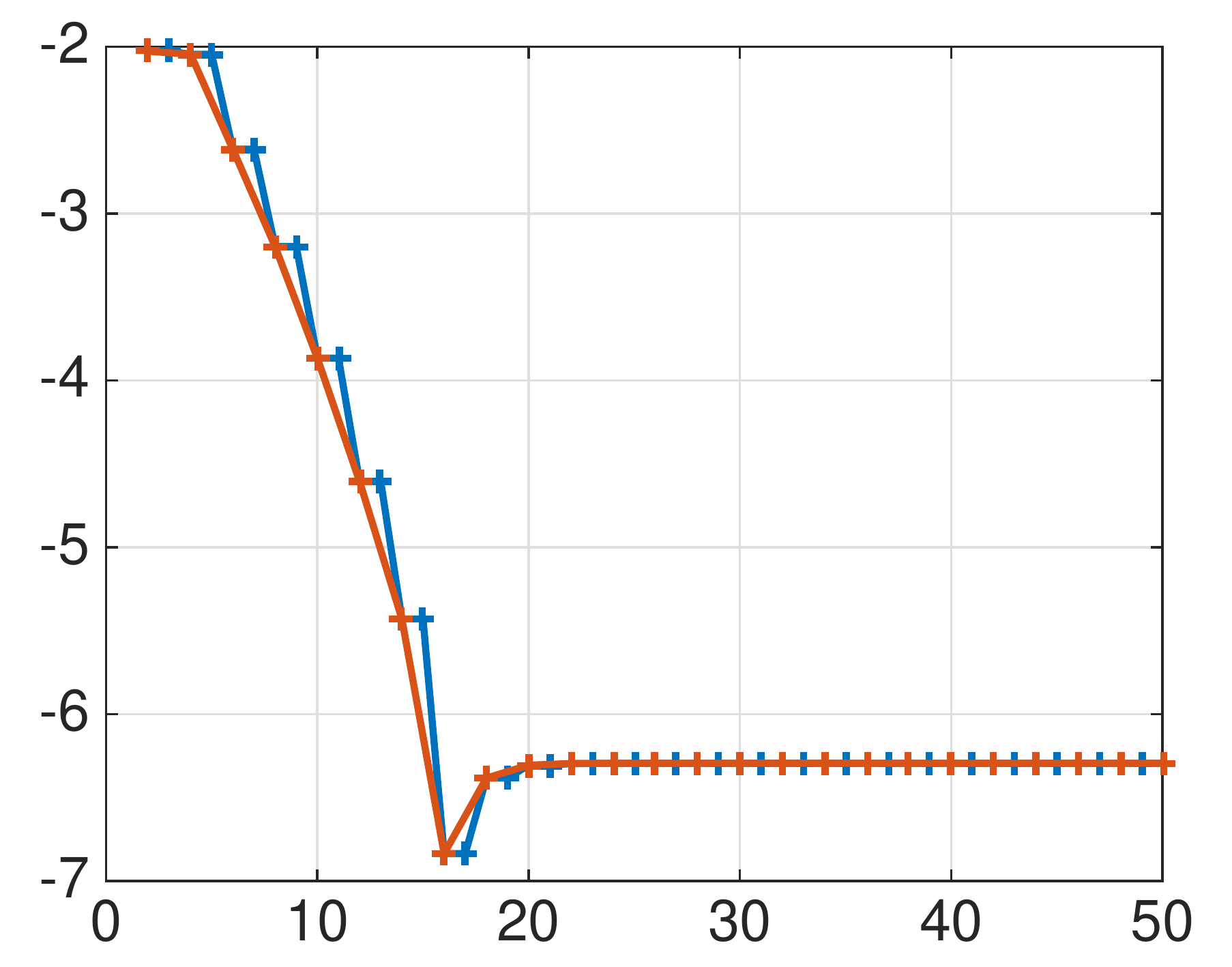}}
\caption{Convergence as a function of $N$ : $D(N)$, $\log_{10}\frac{|{\cal D}_{*}-D(N)|}{{\cal D}_{*}}$, $\kappa(N)$, $\log_{10}\frac{|{\cal K}_{*}-\kappa(N)|}{{\cal K}_{*}}$ vs. $N$. \label{fig.CasBDKn}}
\end{center}
\end{figure}
In Figure~\ref{eq.CasBVP} we present the convergence of the positive eigenvalues. As for the previous potential, we observe the exponential decay of $\lambda_{1}(\gamma)$ to 0 as $\gamma$ increases. 
For large values of $\gamma$ we note a change of slope in the log-log graph and the loss of monotonicity for  the eigenvalues $\lambda_{3},\lambda_{4},\lambda_{5}$. A further investigation would be necessary to decide  whether this behavior 
is a numerical artifact or whether it is the actual behavior of the eigenvalues.

In Figure~\ref{fig.CasBDK}, we display the behavior of the diffusion and drift coefficients as functions of $\gamma$ showing the exponential growth as $\gamma\to\infty$. 
The comparison with \eqref{eq.drift1D} justifies the quality of the approximation. 
In Figure~\ref{fig.CasBDKn} we present the convergence of the drift and diffusion coefficients with respect to the number of modes for $\gamma=1$ and $\gamma=5$: as in Case A, the coefficients are well captured with a few eigenmodes.

\subsection{Case C -- smooth potential with a linear drift\label{Sec.CASC}}
\begin{figure}[h!t]
\begin{center}
\begin{minipage}{.4\linewidth}
\hfill\includegraphics[width=5cm]{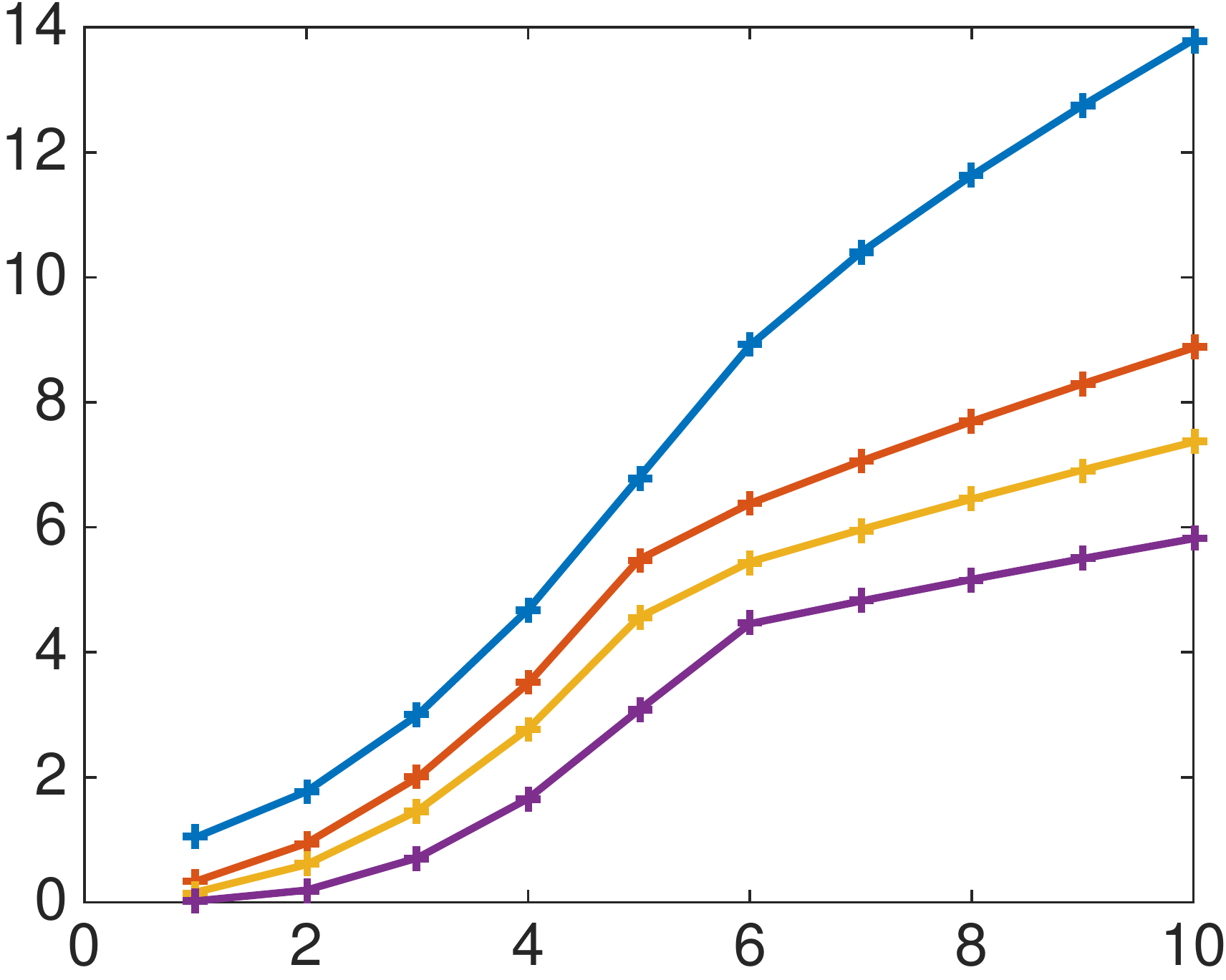}$\qquad$
\end{minipage}
\begin{minipage}{.2\linewidth}
\includegraphics[height=1.5cm]{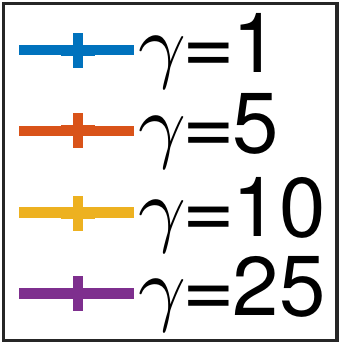}
\end{minipage}
\caption{$\lambda_{1}(\gamma,\delta)$ vs. $\delta$ for $\gamma=1, 5, 10, 25$.\label{fig.CasCVP}}
\end{center}
\end{figure}

\begin{figure}[h!t]
\begin{center}
\begin{tabular}{cc|c}
\includegraphics[height=3.5cm]{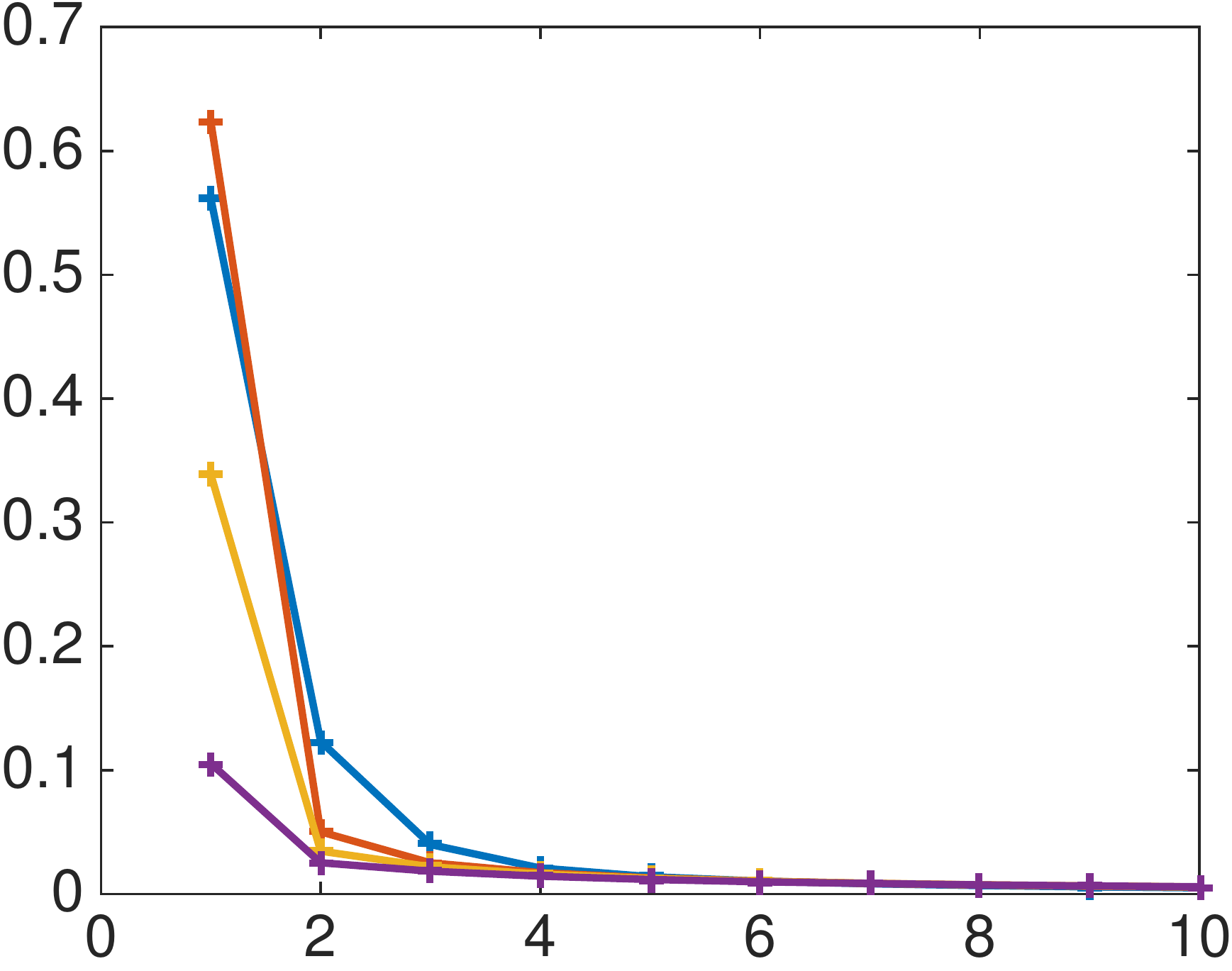}
\qquad&\includegraphics[height=3.5cm]{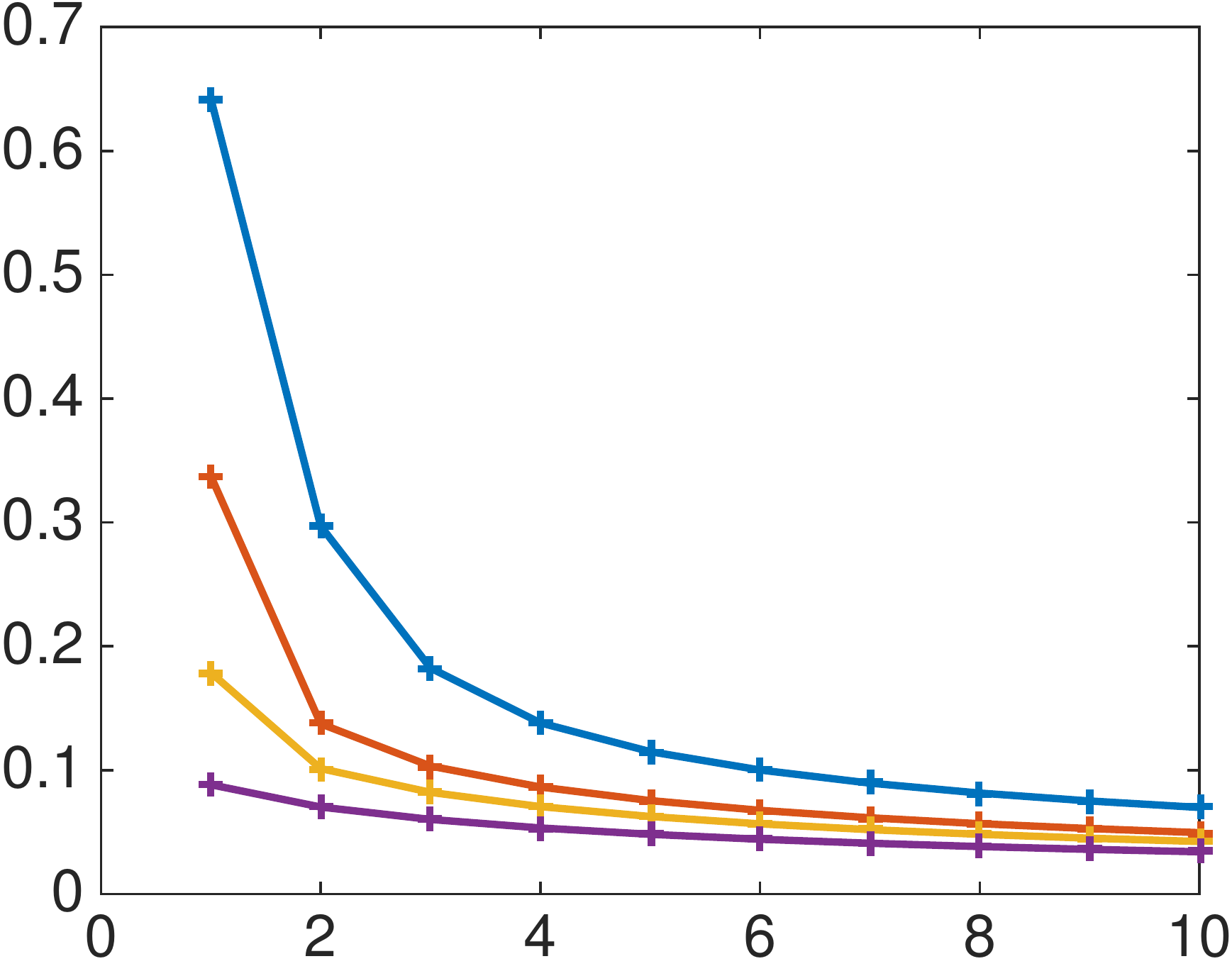}
\qquad&\quad\includegraphics[height=3.5cm]{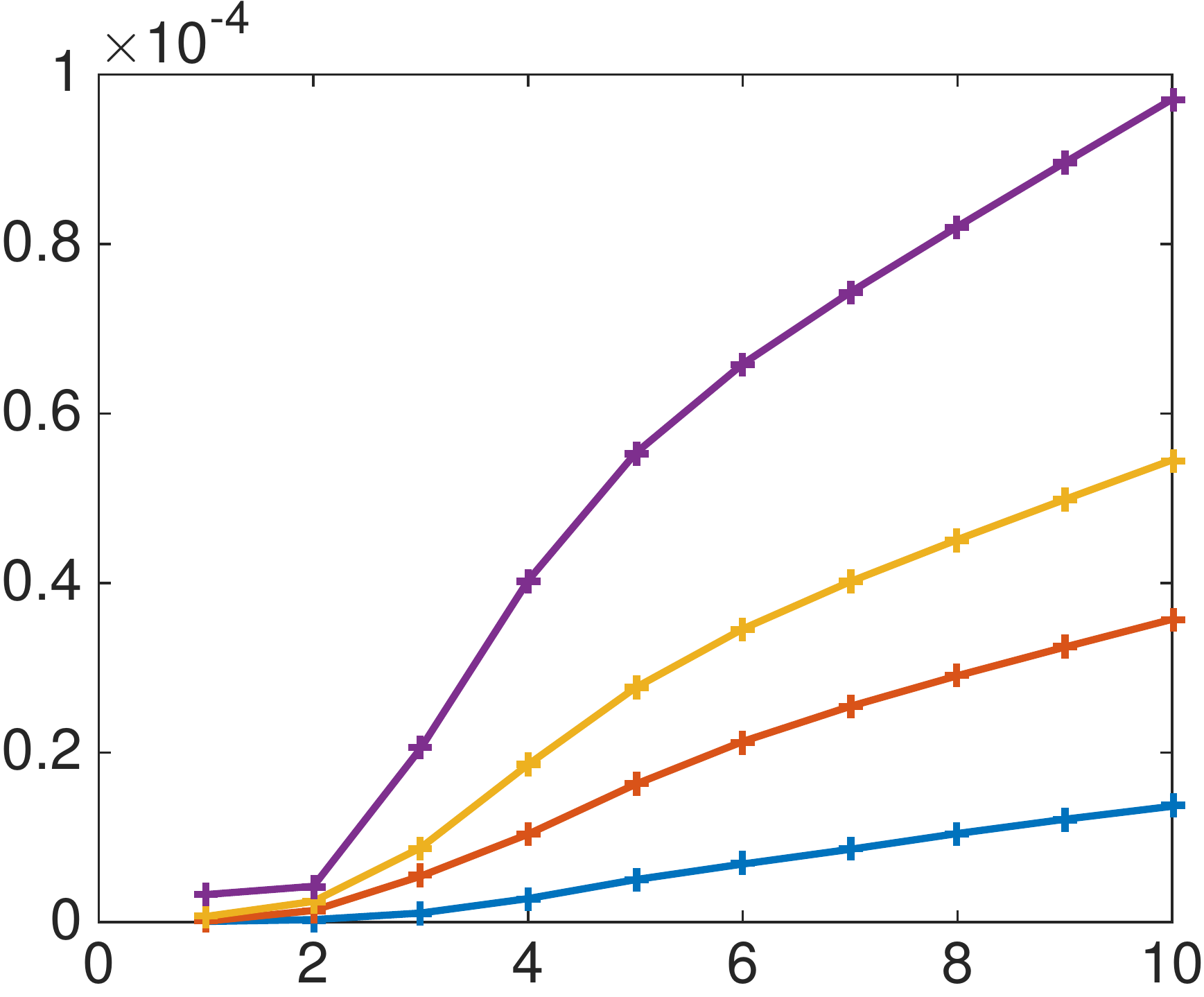}\\
$\delta\mapsto{\mathcal D}(\delta)$ & $\delta\mapsto{\mathcal K}(\delta)$& $\delta\mapsto\frac{|\mathcal K_{*}(\delta)-\mathcal K(\delta)|}{\mathcal K_{*}(\delta)}$ \\[5pt]
\includegraphics[height=3.5cm]{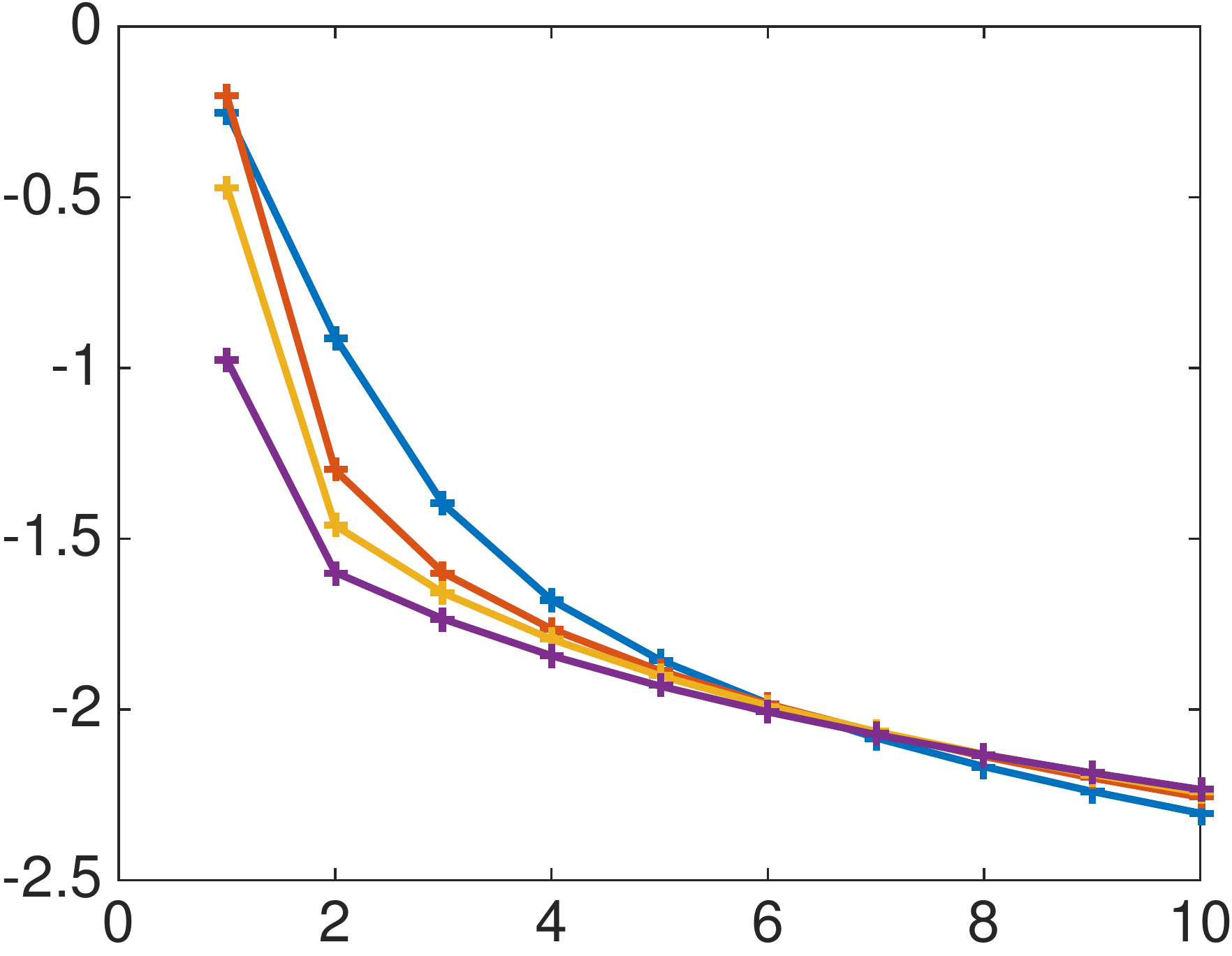}
\qquad&\includegraphics[height=3.5cm]{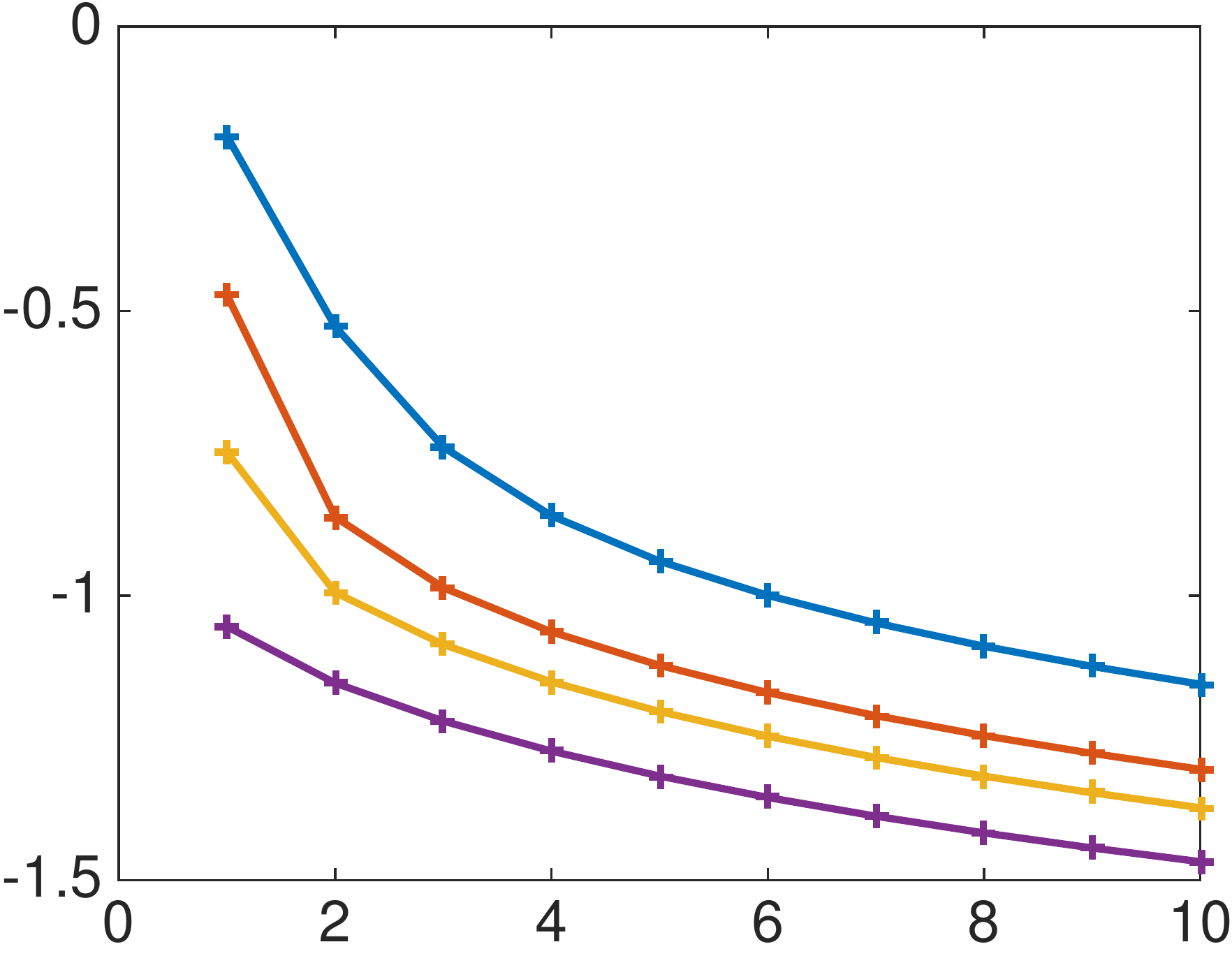}
\qquad&\quad\includegraphics[height=3.5cm]{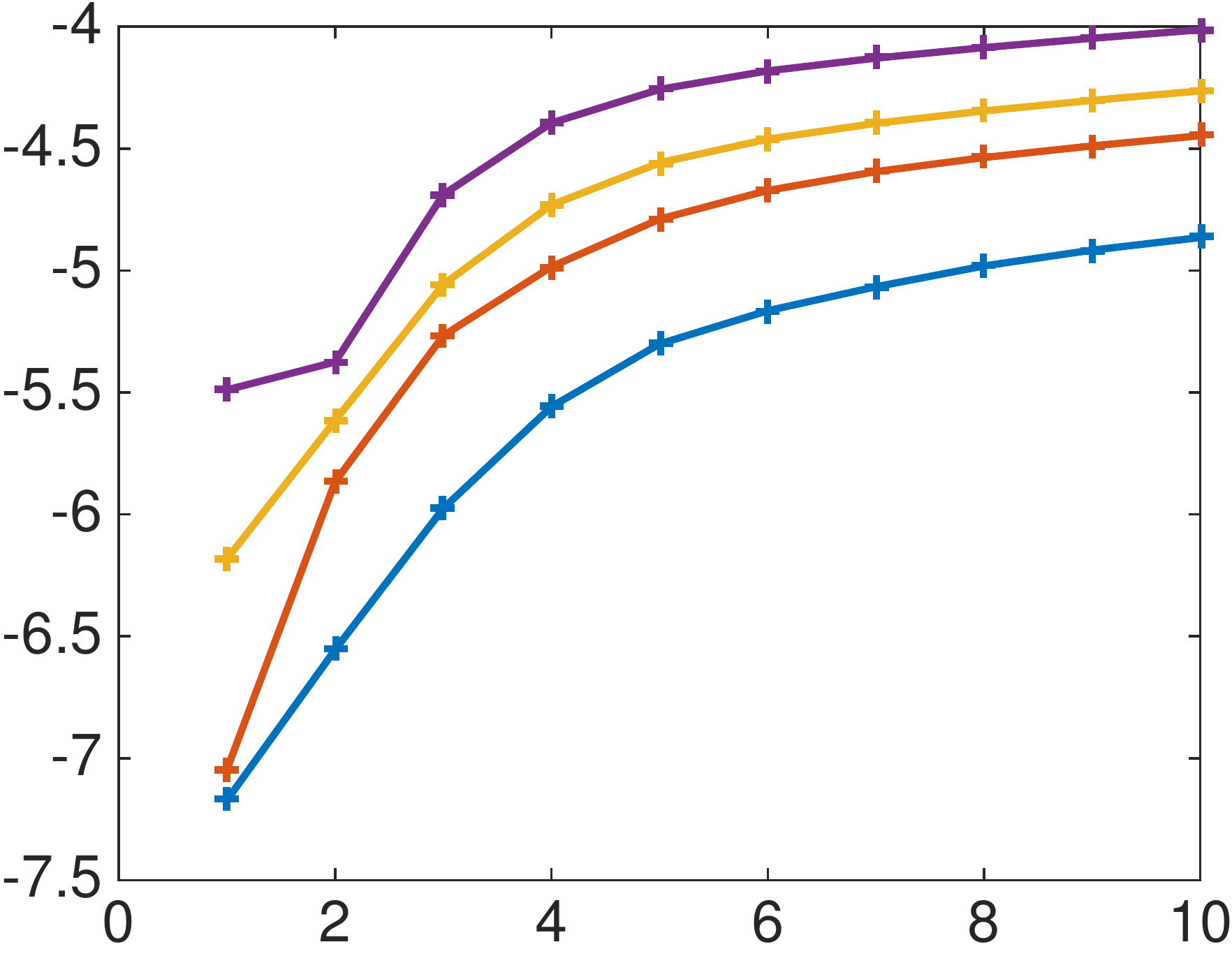}\\
$\delta\mapsto\log_{10}{\mathcal D}(\delta)$ & $\delta\mapsto\log_{10}{\mathcal K}(\delta)$ & $\delta\mapsto \log_{10}\frac{|\mathcal K_{*}(\delta)-\mathcal K(\delta)|}{\mathcal K_{*}(\delta)}$
\end{tabular}
\caption{
Diffusion and drift coefficients as functions of $\delta$ for the nonsymmetric potential~\eqref{e:nonsymm}, with $\gamma=1, 5, 10, 25$ (with legend as in Figure~\ref{fig.CasCVP}).
\label{fig.CasCDK}}
\end{center}
\end{figure}

\begin{figure}[h!t]
\begin{center}
\subfigure[$\gamma=1$, $\delta=1$]{\includegraphics[height=3.3cm]{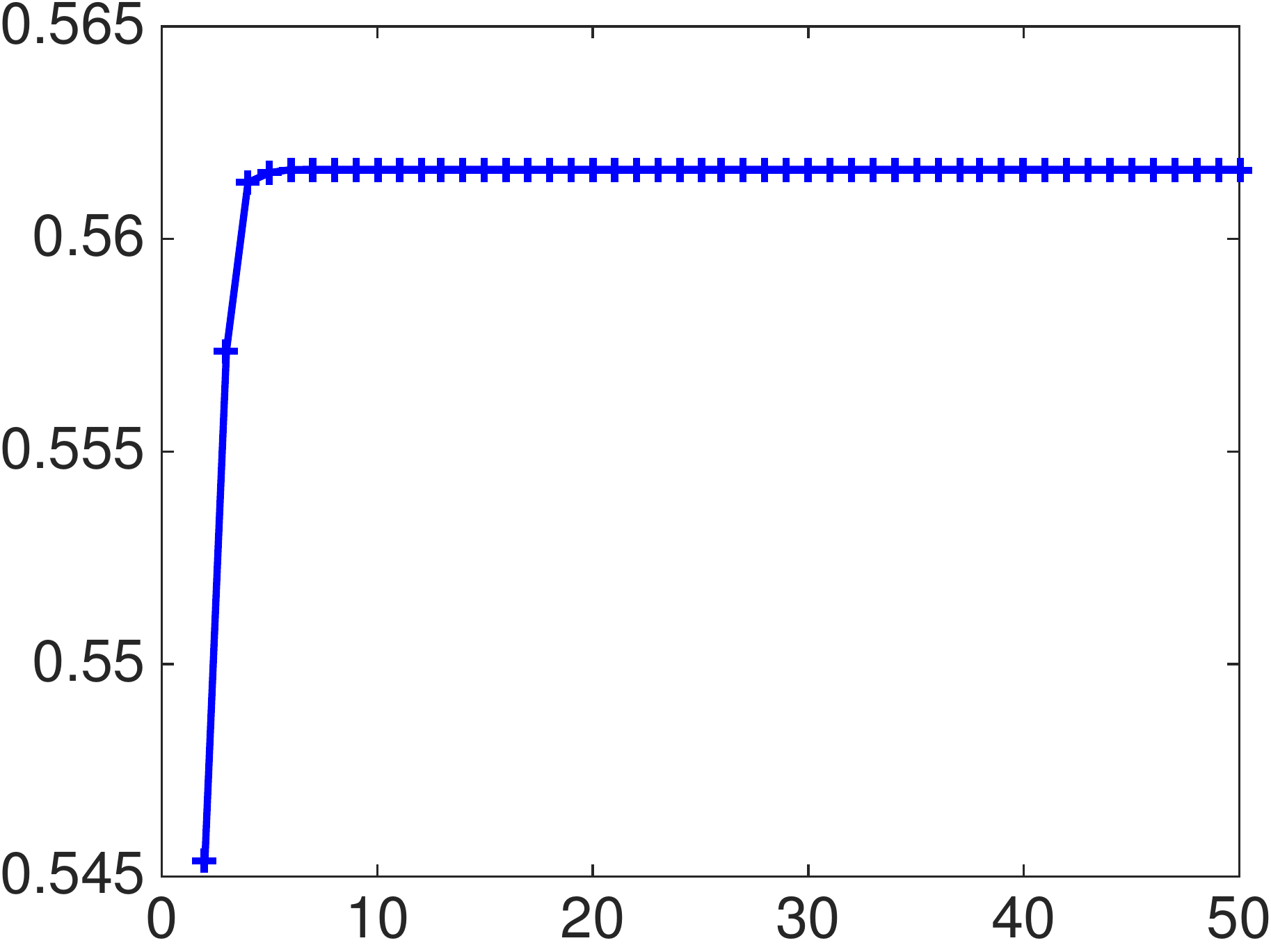}
\includegraphics[height=3.3cm]{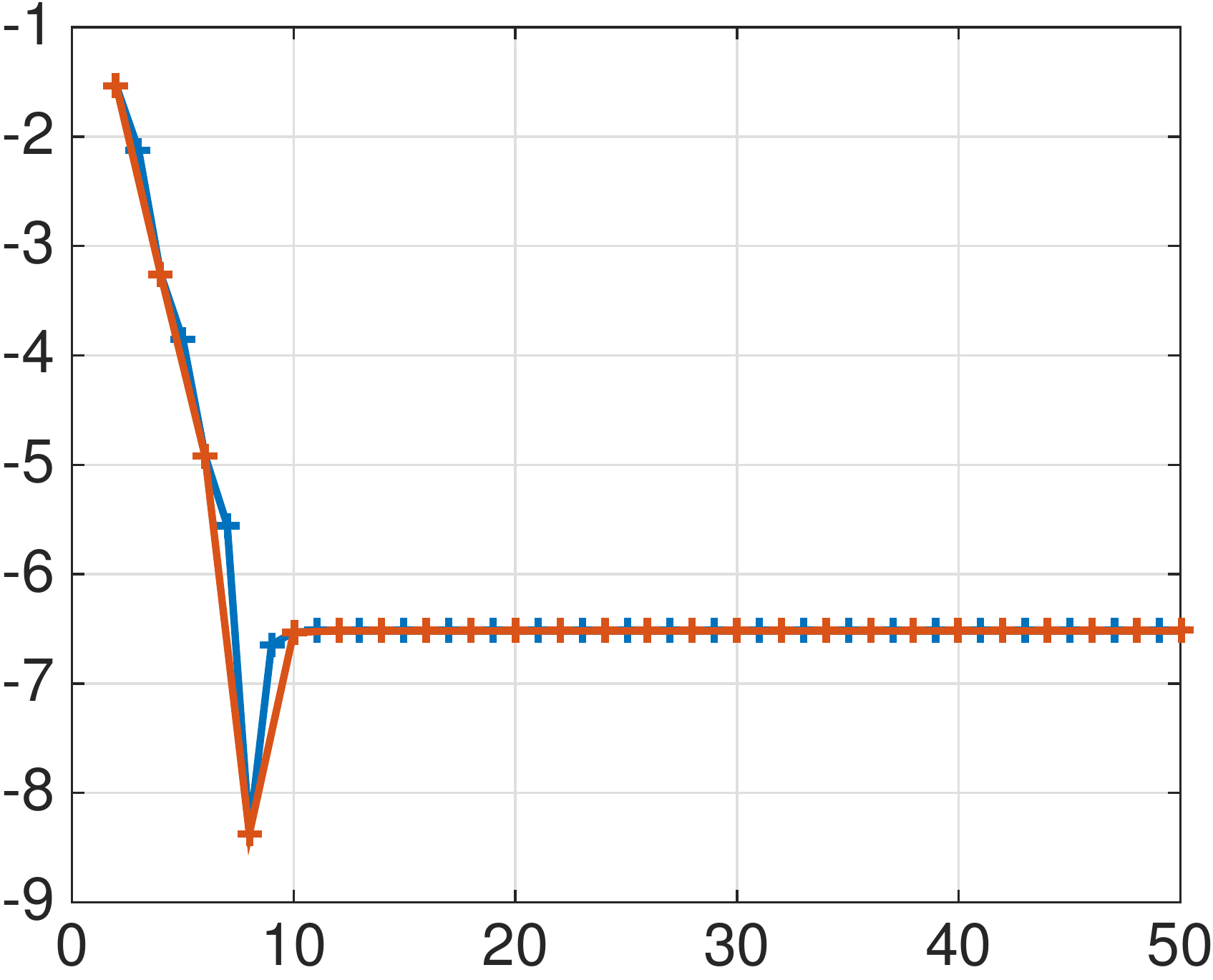}
\includegraphics[height=3.3cm]{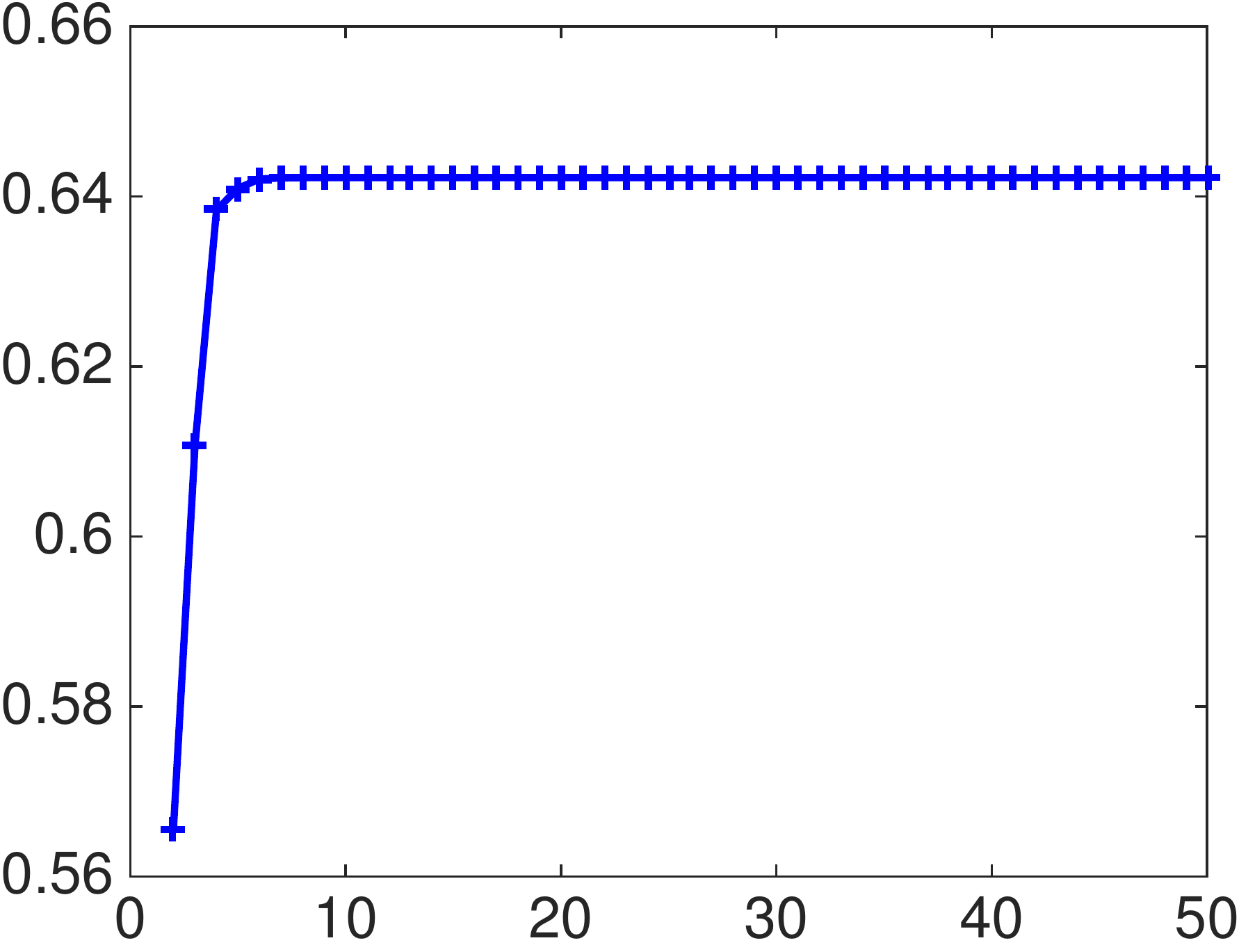}
\includegraphics[height=3.3cm]{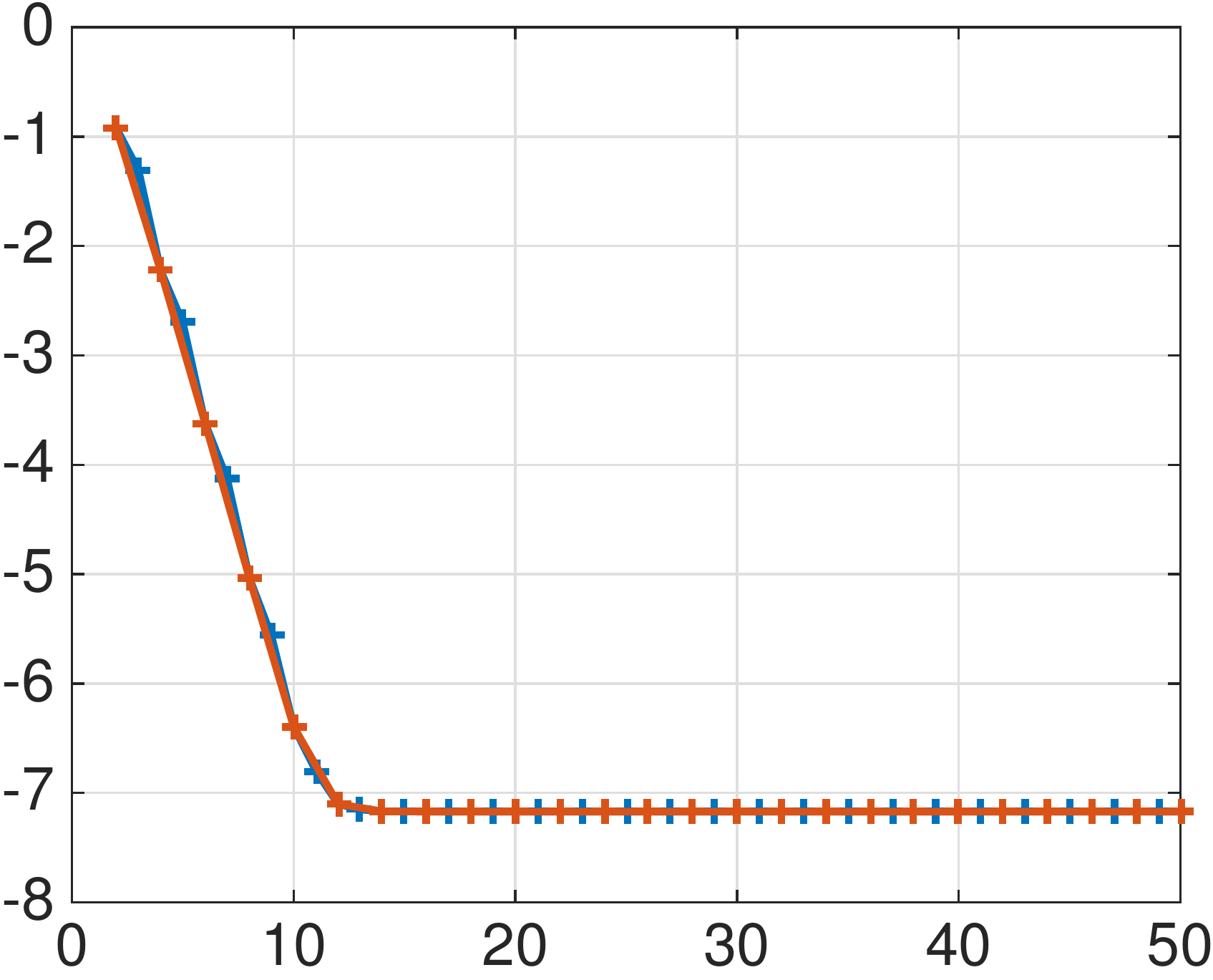}}
\subfigure[$\gamma=1$, $\delta=5$]{\includegraphics[height=3.3cm]{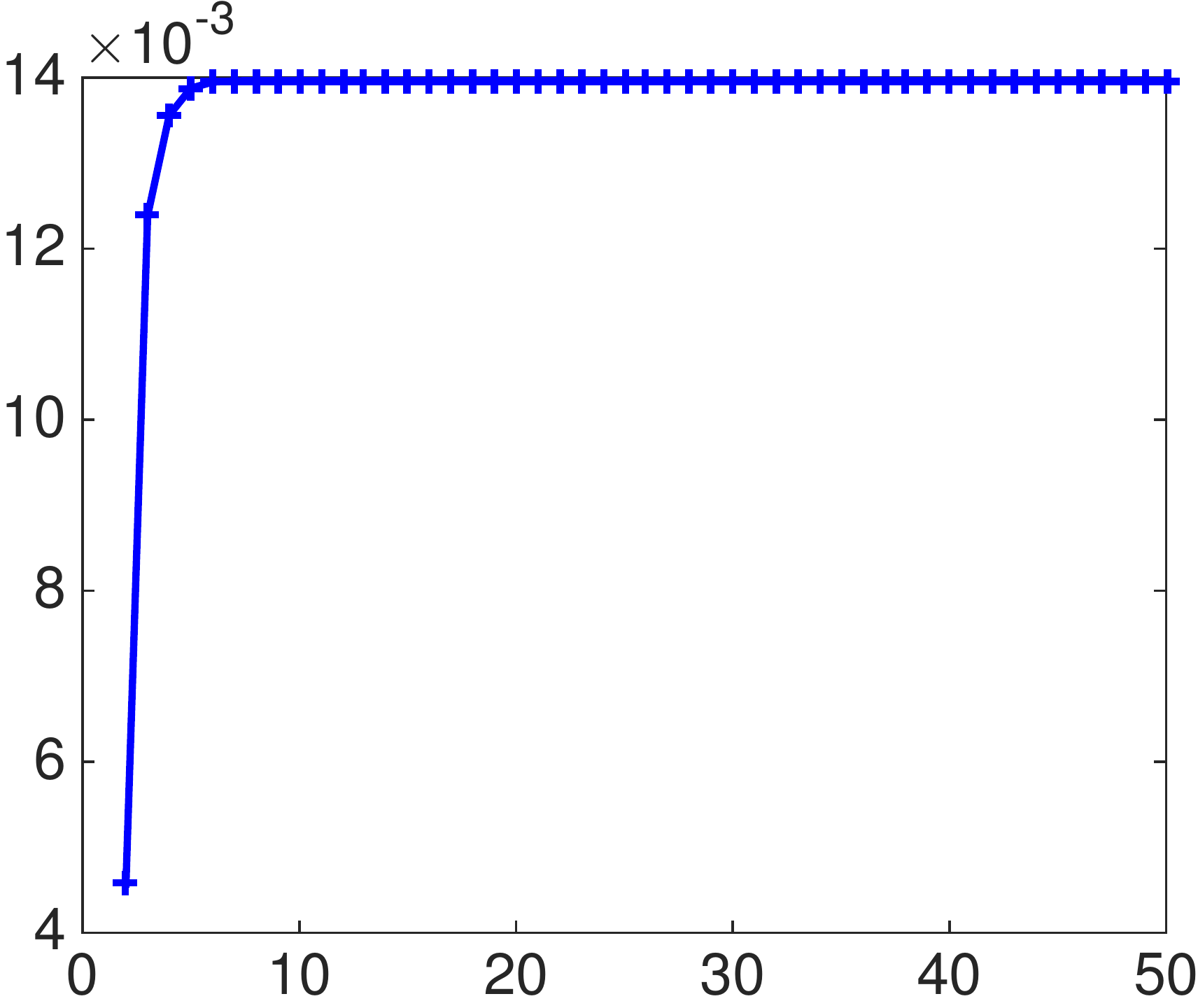}
\;\includegraphics[height=3.3cm]{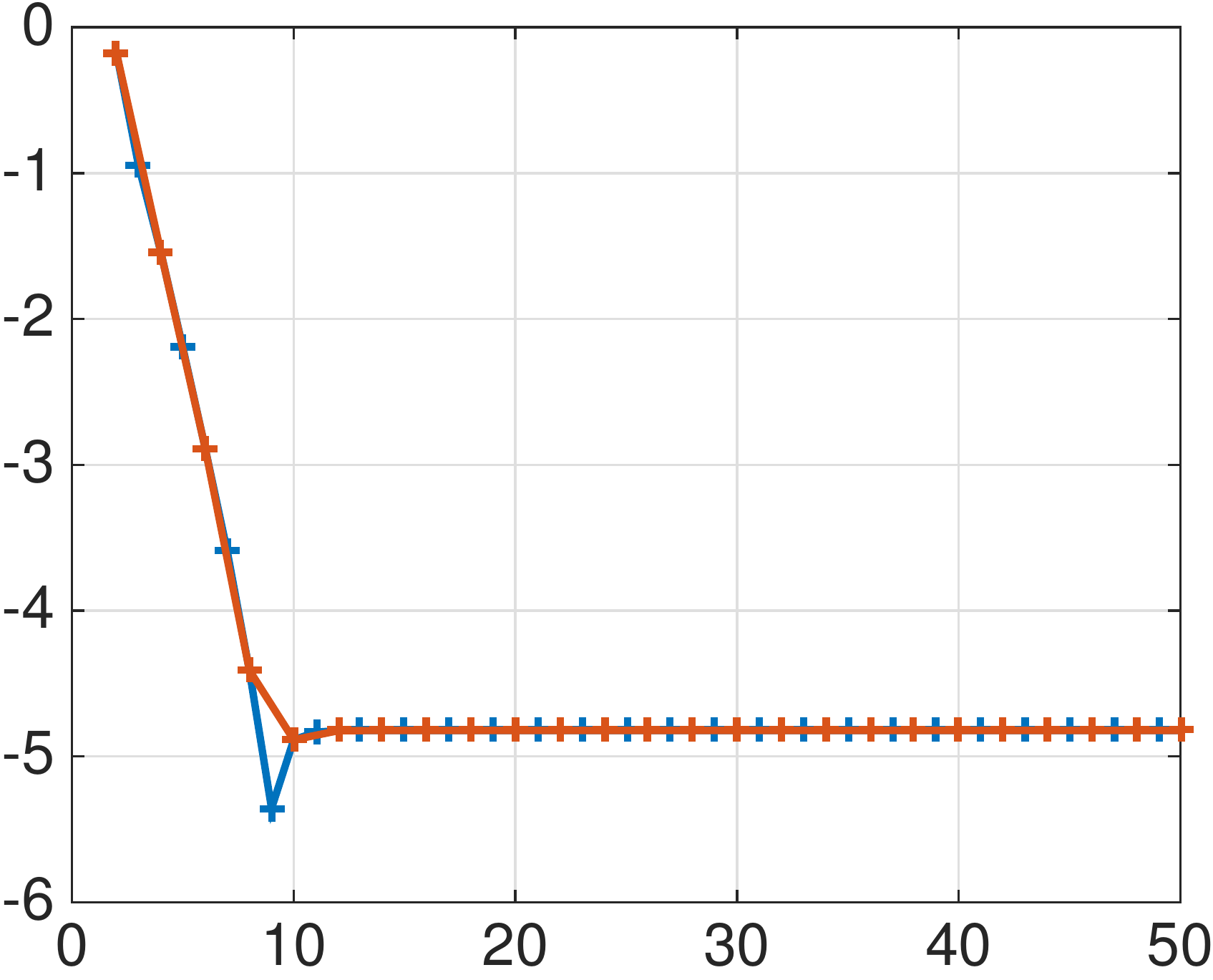}
\;\includegraphics[height=3.3cm]{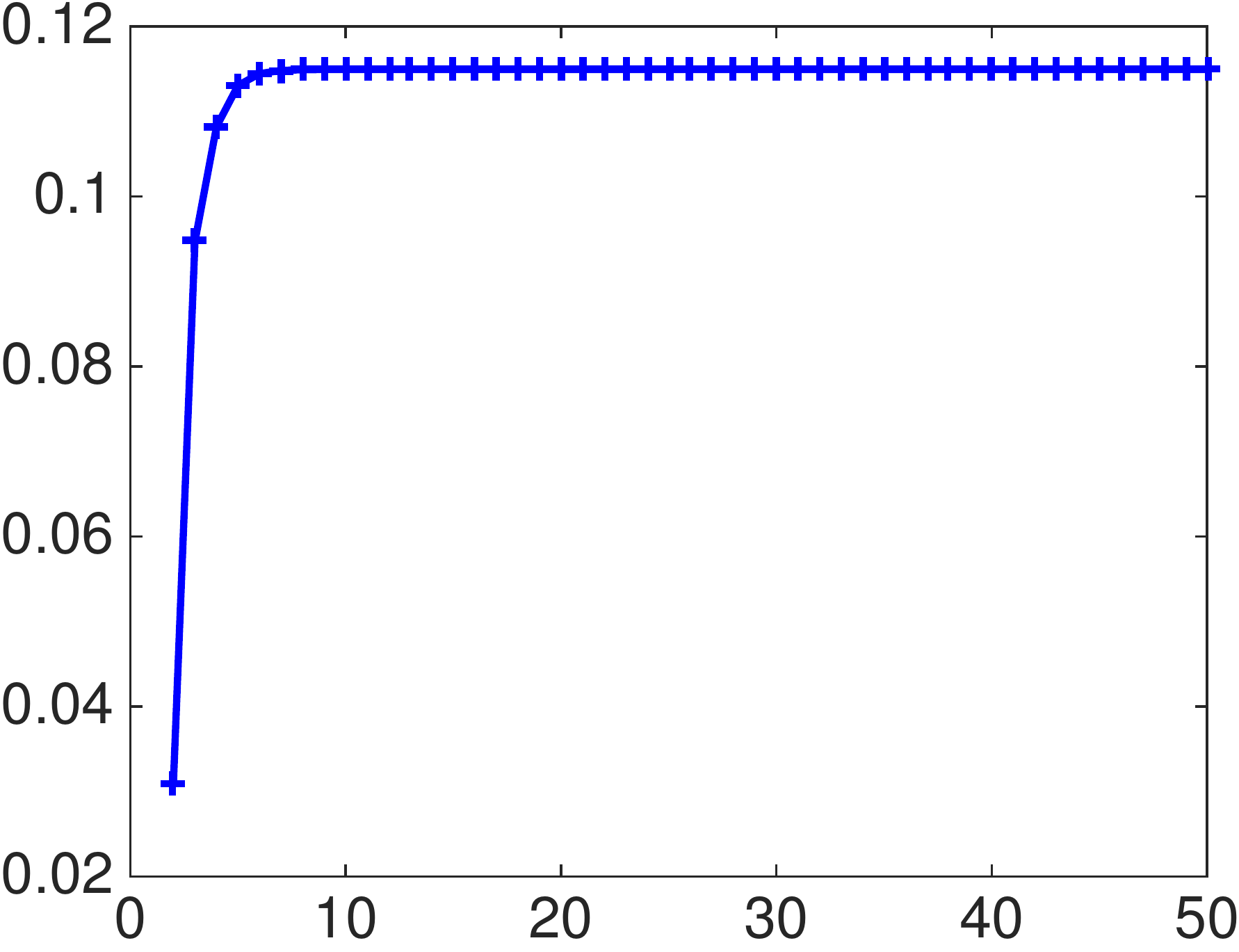}
\;\includegraphics[height=3.3cm]{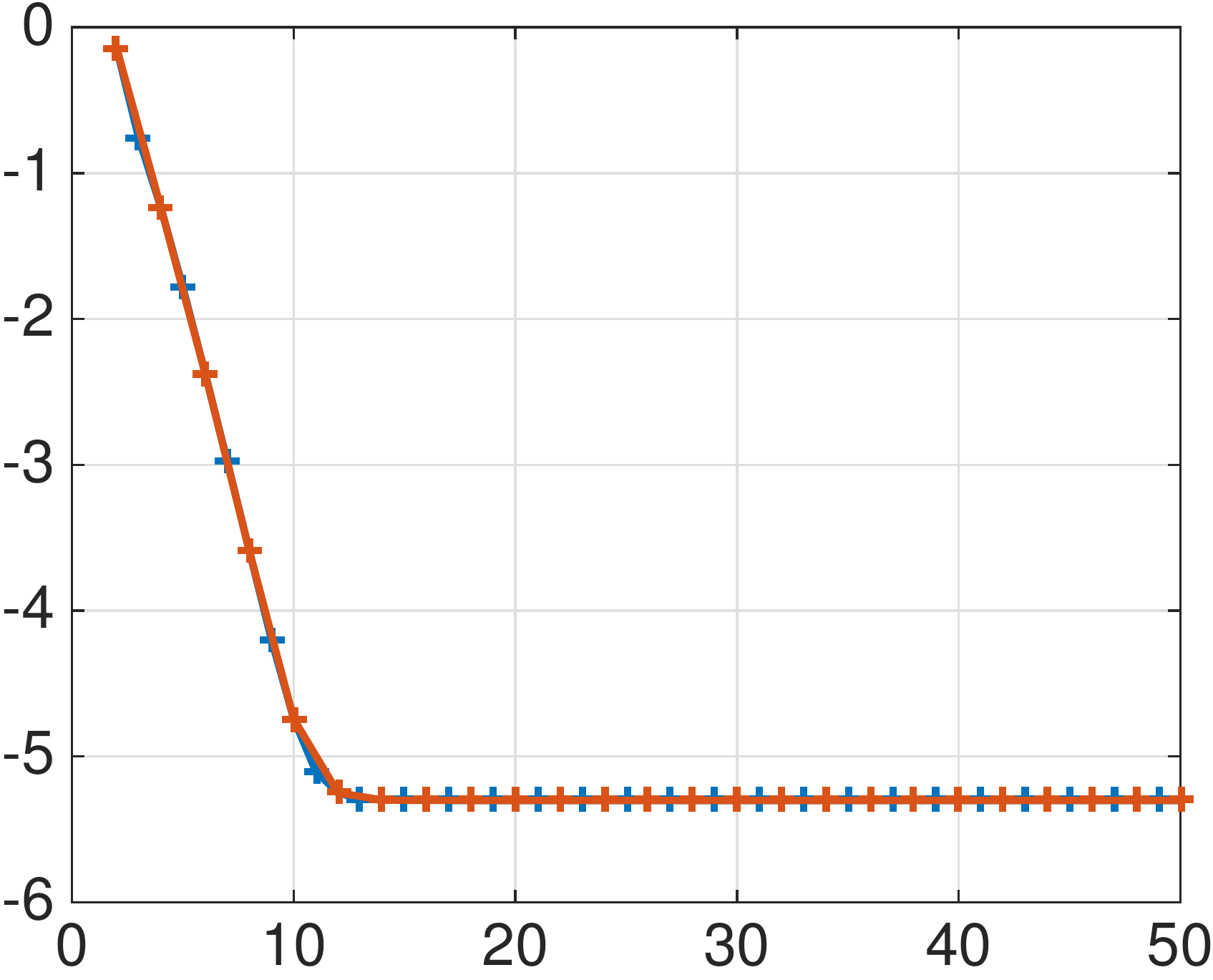}}
\subfigure[$\gamma=1$, $\delta=10$]{\includegraphics[height=3.3cm]{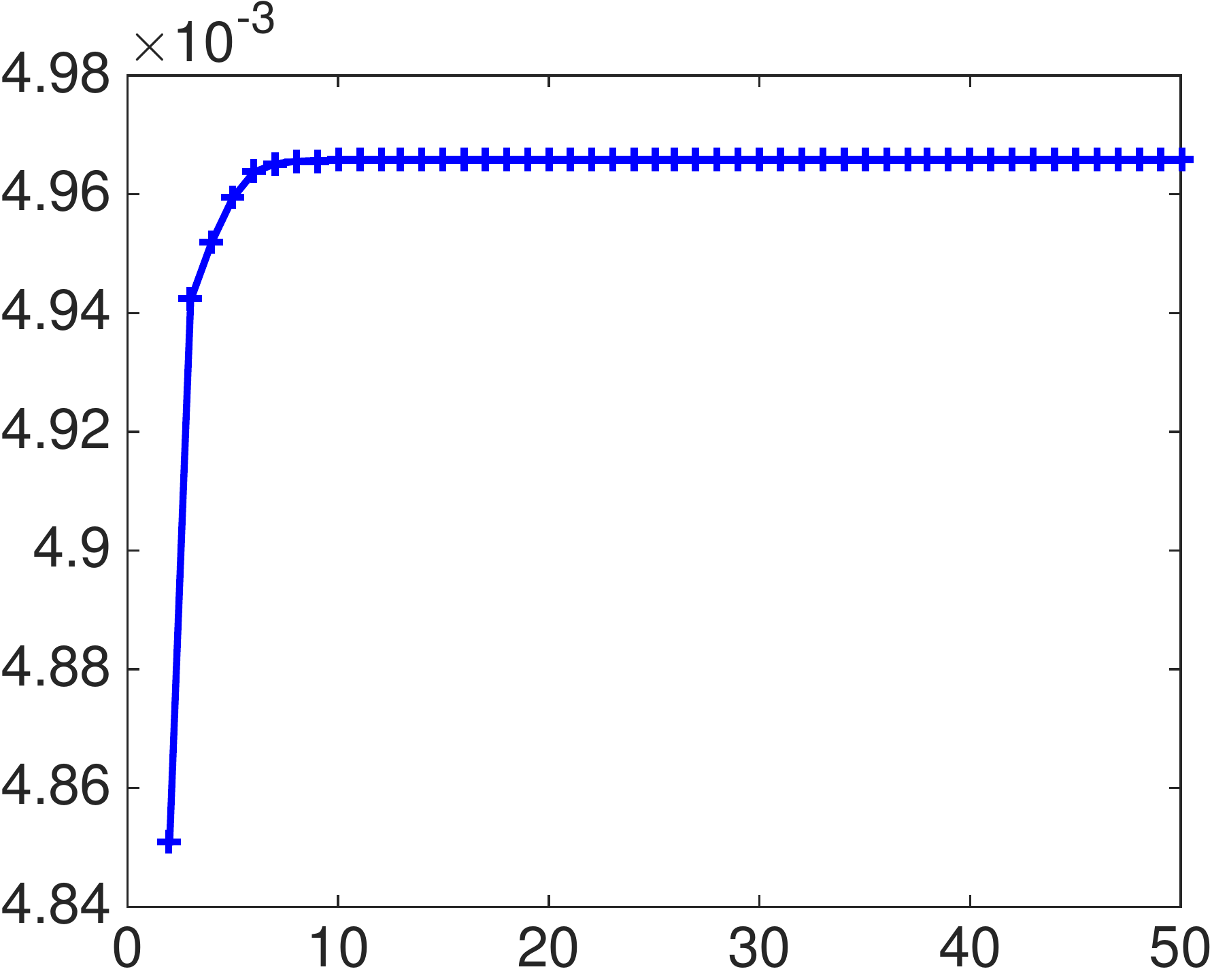}
\includegraphics[height=3.3cm]{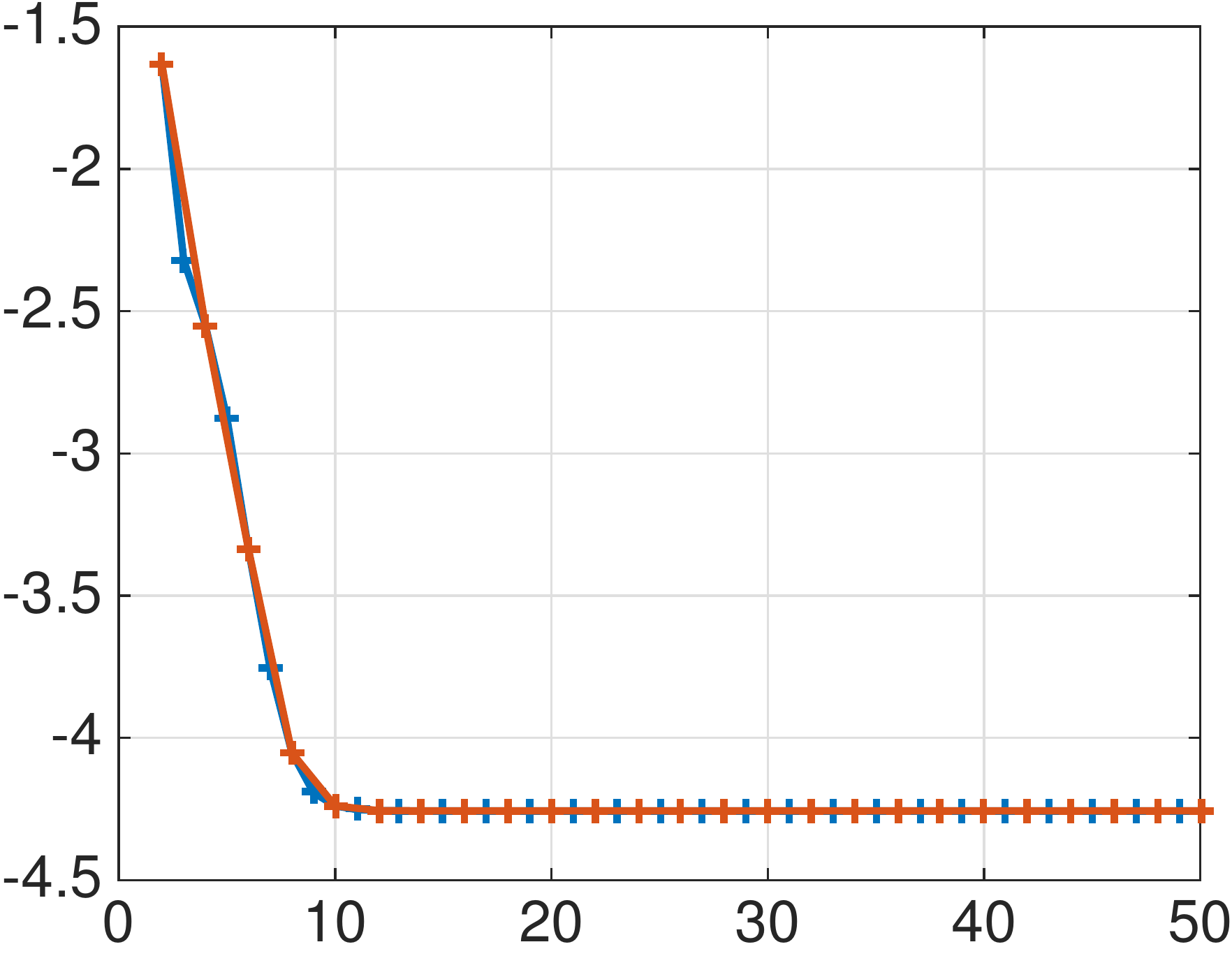}
\includegraphics[height=3.3cm]{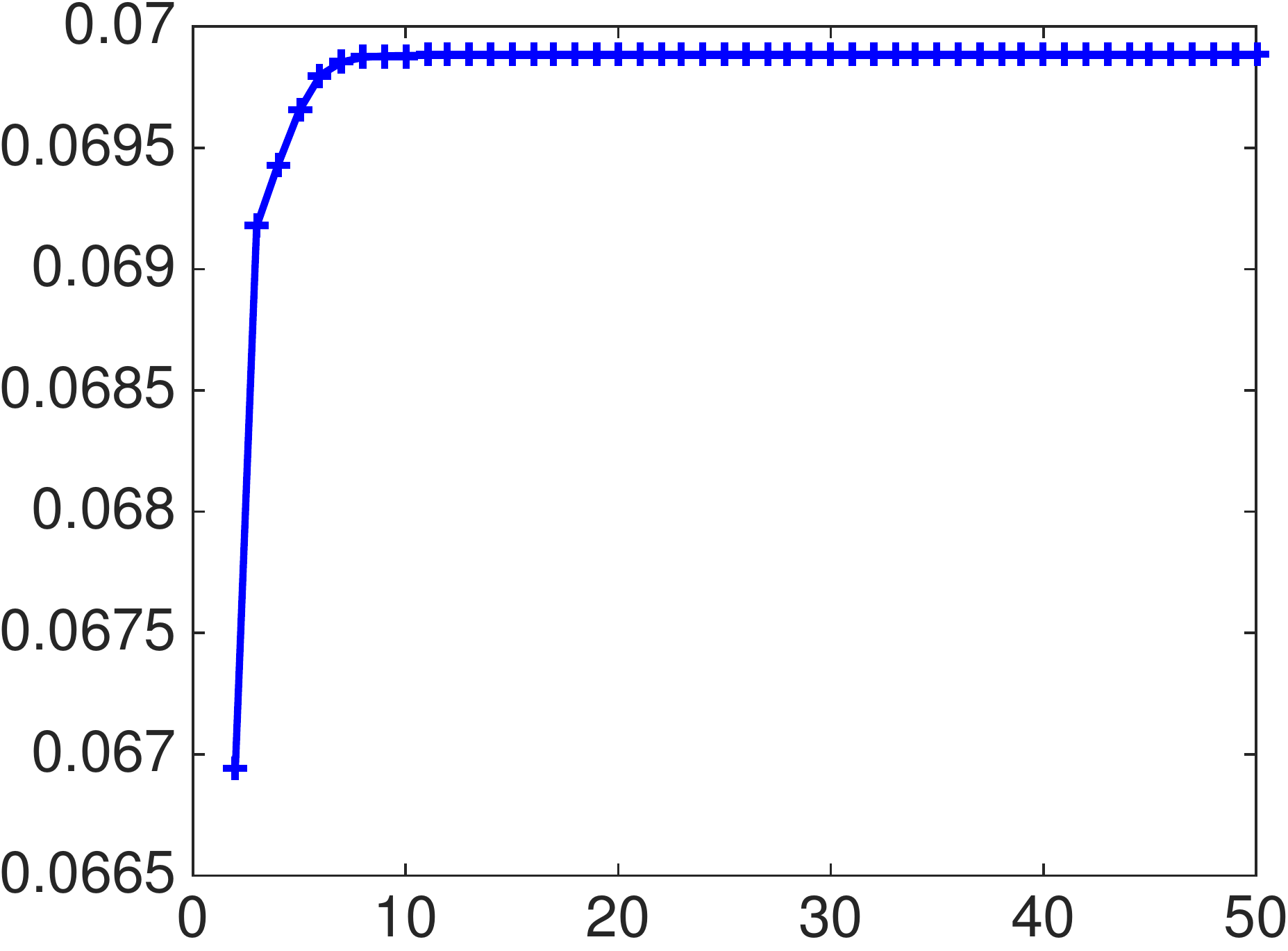}
\includegraphics[height=3.3cm]{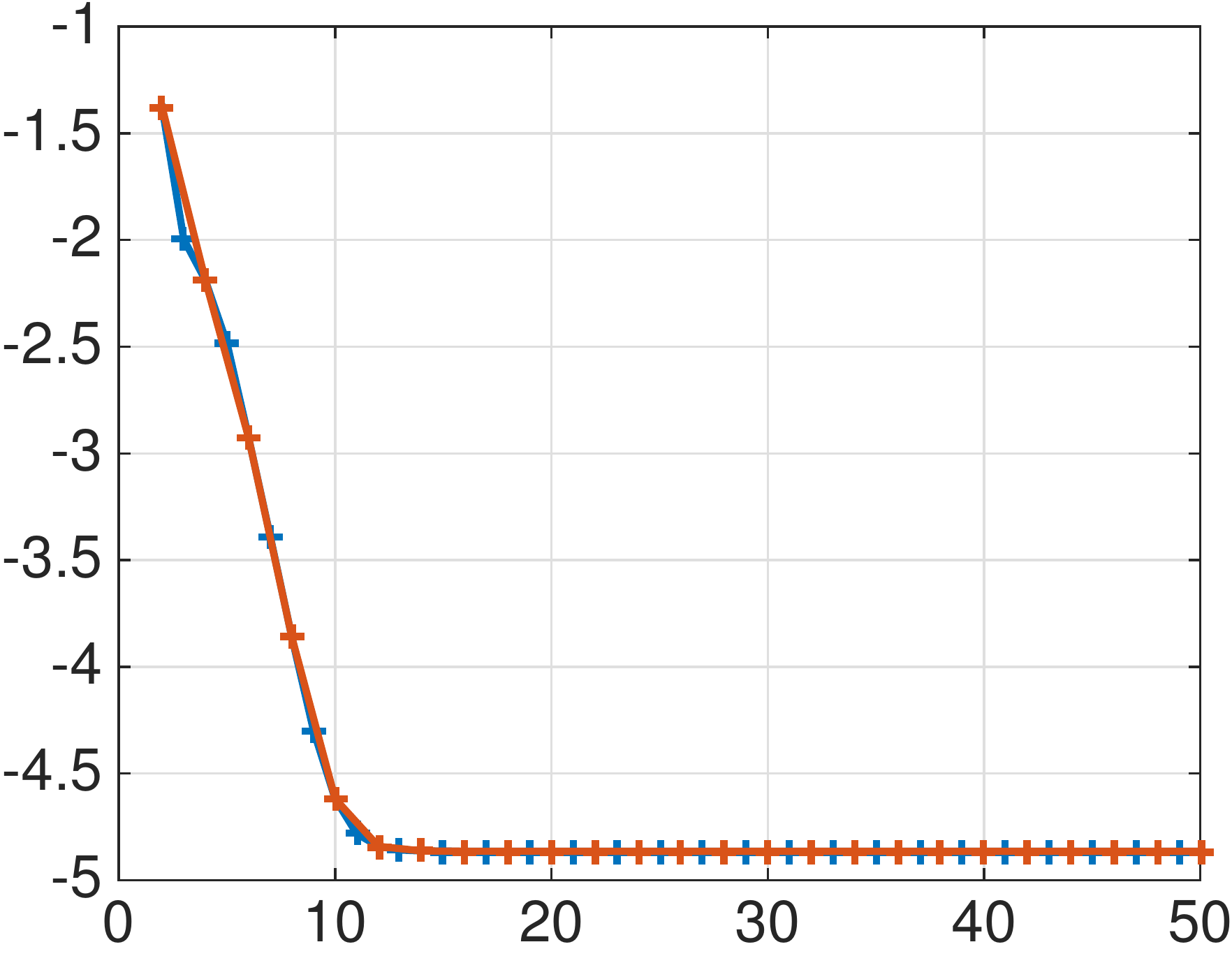}}
\caption{
Convergence as a function of $N$ for the nonsymmetric potential~\eqref{e:nonsymm} : $D(N)$, $\log_{10}\frac{|{\cal D}_{*}-D(N)|}{{\cal D}_{*}}$, $\kappa(N)$, $\log_{10}\frac{|{\cal K}_{*}-\kappa(N)|}{{\cal K}_{*}}$ vs. $N$. \label{fig.CasCDKn}}
\end{center}
\end{figure}

Now we present our computations for the nonsymmetric potential
$$
W(v)=\tfrac{1}{4\gamma} v^4 - \tfrac{1}2 v^2 -\delta v \qquad \mbox{with } \gamma,\delta>0\,.
$$
The compatibility condition, necessary so that we can apply Fredholm's alternative, does not hold: here
$$ \cV := \int_{\R} v M(v)\ud v\neq 0.$$ 
However, we can adapt the asymptotics in order to handle this situation where the flux of the equilibrium state does not vanish. Starting from \eqref{eq.nnn}, we set
$$\widetilde f^{\varepsilon}(t,r,v) = f^{\varepsilon}\left(t,r+\frac{\cV t}{\sqrt\varepsilon},v\right),$$
with $f^\varepsilon$ solution of \eqref{eq.nnn}. 
We check that 
\[
\partial_{t}\widetilde f^{\varepsilon} +\frac{1}{\sqrt\varepsilon}\left((v-\cV)\nabla_{r}\widetilde f^{\varepsilon}
+\div_{v} \big( (\nabla_{r}U\star\widetilde\rho^{\varepsilon}) \widetilde f^{\varepsilon}\big)\right)
=\frac 1\varepsilon Q(\widetilde f^{\varepsilon}),
\]
with $\widetilde\rho^{\varepsilon}=\int \widetilde f^{\varepsilon} \ud v$. 
We perform the Hilbert expansion on $\widetilde f^{\varepsilon}$ as in Section~\ref{Section2}: we still have $Q(\widetilde f^{(0)})=0$ and thus $\widetilde f^{(0)}(t,r,v) = \rho(t,r) M(v)$ at leading order, while $\widetilde f^{(1)}$ is defined by inverting
\begin{align*}
Q(\widetilde f^{(1)}) 
&= (v-\cV)\nabla_{r}\widetilde f^{(0)} - \div_{v}\big((\nabla_{r}U\star \rho)\widetilde f^{(0)}\big)\\
&= (v-\cV)M\nabla_{r}\rho - \frac1\theta\nabla_{v}W\ M\cdot (\nabla_{r}U\star \rho)\rho.
\end{align*}
Therefore, Equations \eqref{e:invert} become
\begin{align*}
Q(\chi_i) &= (v-\cV)_i M(v)  \; , \\
Q(\psi_i) &= \ds\frac1\theta\frac{\partial W}{\partial v_i}(v) M(v).
\end{align*}
The compatibility condition is now satisfied and we get $\widetilde f^{(1)} = \chi\cdot\nabla_{r}\rho + \psi\cdot (\nabla_{r}U\star \rho)\rho$.
Considering the $\varepsilon^0$ terms in the expansion, we finally arrive at
\begin{equation*}
\partial_t \rho-\nabla_r \cdot\big(\mathcal{D} \nabla_r\rho+ \mathcal{K} (\nabla_r U \star
\rho)\rho\big)=0\,,
\end{equation*}
where the diffusion and drift matrices are given by the following analog of \eqref{DD2}:
\begin{equation*}
\mathcal{D} = -\int_{\RR^d} (v-\cV)\otimes\chi \ud v \qquad \mbox{and}
\qquad \mathcal{K} = -\int_{\RR^d} (v-\cV)\otimes\psi \ud v\,.
\end{equation*}
The analysis of this case can be performed as in Appendix~\ref{app.A}, see \cite{GOME} for a similar problem.
\begin{figure}[h!t]
\begin{center}
\subfigure[$\gamma=1$]{\begin{tabular}{ccc}
$\delta=1$ & $\delta=5$ & $\delta=10$\\
\includegraphics[height=3cm]{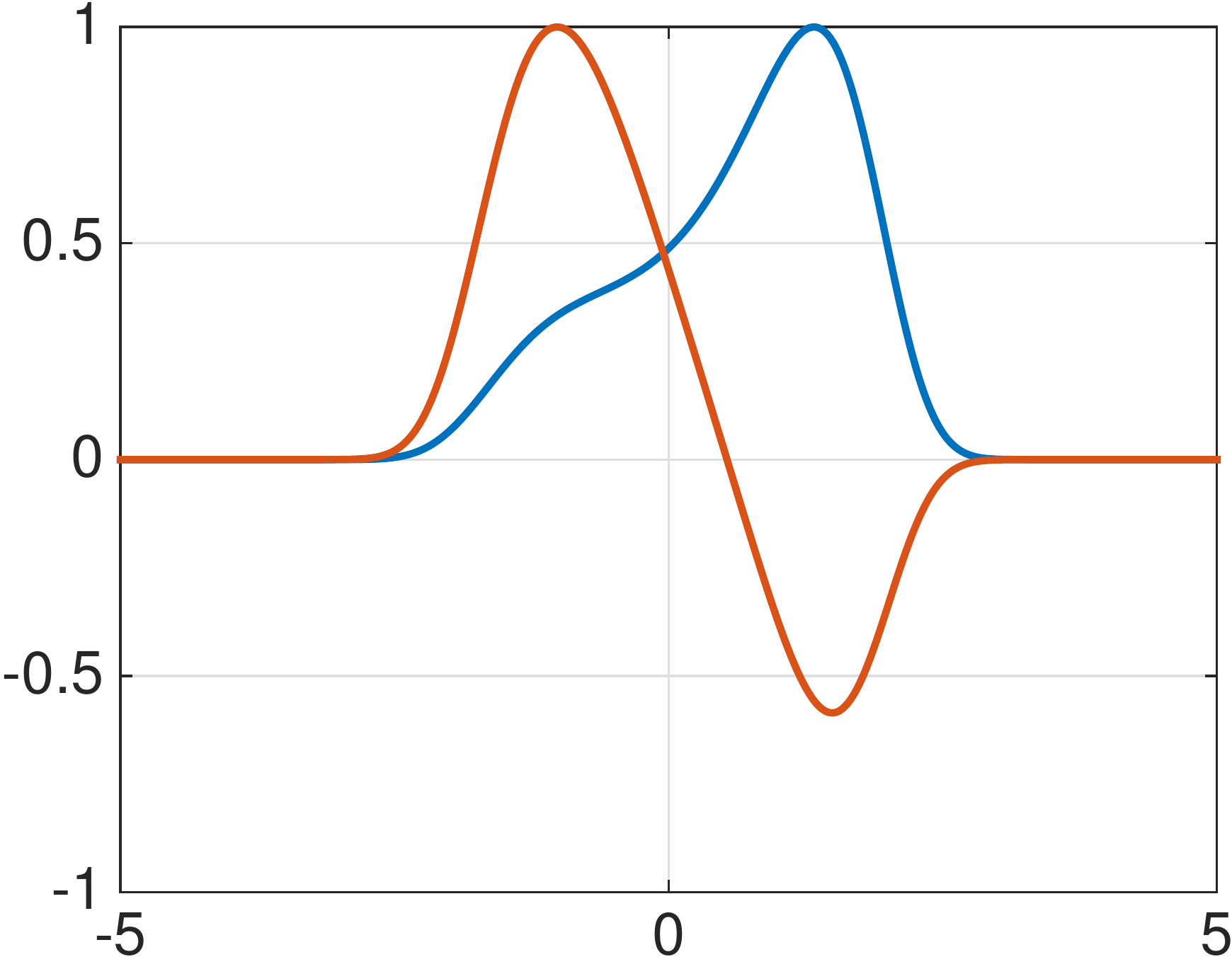}
\qquad&\includegraphics[height=3cm]{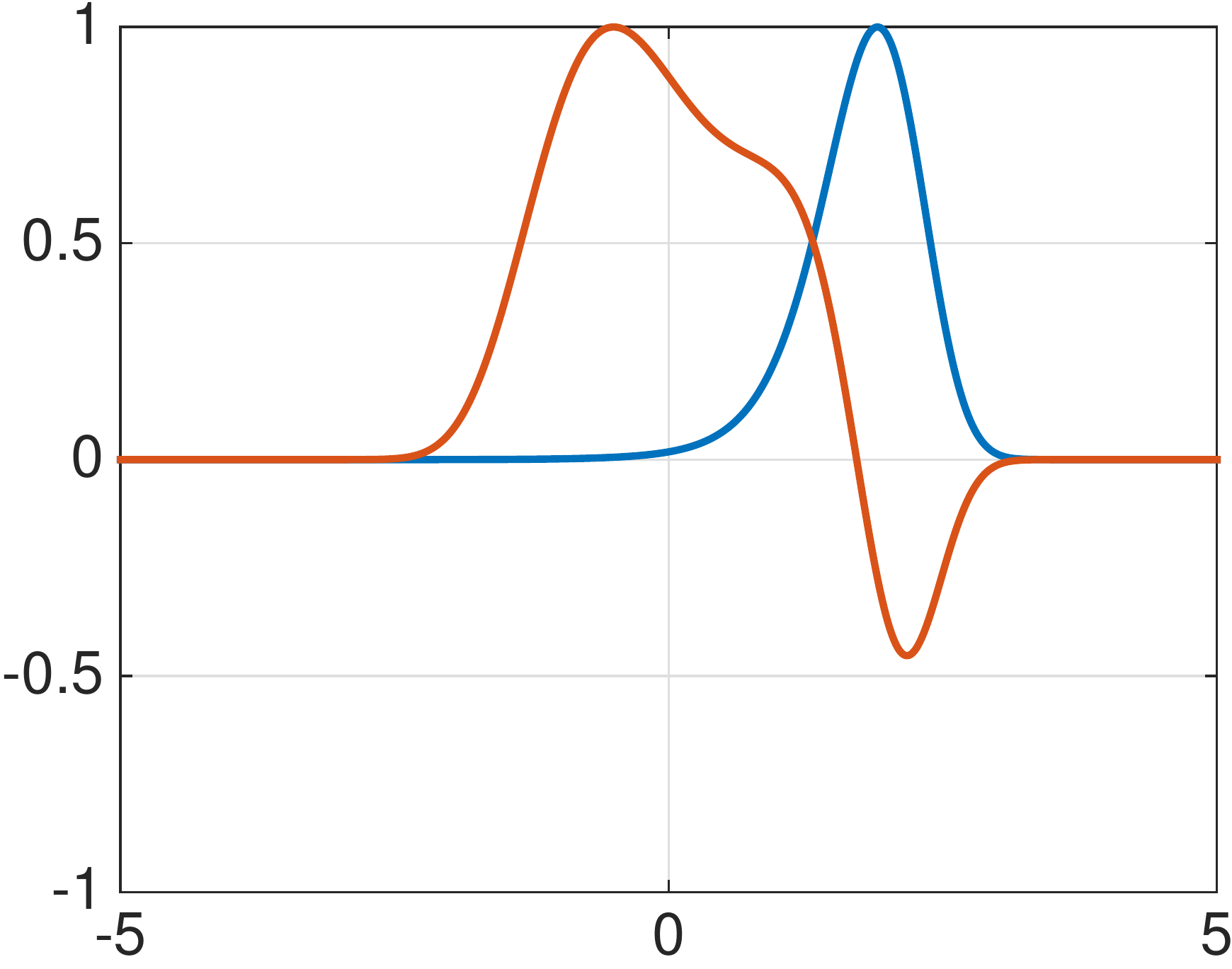}
\qquad&\includegraphics[height=3cm]{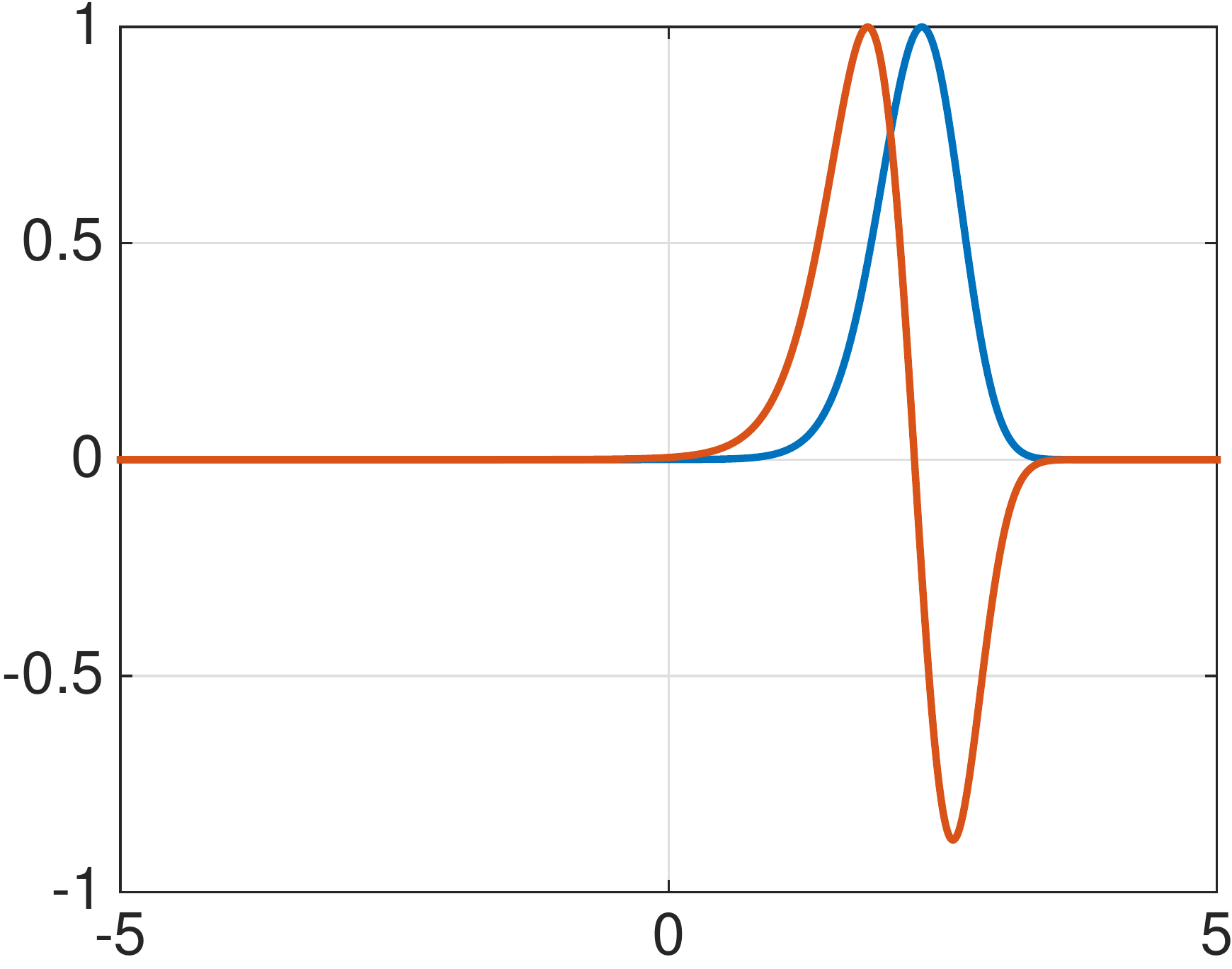}\\
\end{tabular}}
\subfigure[$\gamma=10$]{\begin{tabular}{ccc}
$\delta=1$ & $\delta=5$ & $\delta=10$\\
\includegraphics[height=3cm]{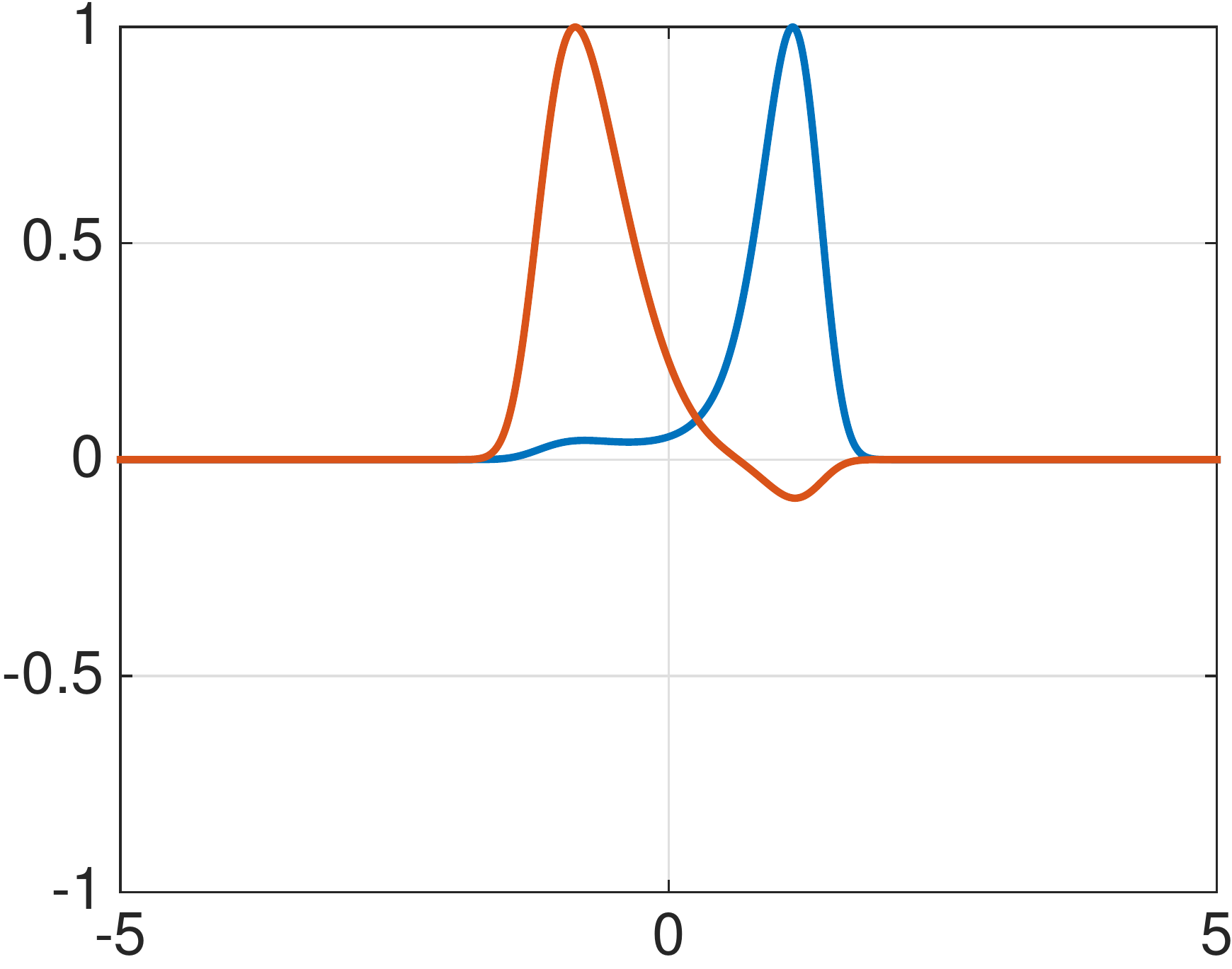}
\qquad&\includegraphics[height=3cm]{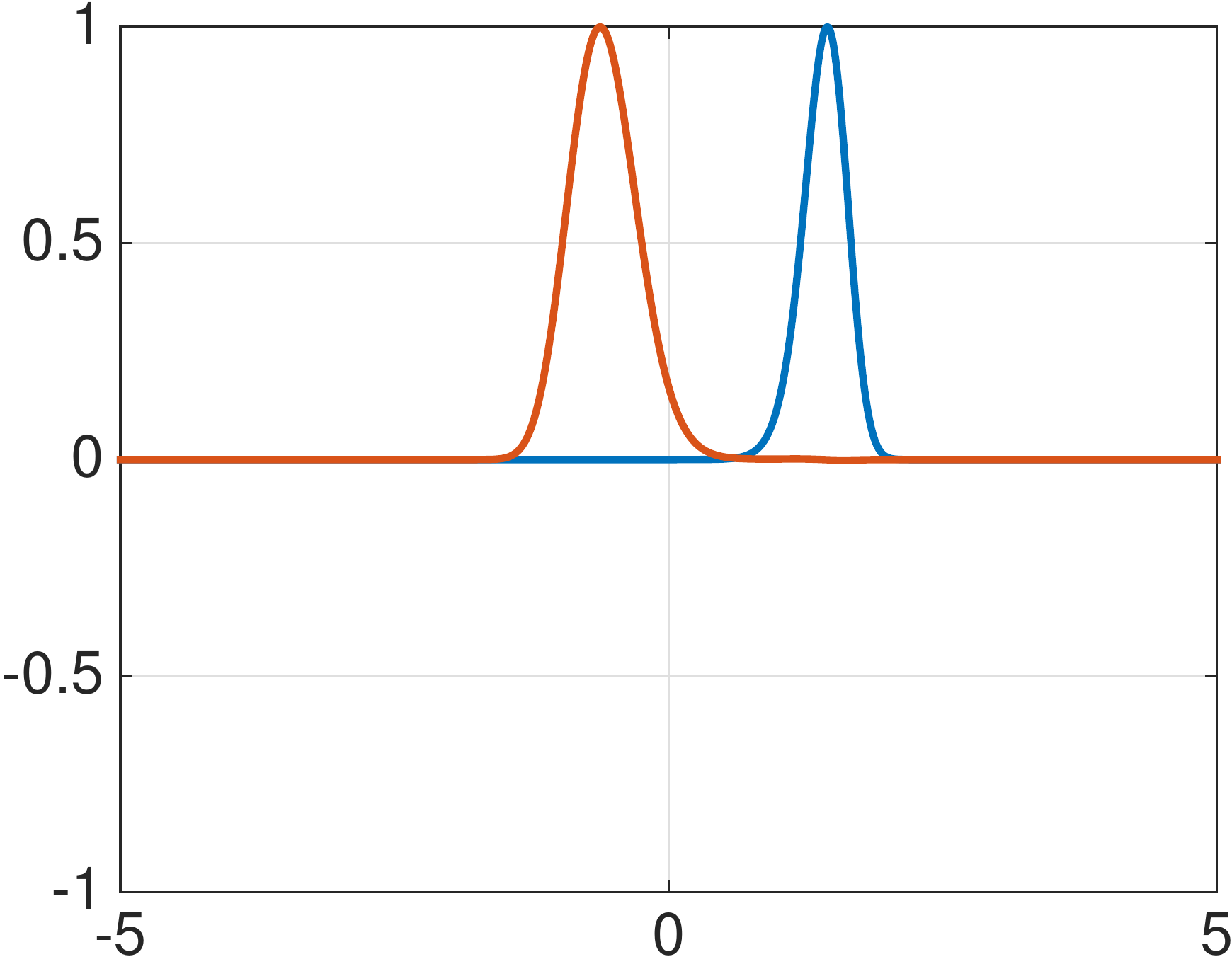}
\qquad&\includegraphics[height=3cm]{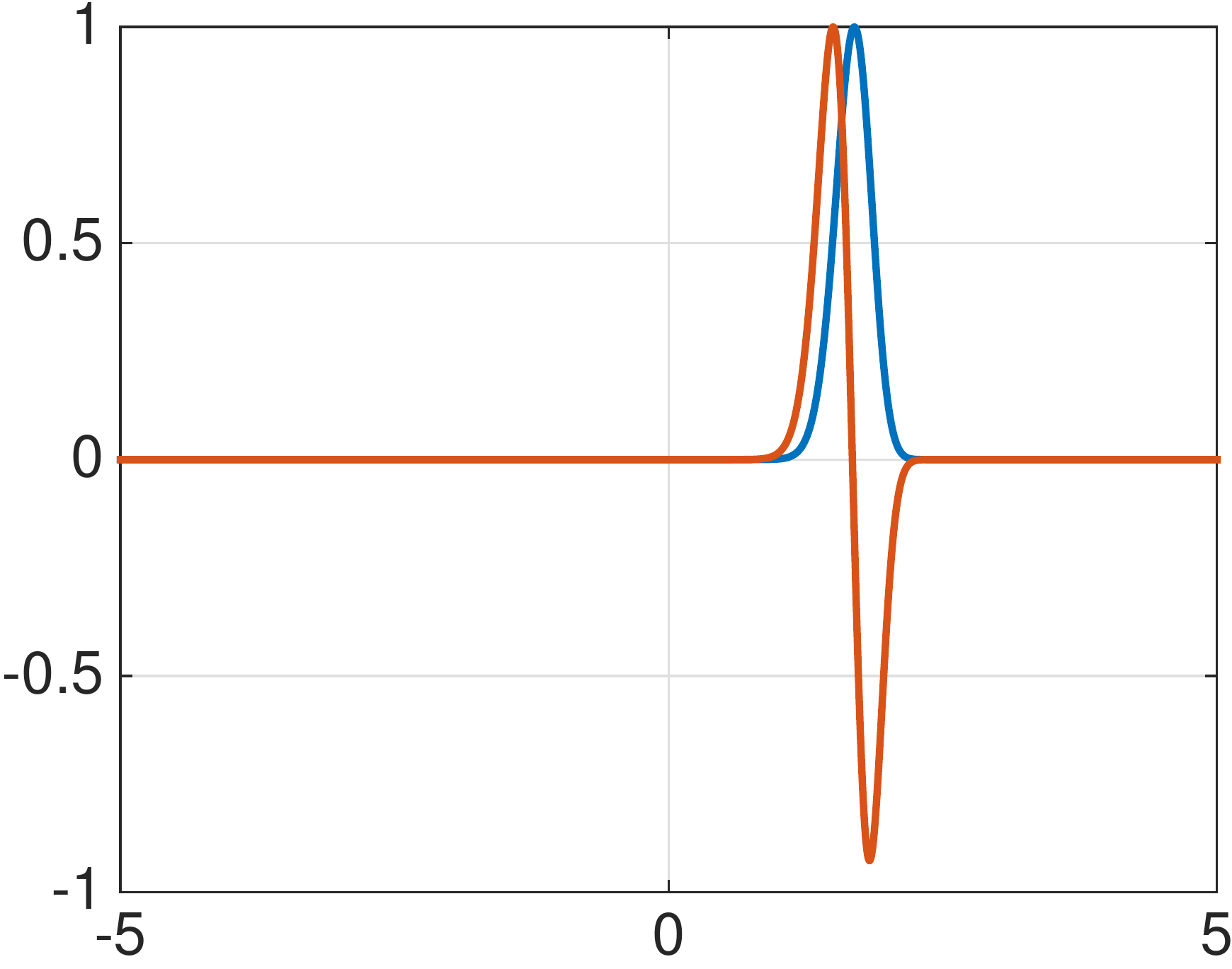}\\
\end{tabular}}\caption{First two eigenfunctions on $[-R/2,R/2]$.\label{fig.CasCVecP}}
\end{center}
\end{figure}
In Figure~\ref{fig.CasCVecP} we present the first two eigenfunctions  for $\gamma=1, 10$ and $\delta=1,5,10$. We observe that
 the eigenfunctions are no longer symmetric. 
In Figure~\ref{fig.CasCVP}
 we plot the first positive eigenvalue $\lambda_{1}(\gamma,\delta)$ as a function of $\delta$ for several values of $\gamma$: $\gamma=1, 5, 10, 25$. 
As in the previous cases, Figure~\ref{fig.CasCDK} gives the approximation of the diffusion and drift coefficients using our algorithm and the comparison with \eqref{eq.drift1D}. We observe, even for this nonsymmetric case, the good approximation of the drift coefficient.  
Figure~\ref{fig.CasCDKn} illustrates the convergence of the algorithm for $\gamma=1$ and $\delta=1,5, 10$. As expected, we need more eigenmodes than for the smooth symmetric potentials in order to obtain an accurate approximation of the drift and diffusion coefficients. 

\section{Conclusions}
\label{SectionC}

A new method for calculating the drift and diffusion coefficients in the diffusion approximation for a swarming model was presented in this paper. Our method is based on the calculation of the eigenvalues and eigenfunctions of an appropriate Schr\"{o}dinger operator. This operator  was obtained after a unitary transformation of the Markov generator that appears in the  Poisson equation arising in the definition of the coefficients of the limiting problem.
 The eigenvalue problem for this Schr\"{o}dinger operator was solved using a high order Finite Elements approximation. Our numerical method was tested to a few simple potentials and the effects of a tilt and of lack of smoothness of the potential on the drift and diffusion coefficients were investigated.
We also investigated the difficulties related to the tunnelling effect that appears in the "semiclassical" limit $\gamma \gg 1$.

We believe that the numerical method developed in this paper can be applied to the calculation of effective coefficients in a wide variety of diffusion approximations, coarse-grained and mean field models that appear in kinetic theory or in homogenization theory. The crucial observation is that in many different settings effective coefficients are given in terms of the solution to an appropriate {\it linear} Poisson equation~\cite{BLP,GP,OllaKomorowskiLandim12,PavlSt08,ResibDeLeen77}. Thus, we believe that the spectral approach advocated in this paper can be of more general interest.
Moreover, this approach also opens to relevant perspectives:

\begin{enumerate}
\item It would be interesting to develop a  detailed analysis of rates of convergence, and a careful study of the computational cost, 
depending on the numerical parameters of the approximation ($R$, $\mu$, etc).
Depending on the underlying operators, the method might also benefit from the use of appropriate preconditioners.

\item It would be interesting to compare the performance of the proposed numerical method with alternative techniques such as Monte Carlo methods~\cite{PavlStBan06} expansions in orthogonal polynomials~\cite{PavlVog08} etc.

\item The method can be  extended to more general types of Poisson equations (and the calculation of the corresponding effective coefficients) including hypoelliptic operators of   
Schr\"{o}dinger type that appear in, e.g.~\cite{HairPavl0,HP07}, or auxilliary equations, possibly depending on the space variable, that appear in homogenization theory \cite{BLP,GP}.

\end{enumerate}  

%%%%%%%%%%%%%%%%%%%%%%%%%%%%%%%%%%%%%%%%%%

\appendix
\section{Analysis of the diffusion asymptotics\label{app.A}}

\subsection{Set up of the problem}

In this Section we provide a few hints about the analysis of the asymptotic regime $\varepsilon\rightarrow 0$ of the problem
\begin{equation*}%\label{a1}
\partial_t f^\varepsilon+\ds\frac1{\sqrt\varepsilon}\big(v\cdot\nabla_rf ^\varepsilon  -\nabla_r \Phi^\varepsilon \cdot\nabla_v f^\varepsilon)=\ds\frac{1}{\varepsilon}Q f^\varepsilon\,,
\end{equation*}
where on the one hand 
\[ 
\Phi^\varepsilon(t,r)=U\star\rho^\varepsilon(t,r),\qquad \rho^\varepsilon(t,r)=\ds\int f^\varepsilon(t,r,v)\ud v,
\]
and on the other hand
\[
Q f=\nabla_v \cdot(\nabla_v Wf+\nabla_v f).
\]
For the sake of simplicity, here and below, the scaling parameter $\theta$ is set to 1.
The system is completed with the initial data 
\[
f^\varepsilon\big|_{t=0}=f^\varepsilon_{\mathrm{Init}}.\]

We set up the technical assumptions on the potentials $W$ and $U$ as follows.
As we mentioned earlier, a typical example motivated from the modelling of swarming leads
to the standard quartic potential with $\nabla_vW(v)=v(\frac1\gamma |v|^2 -1)$.
Hence, in what follows, we assume that $W$ is smooth, radially symmetric and coercive:
\[\begin{array}{l}
\textrm{$v\mapsto W(v) $ 
 depends only on $|v|$,  is of class $C^2$, and there exist 
$C, R, \alpha>0$ such that }
\\
\textrm{
$W(v)\geq C(1+|v|^{2+\alpha})$ for any $|v|\geq R$,}
\\
\textrm{$|\nabla_vW(v)|\leq C(1+W(v))$   for any $v\in\mathbb R^d$}.\end{array}\]
Consequently, $W$ is bounded from below; possibly at the cost of adding to $W$ a constant, we assume that 
\[W(v)\geq 0, \quad \textrm{ for any $v\in\mathbb R^d$}.\]
For $U$ we  
suppose 
\[
\nabla_r U\in L^\infty(\mathbb R^d).
\]
It is likely that this assumption is far from optimal, but the technicalities to relax it are beyond the scope of this paper.
It is convenient to set
\[
M(v)=Z^{-1} e^{-W(v)}\in L^1(\RR^d),\qquad Z= \ds\int e^{-W(v)}\ud v .\]
We further assume that
\[
W(v)M(v)\in L^1(\RR^d).
\]
Note that, due to the fact that $W$ is radially symmetric, 
\begin{equation}\label{nullflux}\ds\int vM(v)\ud v=0\end{equation}
 holds.
This is crucial to the analysis.
As mentioned above, the discussion is strongly inspired by the study of the Vlasov-Poisson-Fokker-Planck system \cite{EGM, MT}. In this Appendix we present a proof of Theorem~\ref{CVTh}. Naturally, the proof of this theorem consists in two parts: 
the derivation of a priori estimates presented in Section~\ref{sec:apriori} and the passage to the limit, presented in Section~\ref{sec:limit}.

\subsection{A priori estimates}
\label{sec:apriori}

We start by observing that
\[\ds\frac{\ud }{\ud t}\ds\iint f^\varepsilon\ud v\ud r=0\]
holds, which means that the total mass is conserved.
Here and below, we always assume 
\[\ds\sup_{\varepsilon>0}\|f_{\mathrm{Init}}^\varepsilon\|_{L^1}<\infty.\]
Owing to the 
regularity of $U$, we observe that
\[\nabla_r \Phi^\varepsilon(t,r)=\ds\int \nabla_r U(r-r')\left(\ds\int f^\varepsilon(t,r',v)\ud v\right)\ud r'\]
verifies
\begin{equation}\label{est_pot_Linf}
|\nabla_r \Phi^\varepsilon(t,r)|\leq \|\nabla_r U\|_\infty \|f^\varepsilon_{\mathrm{Init}}\|_{L^1}
.
\end{equation}
Next, we compute
\[\ds\frac{\ud }{\ud t}\ds\iint \big(W+\ln (f^\varepsilon)\big)f^\varepsilon\ud v\ud r=
-\ds\frac{1}{\varepsilon}\ds\iint 
|\nabla_v W \sqrt {f^\varepsilon}+2\nabla_v\sqrt {f^\varepsilon}|^2
\ud v\ud r- \ds\frac{1}{\sqrt\varepsilon}
\ds\iint \nabla_r \Phi ^\varepsilon f^\varepsilon \cdot\nabla_v W\ud v\ud r.
\]
The last term can be rewritten as 
\[
-\ds\frac{1}{\sqrt\varepsilon}
\ds\iint \nabla_r \Phi^\varepsilon f^\varepsilon \cdot\nabla_v W\ud v\ud r
=
-
\ds\iint \nabla_r \Phi^\varepsilon \sqrt {f^\varepsilon} \cdot \ds\frac{\nabla_v W\sqrt {f^\varepsilon} + 2\nabla_v \sqrt {f^\varepsilon}}{\sqrt\varepsilon}\ud v\ud r
\]
which is dominated by 
\begin{align*}
\ds\frac12\|\nabla_r \Phi^\varepsilon\|_\infty^2 
&\ds\iint f^\varepsilon\ud v\ud r + \ds\frac{1}{2\varepsilon}\ds\iint 
|\nabla_v W \sqrt {f^\varepsilon} +2\nabla_v\sqrt {f^\varepsilon}|^2
\ud v\ud r
\\
& \leq 
\ds\frac12\|\nabla_r U\|_\infty^2 
\|f^\varepsilon_{\mathrm{Init}}\|_{L^1}^3
+ \ds\frac{1}{2\varepsilon}\ds\iint 
|\nabla_v W \sqrt {f^\varepsilon} +2\nabla_v\sqrt {f^\varepsilon}|^2
\ud v\ud r,
\end{align*}
bearing in  mind \eqref{est_pot_Linf}.

In order to control the behavior at infinity of the particle distribution function, we also 
need to evaluate
\begin{align*}
\ds\frac{\ud }{\ud t}\ds\iint |r| f^\varepsilon\ud v\ud r&=
\ds\frac{1}{\sqrt\varepsilon}
\ds\iint v f^\varepsilon\cdot \ds\frac{r}{|r|}\ud v\ud r
\\
&=
\ds\iint v \ds\frac{f^\varepsilon -\rho^\varepsilon M}{\sqrt\varepsilon}
\cdot \ds\frac{r}{|r|}\ud v\ud r
\\
&=
\ds\iint  \ds\frac{\sqrt{f^\varepsilon} -\sqrt{\rho^\varepsilon M}}{\sqrt\varepsilon}
v(\sqrt{f^\varepsilon} +\sqrt{\rho^\varepsilon M})\cdot \ds\frac{r}{|r|}\ud v\ud r
\\
&\leq
\left(
\ds\iint 
v^2\big|\sqrt{f^\varepsilon} +\sqrt{\rho^\varepsilon M}\big|^2\ud v\ud r
\right)^{1/2}
\left(
\ds\iint \ds\frac{|\sqrt{f^\varepsilon} -\sqrt{\rho^\varepsilon M}|^2}{\varepsilon}
\ud v\ud r
\right)^{1/2}
\\
&\leq
\left(
2\ds\iint 
v^2\big(f^\varepsilon +\rho^\varepsilon M\big)\ud v\ud r
\right)^{1/2}
\left(
\ds\int \ds\frac{|\sqrt{f^\varepsilon} -\sqrt{\rho^\varepsilon M}|^2}{\varepsilon}
\ud v\ud r
\right)^{1/2}.
\end{align*}
We note that the dissipation term recasts as
\[
\ds\int 
|\nabla_v W \sqrt {f^\varepsilon} +2\nabla_v\sqrt {f^\varepsilon}|^2
\ud v\ud r
=
\ds\frac{4}{\varepsilon}\ds\int 
\Big|\nabla_v\sqrt{\ds\frac {f^\varepsilon}{M}}\Big|^2
M\ud v\ud r\ud t.\]
By the logarithmic Sobolev inequality, see \cite[Th.~8.14]{LL},
there exists $\lambda>0$ such that
\[
 \ds\int 
\Big(\ds\frac{f^\varepsilon}{\rho^\varepsilon M}\ln\Big(\ds\frac{f^\varepsilon}{\rho^\varepsilon M}\Big)- \ds\frac{f^\varepsilon}{\rho^\varepsilon M} +1\Big)
\rho^\varepsilon M\ud v
\leq 
\lambda\ds\int 
\Big|\nabla_v\sqrt{\ds\frac {f^\varepsilon}{M}}\Big|^2
M\ud v.
\]
Additionally, we shall make use of the elementary inequality
\[|\sqrt s-1|^2\leq s\ln (s)- s+1\]
which holds for any $s\geq 0$.
It follows that
\[
\ds\int \Big|\sqrt{\ds\frac{f^\varepsilon}{\rho^\varepsilon M}}- 1\Big|^2 
\rho^\varepsilon M
\ud v
=
\ds\int \big|\sqrt{f^\varepsilon}-\sqrt{\rho^\varepsilon M}\Big|^2 
\ud v
\leq 
\ds\frac\lambda 4
\ds\int 
|\nabla_v W \sqrt {f^\varepsilon} +2\nabla_v\sqrt {f^\varepsilon}|^2
\ud v.\]
Hence, we obtain 
\[
\ds\frac{\ud }{\ud t}\ds\iint |r| f^\varepsilon\ud v\ud r
\leq 
\lambda \Big(\ds\iint v^2 f^\varepsilon \ud v \ud r+ \ds\iint \rho^\varepsilon M\ud v\ud r\Big)+\ds\frac{1}{4\varepsilon}
\ds\iint 
|\nabla_v W \sqrt {f^\varepsilon} +2\nabla_v\sqrt {f^\varepsilon}|^2
\ud v\ud r.
\]
However, the assumption on the potential $W$ guarantees
that we can find $0< \Xi%K
<\infty$ such that 
\[
\ds\iint v^2 M\ud v\leq 
\Xi %K
,\qquad
\ds\iint v^2 f^\varepsilon\ud v\ud r\leq \Xi %K
\ds\int (1+W) f^\varepsilon\ud v\ud r.
\]
In what follows, we will still denote by $\Xi %K
$ a positive constant which only depends on the data, but which is uniform with respect to $\varepsilon$ and $t$. 
Finally, we arrive at 
\[
\ds\frac{\ud}{\ud t}
\ds\iint f^\varepsilon\big(1+W+|r|+\ln(f^\varepsilon)\big)\ud v\ud r
+\ds\frac{1}{4\varepsilon}
\ds\iint 
|\nabla_v W \sqrt {f^\varepsilon} +2\nabla_v\sqrt {f^\varepsilon}|^2
\ud v\ud r
\leq \Xi %K
\Big(1+\ds\iint Wf^\varepsilon \ud v\ud r\Big). 
\]
We now use the classical trick
\[
s|\ln (s)|=s\ln (s) -2s\ln (s)\mathbf 1_{0\leq s\leq e^{-\Omega}}-2s\ln (s)\mathbf 1_{e^{-\Omega}< s\leq 1}
\leq
s\ln (s)+\Xi %K
e^{-\Omega/2} + 2\Omega s,
\]
with $s=f^\varepsilon(t,x,v)$, $\Omega=\frac14(|r|+W(v))$.
We deduce the following estimate
\[
\begin{array}{l}
\ds\iint f^\varepsilon\Big(1+\ds\frac{W+|r|}{2}+|\ln(f^\varepsilon)|\Big)\ud v\ud r
+\ds\frac{1}{4\varepsilon}
\ds\int_0^t\ds\iint 
|\nabla_v W \sqrt {f^\varepsilon} +2\nabla_v\sqrt {f^\varepsilon}|^2
\ud v\ud r\ud s
\\
\leq 
\ds\iint f^\varepsilon\Big(1+W+|r|+\ln(f^\varepsilon)\big)\ud v\ud r
\\
\qquad
\qquad\qquad
+\ds\frac{1}{4\varepsilon}
\ds\int_0^t\ds\iint 
|\nabla_v W \sqrt {f^\varepsilon} +2\nabla_v\sqrt {f^\varepsilon}|^2
\ud v\ud r\ud s+\Xi %K
\ds\iint \exp\Big(-\ds\frac{|r|+W(v)}{4}\Big)\ud v\ud r
\\
\leq 
\ds\iint f^\varepsilon_{\mathrm{Init}}\Big(1+W+|r|+\ln(f^\varepsilon_{\mathrm{Init}})\big)\ud v\ud r
\\
\qquad\qquad +\Xi %K
\Big(t+\ds\iint \exp\Big(-\ds\frac{|r|+W(v)}{4}\Big)\ud v\ud r\Big)
+\Xi %K
\ds\int_0^t \ds\iint Wf^\varepsilon\ud v\ud r\ud s.
\end{array}
\]
It remains to appeal to the Gr\"onwall lemma to conclude with the following statement.

\begin{prop}
We assume that $f^\varepsilon_{\mathrm{Init}}\geq 0$ satisfies
\[
\ds\sup_{\varepsilon>0}\ds\iint \big(1+|r|+ W+|\ln (f^\varepsilon_{\mathrm{Init}}|)\big)f^\varepsilon_{\mathrm{Init}}\ud v\ud r= M_0<0.
\]
Let $0<T<
\infty$. 
Then, there exists $0<C_T<\infty$, depending only on $W$, $M_0$ and $T$ such that 
\begin{align*}
\ds\sup_{\varepsilon>0}\ds\sup_{0\leq t\leq T}\ds\iint \big(1+|r|+ W+|\ln (f^\varepsilon)|\big)f^\varepsilon\ud v\ud r\leq C_T,
\\
\ds\sup_{\varepsilon>0} \ds\int_0^T\ds\iint \ds\frac{|2\nabla_v\sqrt{ f^\varepsilon}+ \nabla_v W \sqrt f^\varepsilon|^2}{\varepsilon}\ud v\ud x\ud t\leq C_T.
\end{align*}
\end{prop}

This estimate 
can be translated by means of compactness properties, 
as an application of the Dunford-Pettis theorem, see e. g. \cite[Section 7.3.2 \& 7.3.4]{Th}.

\begin{coro}\label{comp_f}
The sequence $\big(f^\varepsilon\big)_{\varepsilon>0}$ is weakly compact in $L^1((0,T)\times \mathbb R^d\times\mathbb R^d)$.
\end{coro}

\begin{coro}\label{compg}
We can write $f^\varepsilon(t,r,v)=\rho^\varepsilon(t,r) M(v)+\sqrt\varepsilon g^\varepsilon(t,r,v)$, with $\int g^\varepsilon\ud v=0$, where
the sequence $\big(\rho^\varepsilon\big)_{\varepsilon>0}$  is weakly compact in $L^1((0,T)\times \mathbb R^d)$ and the sequence $\big(g^\varepsilon\big)_{\varepsilon>0}$ is weakly compact in $L^1((0,T)\times \mathbb R^d\times\mathbb R^d)$.
\end{coro}

\noindent
{\bf Proof.}
The compactness of the macroscopic density $\rho^\varepsilon$ is a direct consequence of Corollary \ref{comp_f}.
Next, let $A\subset (0,T)\times \mathbb R^d\times\mathbb R^d$ be a measurable subset.
Reproducing a manipulation already detailed above we get
\begin{align*}
\ds\iiint _A |g^\varepsilon|\ud v\ud r\ud t
&=
\ds\iiint _A \ds\frac{\sqrt{f^\varepsilon} -\sqrt{\rho^\varepsilon M}}{\sqrt\varepsilon}
(\sqrt{f^\varepsilon} +\sqrt{\rho^\varepsilon M})\ud v\ud r\ud t
\\
&\leq
\left(
2\ds\iiint _A (f^\varepsilon+\rho^\varepsilon M)\ud v\ud r\ud t 
\right)^{1/2}
\left(
\ds\iiint \ds\frac{|\sqrt{f^\varepsilon} -\sqrt{\rho^\varepsilon M}|^2}{\varepsilon}
\ud v\ud r\ud t
\right)^{1/2}
\\
&\leq  
\left(
2\ds\iiint _A (f^\varepsilon+\rho^\varepsilon M)\ud v\ud r\ud t 
\right)^{1/2}
\left(
\ds\int_0^T\ds\iint \ds\frac{|2\nabla_v\sqrt{ f^\varepsilon}+ \nabla_v W \sqrt f^\varepsilon|^2}{\varepsilon}
\ud v\ud r\ud t
\right)^{1/2}
\\
&\leq  \sqrt{2C_T}\ \left(
\ds\iiint _A (f^\varepsilon+\rho^\varepsilon M)\ud v\ud r\ud t 
\right)^{1/2}.
\end{align*}
Coming back to Corollary  \ref{comp_f}, it proves the equi-integrability of the sequence $g^\varepsilon$.
\endproof

Owing to the behavior of $W$ for large $|v|$'s, we shall use the fact that
 the reasoning can be applied to $f^\varepsilon \psi (v)$ 
 with tests functions $v\mapsto \psi(v)$ verifying $\lim_{|v|\rightarrow\infty}\frac{|\psi(v)|}{W(v)}=0$
 and to $g^\varepsilon \psi(v)$ as well, with tests functions $v\mapsto \psi(v)$ verifying $\lim_{|v|\rightarrow\infty}\frac{|\psi(v)|}{\sqrt{W(v)}}=0$
 %,\qquad \ds\lim_{|v|\rightarrow
%\infty} \ds\frac{\psi(v)}{W(v)}=0.\]
In particular, using  \eqref{nullflux}, we have 

\begin{coro}\label{compJ}
The sequence defined by  $$J^\varepsilon(t,r)=\ds\int \ds\frac {v}{\sqrt\varepsilon} \ f^\varepsilon (t,r,v)\ud v=
\ds\int v\ g^\varepsilon (t,r,v)\ud v$$  is weakly compact in $L^1((0,T)\times \mathbb R^d)$.
 \end{coro}

\subsection{Passage to the limit}
\label{sec:limit}

According 
to the discussion above, we can assume, possibly at the price of extracting a subsequence, that 
\[f^\varepsilon\rightharpoonup f\quad\textrm{weakly in $L^1((0,T)\times \mathbb R^d\times \mathbb R^d)$}.\]
It means that 
\[\ds\lim_{\varepsilon\rightarrow 0}\ds\int_0^T\ds\iint f^\varepsilon (t,r,v)\psi(t,r,v)\ud v\ud r\ud t=
\ds\int_0^T\ds\iint f (t,r,v)\psi(t,r,v)\ud v\ud r\ud t
\]
holds for any trial function $\psi\in L^\infty
((0,T)\times \mathbb R^d\times \mathbb R^d)$.
In fact, owing to the uniform estimate on the integral of $Wf^\varepsilon$
the convergence applies for $\psi(t,r,v)=\zeta(t,r)\phi(v)$, with $\zeta\in L^\infty
((0,T)\times \mathbb R^d)$ and $\phi\in C^0(\mathbb R^d)$ such that 
$\lim_{|v|\rightarrow\infty}\frac{\phi(v)}{W(v)}=0$.
Accordingly, we have
\[
\rho^\varepsilon\rightharpoonup \rho=\ds\int f\ud v\quad\textrm{weakly in $L^1((0,T)\times \mathbb R^d)$}.\]
Integrating the equation with respect to $v$, we obtain 
\begin{equation}\label{conteqeps}
\partial_t \rho^\varepsilon+\nabla_r\cdot J^\varepsilon=0.\end{equation}
With Corollary \ref{compg} and Corollary \ref{compJ}, we can also suppose that 
\[\begin{array}{l}
g^\varepsilon\rightharpoonup g\quad\textrm{weakly in $L^1((0,T)\times \mathbb R^d\times \mathbb R^d)$},
\\
J^\varepsilon\rightharpoonup J\quad\textrm{weakly in $L^1((0,T)\times \mathbb R^d)$}.\end{array}\]
Letting $\varepsilon $ go to 0 in \eqref{conteqeps} we are thus led to 
\[\partial_t 
\rho+\nabla_r \cdot J=0.\]
Furthermore, for any trial function $\zeta \in W^{1,\infty}(\mathbb R^d)$, \eqref{conteqeps} and Corollary \ref{compJ}
imply that
\[
\ds\int \rho^\varepsilon(t+h,r)\zeta(r)\ud r- \ds\int \rho^\varepsilon(t,r)\zeta(r)\ud r =
\ds\int_t^{t+h}\ds\int J^\varepsilon(s,r)\cdot\nabla_r \zeta(r)\ud r \ud s\]
can be made arbitrarily small as $h$ goes to 0, uniformly with respect to $\varepsilon$.
Up to a suitable regularization argument, we deduce that $\rho^\varepsilon$ is compact in $C^0([0,T];L^1(\mathbb R^d)-\mathrm{weak})$.
As a consequence, extracting further subsequences if necessary, we 
deduce that
\[
\nabla_r\Phi^\varepsilon(t,r)
 =
 \ds\int \nabla_r U(r-r')\rho^\varepsilon(t,r')\ud r'\xrightarrow[\varepsilon\rightarrow 0]{}
 \nabla_r\Phi(t,r)
 =
 \ds\int \nabla_r U(r-r')\rho(t,r')\ud r'  
\]
holds  a.e. $(0,T)\times \mathbb R^d$, with the uniform 
estimate \eqref{est_pot_Linf}.

It remains to identify 
the limiting current $J$.
Going back to the proof of Corollary \ref{compg}, we justify that 
\[\ds\lim_{\varepsilon\rightarrow 0}\ds\int_0^T\ds\iint g^\varepsilon (t,r,v)\psi(t,r,v)\ud v\ud r\ud t=
\ds\int_0^T\ds\iint g (t,r,v)\psi(t,r,v)\ud v\ud r\ud t
\]
holds for any trial function $\psi\in L^\infty
((0,T)\times \mathbb R^d\times \mathbb R^d)$ as well as  for $\psi(t,r,v)=\zeta(t,r)\phi(v)$, with $\zeta\in L^\infty
((0,T)\times \mathbb R^d)$ and $\phi\in C^0(\mathbb R^d)$ such that 
$\lim_{|v|\rightarrow\infty}\frac{\phi(v)}{\sqrt {W(v)}}=0$.
In particular, we have
\[
J^\varepsilon(t,r)=\ds\int vg^\varepsilon(t,r,v)\ud v\rightharpoonup J(t,r)=\ds\int vg(t,r,v)\ud v\quad
\textrm{weakly in $L^1((0,T)\times\mathbb R^d)$}.
\]

Then, for any $\zeta\in C^\infty_c((0,T)\times\mathbb R^d)$ and $\phi\in C^\infty_c(\mathbb R^d)$, we have
\[\begin{array}{l}
- \varepsilon  \ds\int_0^T\ds\iint
 f^\varepsilon(t,r,v)\ \phi (v)\partial_t \zeta(t,r)
  \ud v\ud r\ud t
 - 
  \ds\int_0^T\ds\iint
 f^\varepsilon(t,r,v)\phi  (v)v\cdot\nabla_r  \zeta(t,r)
  \ud v\ud r\ud t
 \\\qquad+
  \ds\int_0^T\ds\iint
 f^\varepsilon(t,r,v)\ \nabla_v \phi (v) \cdot\nabla_r \Phi^\varepsilon(t,r) \  \zeta(t,r)
  \ud v\ud r\ud t
  \\
  \qquad 
  \qquad
  =
 \ds\frac{1}{\sqrt\varepsilon} \ds\int_0^T
\mathcal{L}\phi (v)\, f^\varepsilon(t,r,v)\zeta(t,r)
 \ud v\ud r\ud t
 =
\ds\int_0^T
\mathcal{L} \phi (v) \,g^\varepsilon(t,r,v)\zeta(t,r)
\ud v\ud r\ud t\,,
\end{array}
\]
where $\mathcal{L}$ is the $L^2$-adjoint operator to $Q$ defined in \eqref{eq.defL}.
Letting $\varepsilon$ go to 0, we are led to 
\begin{align*}
\ds\int_0^T \ds\iint 
\mathcal{L}\phi (v) \, g(t,r,v)\zeta(t,r)
\ud v\ud r\ud t
&=
- \ds\int_0^T\ds\iint
 f(t,r,v)\ \phi (v)v\cdot\nabla_r  \zeta(t,r)
  \ud v\ud r\ud t
 \\ &+
  \ds\int_0^T\ds\iint
 f(t,r,v)\ \nabla_v \phi (v) \cdot\nabla_r \Phi(t,r) \  \zeta(t,r)
  \ud v\ud r\ud t.
\end{align*}
For the nonlinear term, we have combined the weak convergence 
of $f^\varepsilon$ in $L^1$ and the pointwise convergence of the uniformly bounded sequence $\nabla_r \Phi^\varepsilon$, see \cite[Lemma 7.62]{Th}.
By virtue of \eqref{compg}, we have
\[
f(t,r,v)=\rho(t,r) M(v).\]
Therefore, $g$ is characterized by the relation
\begin{multline}\label{eq_g}
\ds\int_0^T\ds\iint
\mathcal{L}\phi (v) \,g(t,r,v)\zeta(t,r)
\ud v\ud r\ud t
=
- \ds\int_0^T\ds\int
 \left(\ds\int \phi (v) v M(v)\ud v\right)\cdot
 \rho(t,r)
 \nabla_r  \zeta(t,r)
 \ud r\ud t
\\
\qquad \,\,\,\,\,\,+
  \ds\int_0^T\ds\int
\left(\ds\int \nabla_v \phi (v)  \ M(v) \ud v\right)\cdot
 \nabla_r \Phi(t,r) \rho(t,r)  \zeta(t,r)
\ud r\ud t,
\end{multline}
for any $\zeta\in C^\infty_c((0,T)\times\mathbb R^d)$, $\phi\in C^\infty_c(\mathbb R^d)$.
Note in particular that
\[
\ds\int \nabla_v \phi (v)  \ M(v) \ud v=
\ds\int  \phi (v)  \nabla_v W(v) M(v) \ud v
.\]
Going back to \eqref{eq_g}, it follows that
\[\begin{array}{l}
\left|
\ds\int_0^T\ds\iint
\mathcal{L}\phi (v) \,g(t,r,v)\zeta(t,r)
\ud v\ud r\ud t
\right|
\leq
\|\phi\|_{L^2(M\ud v)}\ 
\|\rho\|_{L^\infty(0,T;L^1(\mathbb R^d))} T
\\
\qquad
\times
\Big(\|\nabla_r\zeta\|_{L^{\infty}}
\|v\sqrt M\|_{L^2}
+
\|\nabla_r U\|_{L^\infty}\|\rho\|_{L^\infty(0,T;L^1(\mathbb R^d))} 
\|\zeta\|_{L^\infty}
\|\nabla_v W \sqrt M\|_{L^2}
\Big)
\end{array}\]
which makes the regularity of $Qg$ precise.
\\

We recall that $Q$ is a differential operator associated to a natural quadratic form with  domain 
\[D(Q)=\Big\{f\in L^2\Big(\frac{\ud v}{M(v)}\Big),\ \nabla_v\Big(\frac fM\Big)\in L^2(M(v)\ud v)\Big\}\]
and we have 
\[ -\ds\int Q f \ g\ud v = \ds\int \nabla_v\ds\frac fM\cdot \nabla_v\ds\frac gM\ M\ud v. \]
Note that $L^2(\frac{\ud v}{M(v)})\subset L^1$, and  by the Sobolev inequality \cite[Corollary 2.18]{AMTU}, 
we check that $Q$ is coercive on the closed subset of functions with zero-mean.
 Therefore, the Fredholm alternative applies
 as follows.
 
\begin{lem}
For any $h\in L^2(\frac{\ud v}{M(v)})$  such that $\int h\ud v=0$ (resp. $h \in L^2(M\ud v)$ such that $\int hM \ud v=0$), there exists  a unique 
$f\in D(Q)$ such that $Qf=h$ and $\int f\ud v=0$ (resp. $f\in D(\mathcal{L})$ such that $\mathcal{L} f =h$ and $\int f M \ud v=0$).
\end{lem}

Consequently, we can define $\hat \chi$ solution of $\mathcal{L} \hat \chi=v \in L^2(M\ud v)$ as in \eqref{e:invert-L}. Therefore, we can rewrite the limiting current  as follows
\[
J^\varepsilon(t,r)=
\ds\int vg^\varepsilon(t,r,v)\ud v\rightharpoonup J(t,r)=\ds\int vg(t,r,v) \ud v=
\ds\int \mathcal{L}\hat\chi (v)\, g(t,r,v) \ud v.
\]
Identifying limits with \eqref{eq_g}, up to a regularization argument, we arrive at 
\begin{align*}
\ds\int_0^T\ds\int J(t,r)\zeta(t,r)\ud r\ud t
&=
-  \ds\int_0^T\ds\int
 \left(\ds\int \hat\chi (v)\otimes v M(v)\ud v\right)
 \rho(t,r)
 \nabla_r  \zeta(t,r)
 \ud r\ud t
 \\&+ 
  \ds\int_0^T\ds\int
\left(\ds\int  \hat\chi (v)  \otimes \nabla_v W(v) M(v) \ud v\right)
 \nabla_r \Phi(t,r) \rho(t,r)  \zeta(t,r)
\ud r\ud t.
\end{align*}
It corresponds to the dual formulation of the expected relation. 
Indeed, we observe that, on the one hand
\[
\ds\int \hat\chi (v)\otimes v M(v)\ud v=
\ds\int \hat\chi (v)\otimes Q\chi(v)\ud v
=
\ds\int v\otimes \chi(v)\ud v
=-\mathcal D,\]
and, on the other hand
\[
\ds\int  \hat\chi (v)  \otimes \nabla_v W(v) M(v) \ud v=
\ds\int  \hat\chi (v)  \otimes Q\kappa (v) \ud v
=
\ds\int  v \otimes  \kappa (v) \ud v=-\mathcal K.
\]
Hence we have obtained
\[
J=-(\mathcal D\nabla_r\rho+\mathcal K\nabla_r\Phi\ \rho)
.\]
It ends the proof of Theorem \ref{CVTh}.
\\

While this is not necessary for establishing the connection between the kinetic model and the drift-diffusion equation in the regime $\varepsilon\to 0$, it is possible to improve the compactness of the macroscopic density from weak to strong.
The proof relies on the combination of a  renormalization
argument and velocity averaging techniques \cite{GLPS}.
This is detailed in \cite{EGM, MT}.
We start with the following averaging lemma.

\begin{lem}\cite[Prop. 4.1]{MT}\label{av_lem}
Let $h^\varepsilon$ be a uniformly bounded sequence in $L^2((0,T)\times \mathbb R^d\times\mathbb R^d)$.
We assume that 
\[(\sqrt\varepsilon\partial_t +v\cdot\nabla_r )h^\varepsilon=h_0^\varepsilon +\nabla_v \cdot h_1^\varepsilon,\]
with $h^\varepsilon_0$ and $h^\varepsilon_1$ bounded in $L^1 ((0,T)\times \mathbb R^d\times\mathbb R^d)$.
Then, for any function $\psi\in C^\infty_c(\mathbb R^d)$, we have
\[ 
\ds\sup_{\varepsilon>0}
\left\| \ds\int \big( h^\varepsilon(t,r+h,v)- h^\varepsilon(t,r,v)\big)\ \psi(v)\ud v\right\|_{L^1}
\xrightarrow[|h|\rightarrow 0]{} 0.
\]
\end{lem}

Let $0<\delta<1$. We set
\[\beta_\delta (z)=\ds\frac{z}{1+\delta z}.\]
Therefore $h^\varepsilon
=\beta_\delta (f^\varepsilon)$ is 
bounded in $L^1\cap L^\infty((0,T)\times \mathbb R^d\times\mathbb R^d)$, thus in $L^2((0,T)\times \mathbb R^d\times\mathbb R^d)$, uniformly with respect to $\varepsilon$
(but the bound depends on $\delta$).
Hence Lemma \ref{av_lem} applies to $h^\varepsilon
=\beta_\delta (f^\varepsilon)$ for any $0<\delta <1$ fixed with
\[
h^\varepsilon_0=
\beta''_\delta(f^\varepsilon)f^\varepsilon
\Big(
-\ds\frac{(\nabla_v W \sqrt{f^\varepsilon}+2\nabla_v \sqrt{f^\varepsilon})^2}{\sqrt{\varepsilon}}
+
\ds\frac{\nabla_v W \sqrt{f^\varepsilon}+2\nabla_v \sqrt{f^\varepsilon}}{\sqrt{\varepsilon}}\cdot\nabla_v W\sqrt{f^\varepsilon}
\Big)
\]
and \[ 
h^\varepsilon_1=\nabla_r\Phi^\varepsilon\ \beta_\delta(f^\varepsilon)
- \beta'_\delta(f^\varepsilon)\sqrt{f^\varepsilon}\ds\frac{\nabla_v W \sqrt{f^\varepsilon}+2\nabla_v \sqrt{f^\varepsilon}}{\sqrt{\varepsilon}}
.\]
We deduce that
\[ 
\ds\sup_{\varepsilon>0}
\left\| \ds\int \big( \beta_\delta(f^\varepsilon)(t,r+h,v)- \beta_\delta(f^\varepsilon)(t,r,v)\big)\ \psi(v)\ud v\right\|_{L^1}
\xrightarrow[|h|\rightarrow 0]{} 0
\]
holds for any $\delta >0$.
This property passes to $f^\varepsilon$ owing to the equi-integrability of the sequence.
Indeed, we split
\[
 \ds\int \big(  f^\varepsilon(t,r,v)- \beta_\delta(f^\varepsilon)(t,r,v)\big)\ \psi(v)\ud v\ud r\ud t
\]
 by considering separately $z=(t, r,v)$ such that $|z|\geq M$, or $|f^\varepsilon(z)| \geq A$ and the complementary sets.
 The integral of the former can be made arbitrarily small by chosing $M, A$ large enough.
 The integral of the latter is dealt with by using the fact that $\beta_\delta(s)$ converges uniformly on compact sets to $s$.  
 We arrive at 
 \[ 
\ds\sup_{\varepsilon>0}
\left\| \ds\int \big( f^\varepsilon(t,r+h,v)- f^\varepsilon(t,r,v)\big)\ \psi(v)\ud v\right\|_{L^1}
\xrightarrow[|h|\rightarrow 0]{} 0.
\]
Note that owing to the weighted estimates on $f^\varepsilon$ we can remove the restriction of considering compactly supported trial functions.
This is finally combined to \eqref{conteqeps}: the 
a priori estimates tells us that $\partial_t \rho^\varepsilon$ is bounded in $L^1(0,T;W^{-1,1}(\mathbb R^d))$.
We conclude by using a standard approximation argument that $\rho^\varepsilon$ satisfies the Weil-Kolmogorov-Fr\'echet criterion, see \cite[Theorem 7.56]{Th}. Therefore $\rho^\varepsilon$ converges to $\rho$ strongly in $L^1((0,T)\times\mathbb R^d)$.
% the Aubin-Simon compactness argument \cite{Sim}.

\paragraph{Acknowledgments.} 
VBN is partially supported by the ANR (Agence Nationale de la Recherche), projects {\sc Nosevol} n$^{\rm o}$ ANR-11-BS01-0019 and {\sc Aramis}  n$^{\rm o}$  ANR-12-BS01-0021.
JAC acknowledges support from the Engineering and Physical Sciences Research Council of the UK grant number EP/K008404/1, the Spanish project MTM2011-27739-C04-02, and the Royal Society through a Wolfson Research Merit Award.
GP acknowledges financial support from the Engineering and Physical Sciences Research Council of the UK through Grants No. EP/J009636, EP/L024926, EP/L020564 and EP/L025159.

%\bibliography{FPeffCoefbib}
\bibliographystyle{mnachrn}
\small
%%%%%%%%%%%%%%%%%%%%%%%%%%%%%%%%%%%%%%%%%%%%%%%%%%%%%%%%%%%%%%%%%%%%%%%%%%%%%%%%%%%%
%\bibliography{../bibtex_files/mybib}
\bibliography{BCGP}
%\bibliography{mybib}

\end{document}